\documentclass{amsart}
\usepackage{graphicx}
\usepackage{amssymb}
\usepackage{epstopdf}
\usepackage{nicefrac}
\usepackage{enumerate}
\usepackage{hyperref}
\usepackage{float}
\usepackage{color}
\usepackage{pst-all}
\usepackage{tabularx}
\input xy
\xyoption {all}

\newtheorem{thm}{THEOREM}[section]

\newtheorem{cor}[thm]{COROLLARY}

\newtheorem{defn}[thm]{DEFINITION}

\newtheorem{prob}[thm]{PROBLEM}
\newtheorem{prop}[thm]{PROPOSITION}

\newtheorem{remark}[thm]{REMARK}

\newcommand{\ds}{\displaystyle}

\newcommand{\F}{{\mathcal F}}

\newcommand{\G}{\Gamma}

\newcommand{\whtau}{\widehat{\tau}}

\newcommand{\mB}{{\mathbb B}}
\newcommand{\mC}{{\mathbb C}}
\newcommand{\mD}{{\mathbb D}}

\newcommand{\mH}{{\mathbb H}}

\newcommand{\mQ}{{\mathbb Q}}
\newcommand{\mR}{{\mathbb R}}
\newcommand{\mS}{{\mathbb S}}
\newcommand{\mT}{{\mathbb T}}
\newcommand{\mV}{{\mathbb V}}
\newcommand{\mZ}{{\mathbb Z}}

\newcommand{\cA}{{\mathcal A}}

\newcommand{\cC}{{\mathcal C}}

\newcommand{\cG}{{\mathcal G}}

\newcommand{\cI}{{\mathcal I}}

\newcommand{\cK}{{\mathcal K}}

\newcommand{\cM}{{\mathcal M}}

\newcommand{\cP}{{\mathcal P}}

\newcommand{\cS}{{\mathcal S}}

\newcommand{\cV}{{\mathcal V}}

\newcommand{\wtalpha}{{\widetilde{\alpha}}}
\newcommand{\wtbeta}{{\widetilde{\beta}}}
\newcommand{\wtgamma}{{\widetilde{\gamma}}}

\newcommand{\fM}{{\mathfrak{M}}}

\newcommand{\Diff}{{\rm Diff}}
\newcommand{\Maps}{{\rm Maps}}
\newcommand{\Hom}{{\rm Hom}}

\newcommand{\GF}{\Gamma_{\F}}

\newcommand{\oG}{\overline{\Gamma}}

\begin{document}

\title{Characteristic classes for flat Diff(M)-foliations} 

\author{Steven Hurder}
\address{Steven Hurder, Department of Mathematics, University of Illinois at Chicago, 322 SEO (m/c 249), 851 S. Morgan Street, Chicago, IL 60607-7045}
\email{hurder@uic.edu}
\thanks{Version date:  November 15, 2023}

\date{}


  
  \begin{abstract}
In this work we relate the known results about the homotopy type of   classifying spaces for smooth foliations, with the homology and cohomology of the discrete group of diffeomorphisms of a smooth compact connected oriented manifold. The   Mather-Thurston Theorem  forms a bridge between the results on the homotopy types of classifying spaces and the various   homology and cohomology groups that are studied. We introduce   the algebraic K-theory of a manifold M that is derived from the discrete group of diffeomorphisms of M, and observe that the calculations of homotopy groups in this work are about these K-theory groups. We include a variety of remarks and open problems related to   the study of the diffeomorphism groups and their homological invariants using the Mather-Thurston Theorem.
\end{abstract}

\maketitle
 
  \tableofcontents

\vfill
\eject

\section{Introduction} \label{sec-intro}

   Let $M$ be a smooth connected  orientable manifold  of dimension $q$.  Let $\Diff_c(M)$ denote the subgroup of compactly supported diffeomorphisms of the topological group $\Diff(M)$,  and let $\Diff_c^+(M)$ and $\Diff^+(M)$ be the respective subgroups of orientation-preserving diffeomorphisms.

For a topological group $G$, we let $G^{\delta}$ denote the group with the discrete topology. The cohomology groups $H^*(BG^{\delta}; \cA)$ are the universal characteristic classes for flat bundles associated with  $G^{\delta}$. For the identity map $\iota \colon G^{\delta} \to G$, there is a natural induced map $\iota^*\colon H^*(BG; \cA) \to H^*(BG^{\delta}; \cA)$. When $G$ is a connected Lie group, 
Milnor calculated this map   for many cases in \cite{Milnor1983}. For example, when the coefficient ring $\cA$ is finite and $G$ is solvable,   the map is an isomorphism. In contrast, when $\cA$ has characteristic $0$ and   $G$ is a compact or complex semisimple Lie group, Milnor  showed that  the map is trivial. One can ask for analogous results in the case of the groups $\Diff_c(M)$  and   $\Diff_c^+(M)$.  

   The classifying space $B\Diff(M)$ classifies $M$-bundles over a base space $X$, and the classifying space $B\Diff(M)^{\delta}$ classifies $M$-bundles with discrete structure group, that is, with a foliation on the total space over $X$  that is transverse to the fibers. The cohomology groups $H^*(B\Diff(M)^{\delta} ; \cA)$, for a coefficient group $\cA$,  give invariants for foliated $M$-bundles, which is one motivation for the study of the topological properties of the spaces $B\Diff(M)^{\delta}$.

The   inclusion map $\iota_{\delta} \colon \Diff_c^+(M)^{\delta} \to \Diff_c^+(M)$ induces a map of the classifying spaces, which has a homotopy theoretic fiber  $B\overline{\Diff_c(M)}$, and there  is  a fibration 
\begin{equation}\label{eq-basicfibration}
B\overline{\Diff_c(M)} \stackrel{\iota}{\longrightarrow} B\Diff_c^+(M)^{\delta} \stackrel{\iota_{\delta}}{\longrightarrow} B\Diff_c^+(M)  \  .
\end{equation}
The space $B\overline{\Diff_c(M)}$ classifies foliations on a product space $X \times M$ which are transverse to the factor $M$, and are a product outside of a compact subset of $X \times M$.
The cohomology of the fiber space $B\overline{\Diff_c(M)}$ is viewed as comparing the cohomology  of the group $\Diff_c^+(M)$ equipped with two topologies, one continuous and the other discrete (see \cite{MostowM1976,Stasheff1978}.)  
The following is the most basic problem,  which remains for the most part unsolved:

\begin{prob}\label{prob-mostbasic}
Calculate the groups
 $H_*(B\overline{\Diff_c(M)}; \cA)$ and $H_*(B\Diff_c^+(M)^{\delta}; \cA)$, and the map  between them, for   $M$  a connected oriented manifold, and $\cA$   an abelian group. 
 \end{prob}

The direct approach to the study of $H_*(B\Diff_c^+(M)^{\delta} ; \mQ)$ is to construct a group action of a finitely generated group $\G$ on $M$, that is, a homomorphism $\Phi \colon \G \to \Diff_c^+(M)$,   and show that the induced map $\Phi_* \colon H_*(B\G ; \mQ) \to H_*(B\Diff_c^+(M)^{\delta} ; \mQ)$ has non-trivial classes in its image. The celebrated work of Thurston in \cite{Thurston1972} did this for the case $M=\mS^1$. The various constructions for higher dimensional manifolds described in Section~\ref{sec-NVsecondary} also use this approach for $M = \mS^q$ a higher dimensional sphere, and for a few other cases where $M$ is a homogeneous space. There are very few other examples where this approach has been carried through. 
Nariman considered in   \cite{Nariman2022} the case where $M = G$ is a connected Lie group, and related $H_2(B\Diff_c^+(G)^{\delta} ; \mQ)$ with the algebraic  K-theory group $K_2(\mC)^+$. 

Another approach to Problem~\ref{prob-mostbasic} is to study the  spectral sequence associated to the fibration \eqref{eq-basicfibration}, which relates the cohomology groups $H^*(B\Diff_c^+(M)  ; \cA)$ and $H^*(B\overline{\Diff_c(M)} ; \cA)$ with the groups $H^*(B\Diff_c^+(M)^{\delta} ; \cA)$.
Nariman   studied this spectral sequence   in \cite{Nariman2022},  with a focus on the induced map on cohomology $\iota_{\delta}^* \colon H^*(B\Diff_c^+(M)  ; \mQ) \to H^*(B\Diff_c^+(M)^{\delta} ; \mQ)$. 
 
 The approach  in this work uses the Mather-Thurston Theorem~\ref{thm-MTT} to reformulate the problem in terms of the homotopy groups of a   space associated to $\Diff_c^+(M)^{\delta}$.  
 The  works of Mather \cite{Mather1975,Mather1979} and Thurston \cite{Thurston1974}  show that there is a  map $\tau \colon B\overline{\Diff_c(M)} \to \cS(M)$   which induces isomorphisms in homology, where   $\cS(M)$ is   the space of compactly supported sections of a bundle with fiber the classifying space $B\oG_q$ of framed smooth codimension-q foliations. 
 The following notion is key to our study of Problem~\ref{prob-mostbasic}:
\begin{defn}\label{def-Ktheory}
The algebraic K-theory groups of $M$ are defined as
\begin{equation}
\cK_{*}(M) = \pi_{*} \left (\left\{B\overline{\Diff_c(M)} \right\}^+ \right) \cong \pi_{*} \left ( \cS(M) \right) \ .
\end{equation}
\end{defn}
Here,  $\ds \left\{B\Diff_c^+(M)_{\delta} \right\}^+$  denotes the space obtained by applying the Quillen +-construction. These are called algebraic K-theory groups, in analogy with  the definition of the algebraic K-theory groups of a ring \cite{Weibel2013}. 
We denote these groups by $\cK_*(M)$ to distinguish them from the usual K-theory of $M$ denoted by $K^*(M)$, and the K-homology denoted by $K_*(M)$.

For $M = \mR^q$, one has    $\cS(\mR^q) \cong \Omega^q B\oG_q$  is a $q$-fold  \emph{based loop} space, and   $\cK_n(\mR^q) \cong \pi_{q+n}(B\oG_q)$. These groups are absolutely enormous, as seen for example in Theorem~\ref{thm-lotsaclasses}. The Mather-Thurston Theorem~\ref{thm-MTT}  says there is an isomorphism $H_*(B\overline{\Diff_c(\mR^q)} ; \mZ) \cong H_*(\cS(\mR^q) ; \mZ)$, and as $\cS(\mR^q)$ is an H-space, its rational homology is a free graded   algebra over the generating set $\cK_*(\mR^q) \otimes \mQ$.

For $M = \mS^q$, one has    $\cS(\mS^q) \cong \Lambda^q B\oG_q$  is a $q$-fold  \emph{free loop} space,  so   $\cK_n(\mS^q) \cong \pi_n(B\oG_q) \oplus \pi_{q+n}(B\oG_q)$. There is an inclusion of the based loops into the free loops, $i \colon \Omega^q B\oG_q \to \Lambda^q B\oG_q$ which induces an inclusion  map $i_* \colon \cK(\mR^q) \to \cK(\mS^q)$. The calculation of the homology groups $H_*(\cS(\mS^q) ; \cA)$ is no longer straightforward, even with rational coefficients $\cA = \mQ$, as it requires unknown properties of the homotopy theory of $B\oG_q$ to calculate the differentials in the minimal model for the space. Haefliger describes the construction of this model in \cite{Haefliger1978a}.

For $M$ a  connected orientable manifold in general, the approach to studying $\cK(M)$ and   $H_*(\cS(M) ; \cA)$ using minimal models encounters even more difficulties.

The strategy of our approach to Problem~\ref{prob-mostbasic} is to construct elements of $\cK_*(M)$ which pair non-trivially with cohomology classes  in $H^*(\cS(M); \mR)$   that are derived from the secondary characteristic classes in $H^*(B\oG_q ; \mR)$. This is used to show that their images under the Hurewicz homomorphism $h \colon \cK_*(M) \to H_*(\cS(M); \mQ) \cong H_*(B\overline{\Diff(M)}  ; \mQ)$ are non-trivial.  Additional arguments in special cases are used to show that some of these classes then map non-trivially into $H_*(B\Diff_c^+(M)^{\delta}; \mQ)$.

This method yields uncountably many independent cycles  in $H_*(B\overline{\Diff_c(M)} ; \mQ)$, each of which is  represented by a homomorphism $\rho \colon \G \to \Diff_c(M)$ of a finitely generated group $\G$, as discussed in Section~\ref{sec-flatbundles}. The constructions  reveal no information on the group $\G$ or how it acts on $M$. In fact,   Freedman estimates in \cite{Freedman2020}   that using the approach   to the proof   of the Mather-Thurston Theorem via foliated surgery  in Meigniez \cite{Meigniez2021}, the number of generators for  $\G$ can be ``exponentially large''.

The results obtained by the approach using the secondary invariants of foliations depend upon the dimension $q$ and the topology of $M$, and are given in detail in later sections of this work. We formulate here some general results, as an overview of the results obtained.
 
 \begin{thm}\label{thm-main1}
  Let $M$ be a smooth   connected  orientable manifold  of dimension $q$. Then there exists a sequence of integers $v_q \geq 1$ which tend to infinity as $q$ tends to infinity, and there exists an uncountably generated  subgroup   $\mV_q \subset \cK_{q+1}(M)$, and  for $1 \leq i \leq v_q$,    linear functionals $\phi_i \colon \cK_{q+1}(M) \to \mR$,  such that  the vector-valued homomorphism 
  \begin{equation}
\vec{\phi} = (\phi_1, \ldots , \phi_{v_q}) \colon \cK_{q+1}(M) \otimes \mQ \ \longrightarrow \ \mR^{v_q}
\end{equation}
is a surjection when restricted to the rational subspace $\mV_q \otimes \mQ$. 
 \end{thm}
\begin{cor}\label{cor-main1}
  Let $M$ be a smooth   connected  orientable manifold  of dimension $q$. Then there exists a sequence of integers $v_q \geq 1$ which tend to infinity as $q$ tends to infinity, and there exists an uncountably generated  subgroup    $\mH_q \subset H_{q+1}(B\Diff_c^+(M)^{\delta}; \mQ)$, and  for $1 \leq i \leq v_q$, cohomology classes $\alpha_i \in H^{q+1}(B\Diff_c^+(M)^{\delta}; \mR)$,  such that  the vector-valued homomorphism 
  \begin{equation}
\vec{\alpha} = (\alpha_1, \ldots , \alpha_{v_q}) \colon H_{q+1}(B\Diff_c^+(M)^{\delta}; \mQ) \ \longrightarrow \ \mR^{v_q}
\end{equation}
is a surjection when restricted to the rational subspace $\mH_q \otimes \mQ$. 
 \end{cor}

When the manifold $M$ has trivial tangent bundle, there are homotopy operations on $\cK_*(M)$ which generate many more non-trivial classes. These are the brace products described in Proposition~\ref{prop-brace}, and the results now apply to the space $B\overline{\Diff_c(M)}$. In Section~\ref{sec-q>3},  integers $\widehat{v}_{q,k}$ and $\kappa_q$   are defined based on the number of independent secondary classes which are variable. We always have $\widehat{v}_{q,k} = v_q$ when $k=2q+1$, and for $q \geq 3$ we have $\widehat{v}_{q,k} = v_q$ for $k=2q+4$, and both indices tend to infinity as $q$ tends to infinity. Other values of $k$ are also tending to infinity, as the dimension $q$ increases.
Note that in the following theorem, we now assume that $M$ is compact.
 \begin{thm}\label{thm-main2}
  Let $M$ be a smooth compact connected   manifold  of dimension $q$ with trivial tangent bundle. Then there exists a sequence of integers $\widehat{v}_{q,k}$, for $k \geq  2q+1$,     
  and there exists an uncountably generated  subgroup   $\mV_{q,k} \subset \cK_{k}(M)$, and  for $1 \leq i \leq \widehat{v}_{q,k}$,    linear functionals $\widehat{\phi}_i \colon \cK_{k}(M) \to \mR$,  such that  the vector-valued homomorphism 
  \begin{equation}
\vec{\phi} = (\widehat{\phi}_1, \ldots , \widehat{\phi}_{v_{\widehat{v}_{q,k}}}) \colon \cK_{k}(M) \otimes \mQ \ \longrightarrow \ \mR^{\widehat{v}_{q,k}}
\end{equation}
is a surjection when restricted to the rational subspace $\mV_{q,k} \otimes \mQ$. 
 \end{thm}
  Then again using the  Mather-Thurston Theorem, we obtain:
 \begin{cor}\label{cor-main2}
  Let $M$ be a smooth compact connected   manifold  of dimension $q$ with trivial tangent bundle.
  Then there exists a sequence of integers $\widehat{v}_{q,k}$, for $k \geq  2q+1$,  and there exists an uncountably generated  subgroup    $\mH_{q,k} \subset H_{k}(B\overline{\Diff(M)} ; \mQ)$, and  for $1 \leq i \leq \widehat{v}_{q,k}$, cohomology classes $\widehat{\alpha}_i \in H^{k}(B\overline{\Diff(M)} ; \mR)$,  such that  the vector-valued homomorphism 
  \begin{equation}
\vec{\alpha} = (\widehat{\alpha}_1, \ldots , \widehat{\alpha}_{\widehat{v}_{q,k}}) \colon H_{k}(B\overline{\Diff(M)} ; \mQ) \ \longrightarrow \ \mR^{\widehat{v}_{q,k}}
\end{equation}
is a surjection when restricted to the rational subspace $\mH_{q,k} \otimes \mQ$. 
 \end{cor}

There are many further constructions   of non-trivial families of classes in $H_{*}(B\overline{\Diff_c(M)} ; \mQ)$ beyond those above, obtained using the techniques developed in Theorems~\ref{thm-method2} and \ref{thm-method3}. A selection of these constructions is given in later sections.

We give a brief overview of the contents of the work. 
Section~\ref{sec-classifying} recalls the definition of the spaces $B\oG_q$ and   $B\G_q$ which classify smooth foliations of codimension $q$.   

Section~\ref{sec-MT}   introduces the space $\cS(M)$ and discusses some basic properties about its homotopy type.
Then in  Section~\ref{sec-constructions}   we give some general constructions of non-trivial   classes in $\cK_*(M)$, assuming the existence of spherical cohomology classes in  $H^*(B\oG_q ; \mR)$.   
 
The construction of the secondary characteristic classes in $H^*(B\G_q^+ ; \mR)$ and $H^*(B\oG_q ; \mR)$ is reviewed  in Section~\ref{sec-secondary}, and in Section~\ref{sec-NVsecondary} we recall the examples from the literature showing that many of these classes are non-trivial. 
 Then in Section~\ref{sec-spherical}, we apply the results of Sections~  \ref{sec-constructions}, \ref{sec-secondary}  and \ref{sec-NVsecondary}  to the show the existence of non-trivial spherical cohomology classes in  $H^*(B\oG_q ; \mR)$.
 The current understanding of the topological types of $B\oG_q$ and $B\G_q$ as described in Section~\ref{sec-spherical}  remains  essentially unchanged since the 1980's, as described in  \cite{Hurder1985a}, and also the  survey \cite{Hurder2008}.

 Constructions of flat $\Diff_c^+(M)$-bundles using the results of Section~\ref{sec-spherical} and the ideas from the proof of the Mather-Thurston Theorem are given in Section~\ref{sec-flatbundles}. In particular, we consider which of these homology classes in $H_*(B\overline{\Diff_c(M)} ; \mQ)$ map non-trivially into $H_*(B\Diff_c^+(M)^{\delta} ; \mQ)$.
 
In Section~\ref{sec-q=1}, we consider what is known for Problem~\ref{prob-mostbasic}  when $q=1$, so either $M = \mR$ or $M = \mS^1$. Morita  has obtained the most comprehensive results in this case  in the publications \cite{Morita1984a,Morita1985}.

In Section~\ref{sec-q=2}, we consider the case when $q=2$, and the special cases when $M=\mR^2$, $\mS^2$,    $\mT^2$ and for closed surfaces $\Sigma_g$ with genus $g \geq 2$. Our analysis uses the results of Rasmussen in \cite{Rasmussen1980}.  The case when $M = \Sigma_g$ is a compact surface of genus $g \geq 2$ has been studied in depth in the works of Bowden \cite{Bowden2012} and Nariman \cite{Nariman2017b,Nariman2020}, and we   recall some of their results as well.

In Section~\ref{sec-q=3}, we consider the case for orientable $3$-manifolds, where the use the the James brace product descibed in Appendix~\ref{subsec-brace} now enters into the constructions.   
In Section~\ref{sec-q>3}, we describe the results for higher dimensional manifolds. 

For dimensions for $q=2$ and $3$, the structure of the groups $\Diff_c^+(M)$ is mostly well-understood, but for higher dimensions these topological groups become increasingly mysterious. 
There have also been many   recent developments in the study  of the  homology and cohomology     of the group of diffeomorphisms for manifolds with dimensions $q \geq 5$, as discussed by   Kupers  \cite{Kupers2019a}. Understanding the relations between  these classes and the groups $H^*(\Diff_c^+(M)^{\delta} ; \mQ)$ is a fundamental question, which is mostly unsolved except in special cases.  Problems~\ref{prob-FH} and \ref{prob-Waldhausen} at the end of this section were one of the  original motivations for this work.

  In the author's thesis \cite{Hurder1980}, a variety of basic techniques from homotopy theory of spaces and fibrations were used to extract homotopy results from the non-vanishing of the secondary invariants. These were recalled and extended in the work \cite{Hurder1985a}, and the applications to the study the homology of diffeomorphism groups was developed in   the unpublished  manuscript \cite{Hurder1985b}. The formulation  of the K-theory of a smooth manifold   $\cK_*(M)$  was   given    in the author's 1987 lectures in Berlin \cite{Hurder1987} as applications of these calculations of $\pi_*(B\oG_q)$. 
  Throughout this text, we liberally cite references to works on foliation theory and related topics, so that the bibliography is extensive as a result.

 The author thanks Sam Nariman for helpful comments on this manuscript.

\section{Foliation classifying spaces}\label{sec-classifying}

We give an overview of the definitions and basic properties of the classifying spaces for foliations introduced by Haefliger in  \cite{Haefliger1970,Haefliger1971}. All foliations and maps   are assumed to be smooth.

 Let $\Gamma_q$ denote the pseudogroup generated by local $C^{\infty}$-diffeomorphisms of $\mR^q$. A codimension $q$-foliation $\F$ on $M$, along with a good open covering of  a paracompact manifold $M$, defines a 1-cocycle with values in $\Gamma_q$. 
Haefliger observed in \cite{Haefliger1970} that there is a classifying space $B\Gamma_q$ for these $\Gamma_q$-cocycles, so 
that the foliation $\F$ induces   a continuous map $\chi_{\F} \colon M \to B\Gamma_q$ whose homotopy class is independent of the choice of the good covering.  One says that  $\chi_{\F}$ is the classifying map for $\F$. 
The Jacobian map of a local diffeomorphism induces a universal map $\nu \colon B\Gamma_q \to B{\rm O}(q)$, which classifies the normal bundle on the universal foliated bundle over $B\Gamma_q$.

There are   alternate models for the classifying space $B\Gamma_q$ which are all weak homotopy equivalent,      by   Greenberg \cite{Greenberg1985},   Jekel \cite{Jekel1984},   McDuff \cite{McDuff1979a},   Moerdijk \cite{Moerdijk1991}, and   Segal \cite{Segal1978}. 
Any of these models suffice for ``classifying'' foliations of codimension-q on paracompact manifolds.

The work by Bott, Shulman and Stasheff \cite{BSS76} (see also the text by Dupont \cite{Dupont1978})
  constructed a semi-simplicial model for $B\Gamma_q$ and a semi-simplicial de~Rham complex of forms whose   cohomology groups are isomorphic to  $H^*(B\Gamma_q , {\mathbb R})$. The existence of a de~Rham complex for $B\Gamma_q$ allows one to use de~Rham homotopy methods in the study of the topology of the space, as for example in \cite{Bott1978}.
 
 There are various alternate subclasses of groupoids that arise for foliations with some additional structure,  as discussed by Haefliger in \cite{Haefliger1971}. We consider here 
   the pseudogroup $\Gamma_q^+$ of orientable   local $C^{\infty}$-diffeomorphisms of $\mR^q$ and its classifying space $B\Gamma_q^+$. There is a natural   map   $\nu \colon  B\Gamma_q^+ \to B{\rm SO}(q)$ classifying the normal bundle  on $B\Gamma_q^+$. Also considered in the literature are groupoids of germs of maps preserving a volume form on $\mR^q$, or preserving a symplectic form, or a Riemannian metric, to cite some of the cases which have been studied.
   
  There is a canonical basepoint $w \in B\Gamma_q^+$ given by the image of the point foliation of ${\mathbb R}^q$. 
   Let    $B\overline{\Gamma}_q$ denote the homotopy fiber of  $\nu \colon  B\Gamma_q^+ \to B{\rm SO}(q)$ over the basepoint $w_0 = \nu(w)  \in B{\rm SO}(q)$.   Then $B\overline{\Gamma}_q$ is   realized as the pull-back in the diagram \eqref{eq-commutingdiagram}:
\begin{align}\label{eq-commutingdiagram}
\xymatrix{
{\rm SO}(q) \ar[d]^{\simeq} \ar[r]^{\iota} &  B\overline{\Gamma}_q\ar[d]^{\overline{\nu}} \ar[r]^{\iota}
& \hspace{2mm} B\Gamma_q^+  \ar[d]^{\nu}  \\
{\rm SO}(q) \simeq \Omega(B{\rm SO}(q)) \ar[r]^{\iota}& {\mathcal P}B{\rm SO}(q)   \ar[r]^{eval}
& \hspace{2mm}    B{\rm SO}(q)    
} 
\end{align}
The map $\overline{\nu}$ induces a trivialization of the universal normal bundle over $B\overline{\Gamma}_q$ so that $B\overline{\Gamma}_q$ can also be interpreted as the classifying space of foliations with trivialized normal bundles. Haefliger observed in \cite[Theorem~3]{Haefliger1970} that this   implies   the space $B\overline{\Gamma}_q$ is $(q-1)$-connected, as the Gromov-Phillips immersion theory \cite{Gromov1969,Phillips1968,Phillips1969} implies there is a unique foliation by points on ${\mathbb R}^q$. More is true, as Haefliger showed that $\pi_q(B\overline{\Gamma}_q)$ is trivial, that there is a unique homotopy class of $\Gamma_q$-structures on a foliated microbundle over ${\mathbb S}^q$. 
\begin{prob}\label{prob-connected}
Show that   the space  $B\overline{\Gamma}_q$ is $2q$-connected, or  
  determine the least value  of $k \leq 2q$ for which $\pi_k(B\overline{\Gamma}_q)$ is non-trivial. 
\end{prob}
The best partial result  to date  on this problem is Theorem~\ref{thm-MT} below.

 \section{Mather-Thurston Theory}\label{sec-MT}
 
    Let $M$ be a connected smooth orientable manifold  of dimension $q$.  We recall the Mather-Thurston Theory, and then introduce the algebraic K-theory groups associated to this theory.
    
     The space $B\overline{\Diff_c(M)}$ classifies foliations on a product space $X \times M$, for $X$ compact,  which are transverse to  the second factor $M$, and are horizontal outside of a compact subset of $X \times M$.
Thus, there is a natural classifying map 
\begin{equation}\label{eq-fibration2}
 \tilde{\nu}_M \colon B\overline{\Diff_c(M)} \times M \to B\Gamma_q^+
\end{equation}
Let $\nu_{M} \colon M \to B{\rm SO}(q)$ denote the classifying map of the tangent bundle to $M$. Then we have a commutative diagram:
 
\begin{picture}(200,60) 
\put(240,50){$B\Gamma_q^+ $}
\put(250,40){\vector(0,-1){17}}
\put(188,23){\vector(3,2){30}}
\put(260,30){$\nu$}
\put(192,42){$\tilde{\nu}_M$}
\put(210,20){$\nu_{M}$}
\put(130,10){$B\overline{\Diff_c(M)} \times M$}
\put(210,10){$ \longrightarrow$}
\put(235,10){$B{\rm SO}(q)  $}
\end{picture}

Let $\sigma_M \colon M \to B\G_q^+$ be the classifying map of the point foliation on $M$. Then   $\cS(M)$ denotes the space of   compactly supported liftings of the map $\nu_{M} \colon M \to B{\rm SO}(q)$ to $B\G_q^+$. 
That is, $ g \colon X \to \cS(M)$ is a map $g \colon X \times M \to B\G_q^+$ such that $\nu \circ g = \nu_M$, and $g = \sigma_M \colon X \times M \to B\G_q^+$   outside of a compact subset of $X \times M$. (See \cite[Section~6]{Mather1979}.)
 Then by the works of Mather and Thurston \cite{Mather1975,Mather1979,Thurston1974}:

 \begin{thm} \label{thm-MTHT}
 The   adjunct   map $\tau \colon B\overline{\Diff_c(M)} \to \cS(M)$ induces an isomorphism in homology.
 \end{thm}
 For the special case $M=\mR^q$, the mapping $\widetilde{\nu}_M$ can be considered as being  defined on the $q$-fold  suspension $\Sigma^q B\overline{\Diff_c(\mR^q)}$, and so its adjoint induces a map from $B\overline{\Diff_c(\mR^q)}$ to the iterated loop space on $B\overline{\Gamma}_q$. As a special case, there is the remarkable consequence of Theorem~\ref{thm-MTHT}:
    \begin{thm}\label{thm-MTT}
The   adjunct map   $\tau \colon B\overline{\Diff_c(\mR^q)} \to \Omega^q B\overline{\Gamma}_q$   is a homology isomorphism. 
\end{thm}

 There are several   proofs of   Theorem~\ref{thm-MTHT} in the literature, beyond the initial works of Mather and Thurston \cite{Mather1975,Mather1979,Thurston1974}.  McDuff gave a proof using categorical methods in \cite{McDuff1979a}, following ideas of Segal in \cite{Segal1978}. Mitsumatsu and Vogt  gave a geometric proof in \cite{MitsumatsuVogt2017} for foliations whose leaves have dimension $2$.  Meigniez gave a    proof      in  \cite{Meigniez2021} using   bordism  techniques, for $q \geq 2$.

The Mather-Thurston  Theorem~\ref{thm-MTT} is key to the proof of the following fundamental result, due to Mather \cite{Mather1973,Mather1975} for $q=1$, and to  Thurston \cite{Thurston1974} for $q > 1$.   

\begin{thm}[Mather \cite{Mather1975,Mather1979}, Thurston \cite{Thurston1974}]\label{thm-MT}
The space $B\oG_q$ is $(q+1)$-connected. That is,  for   $q \geq 1$,  $\pi_k(B\overline{\Gamma}_q) = 0$ for all $0 \leq k \leq q+1$.
\end{thm}
\begin{cor}\label{cor-simplyconnected}
 The space  $\cS(M)$ is simply connected.
\end{cor}

  Let $X$ be a connected topological space, and assume that $\pi_1(X)$ is a perfect group.  Let $X^+$ denote the space obtained by applying the Quillen +-construction to $X$; see \cite{Grayson2013} or \cite[Chapter IV]{Weibel2013}. Then $X^+$ is a   connected space with  a natural map $Q \colon X \to X^+$ which induces isomorphisms on all homology groups.   The space $X^+$ with these properties is unique up to weak homotopy equivalence.
  
Given a smooth connected  manifold $M$ of dimension $q$, the adjunct map $\tau \colon B\overline{\Diff_c(M)} \to \cS(M)$   is a homology isomorphism, and $\cS(M)$ is simply connected.
Thus, there is a weak homotopy equivalence between  $\cS(M)$  and $\{ B\overline{\Diff_c(M)}\}^+$, and so these spaces have isomorphic homotopy groups. 
We thus have the equality stated in Definition~\ref{def-Ktheory}, 
 \begin{equation}
\cK_{\ell}(M) \ = \ \pi_{\ell} \left (\left\{B\overline{\Diff_c(M)} \right\}^+ \right) \cong \pi_{\ell} ( \cS(M) )  \ , \ {\rm for \ all} \ \ell \geq 1 \ .
\end{equation}
 
 We thus can formulate the problem addressed in this work:
 \begin{prob}
 For a smooth connected  manifold $M$ of dimension $q$, construct classes which map to non-zero classes in the image of the Hurewicz map
 \begin{equation}\label{eq-rhtadjunct}
h_* \colon  \pi_{\ell} ( \cS(M) ) \longrightarrow H_{\ell}(\cS(M) ; \mQ) \ .
\end{equation}
 \end{prob}
 We will first introduce techniques for constructing non-trivial classes in  $\cK_*(M)$, and then study which of these classes pair non-trivially with    cohomology classes  in   $H^*(\cS(M) ; \mQ)$. The second step is critical, as the results of  Section~\ref{sec-spherical} imply that  the groups $\pi_{\ell} ( \cS(M) )$ are uncountable for an infinite range of $\ell \to \infty$, but only some of these classes are known to yield non-trivial homology classes.  
  Thus, the refined strategy    for constructing non-trivial classes in $H_*(B\overline{\Diff(M)} ; \mQ)$ to first construct elements of $\cK_*(M)$ using techniques from rational homotopy theory and properties of the homotopy type of $B\oG_q$, then showing the  secondary classes in $H^*(B\oG_q ; \mR)$    detect these classes.   Theorem~\ref{thm-MTHT} implies these calculations yield  non-trivial  homology classes defined by flat $\Diff_c^+(M)$-bundles, and as discussed in Section~\ref{sec-flatbundles}, one can show that some of these classes are nontrivial in the image of the map $\iota_* \colon H_*(B\overline{\Diff(M)} ; \mQ)  \to H_*(B\Diff_c^+(M)^{\delta} ; \mQ)$.

 \section{Calculations of the K-theory groups}\label{sec-constructions}

   We first make some general remarks about the calculation of the groups $\cK_*(M)$, then give three methods for constructing classes homotopy classes in $\pi_*(\cS(M))$ which are detected by cohomology classes in $H^*(\cS(M) ; \mR)$.
   A major difficulty with the study of $\cS(M)$ is  that there is very incomplete knowledge of the topological type of the space $B\oG_q$, so one typically works with an approximation to the fiber group $B\oG_q$ which suffices to construct non-trivial classes in $H_*(\cS(M) ; \mQ)$..

    If the classifying map $\nu_{M} \colon M \to B{\rm SO}(q)$ is not trivial,   then the calculation of the   cohomology  of the mapping space $\cS(M)$ can be quite complicated. The works by Haefliger \cite{Haefliger1978a,Haefliger1982}, by    Shibata \cite{Shibata1984a,Shibata1984b} and by Tsujisita \cite{Tsujishita1977,Tsujishita1981} use rational homotopy techniques to study analogs of $\cS(M)$.

If the tangent bundle to $M$ is trivial, so the map classifying $\nu_M$ can be assumed to be  the constant map to the basepoint $w_0 \in BSO(q)$, 
 then  $ \cS(M) \cong \Maps_c(M ,  B\oG_q )$ which greatly simplifies the study of $\cS(M)$.
The study of a space of maps is a classical topic in topology, as discussed by Cohen and Taylor in \cite{CohenTaylor1978,CohenTaylor1983} and   by Bendersky and Gitler in \cite{BenderskyGitler1991}; see also F\'elix and Tanr\'e \cite{FelixThomas2004} which compares these two approaches.     
 
 The case when $M$ is an open manifold, so not compact, then there is a weak homotopy equivalence $\cS(M) \cong \cS(\widehat{M}, x_{\infty})$ where $\widehat{M}$ denotes the 1-point compactification of $M$, and $x_{\infty}$ is the point at infinity, which is taken as the basepoint. Then $\cS_{\infty}(\widehat{M}, x_{\infty})$ denotes the sections which map the basepoint $x_{\infty}$ to the basepoint $\sigma_M \in B\G_q^+$. That is, for a compact space $X$, the $\G_q^+$ structure induced on $X \times M$ by a map $f \colon X \to \cS(\widehat{M}, x_{\infty})$ is that of the horizontal product foliation outside of a compact subset of $X \times M$. In the case when $M= \mR^q$, which was discussed in the introduction, we obtain $\cS(\mR^q) \cong \Maps(\mS^q, B\oG_q, x_{\infty}) \cong \Omega^q B\oG_q$. For other manifolds, such as when $M$ is the interior of a compact manifold with boundary, then other interesting calculations arise.
 
We now give some methods for analyzing the space $\cS(M)$. 
Given a non-trivial class $\phi$ in $H^*_{\pi}(B\G_q^+ ; \mR)$ or $H^*_{\pi}(B\oG_q ; \mR)$, we construct classes in $\cK_{*-k}(M)$ which are detected by  cohomology classes in $H^{*-k}(\cS(M) ; \mR)$,   where the degree is decreased by an integer $k$ which depends   the construction.

Recall the following  standard notion.

 \begin{defn}\label{def-spherical}
 Let $X$ be a topological space. A class $ \phi \in H^n(X; \mR)$ is said to be \emph{spherical} if there exists a map $f \colon \mS^n \to X$ such that $f^*(\phi) \in H^n(\mS^n ; \mR)$ is non-trivial. Let $H^*_{\pi}(X ; \mR)$ denote the subspace of spherical cohomology classes.
 \end{defn}
 
 Here is the first method of constructing classes.
  \begin{thm}\label{thm-method1}
  Let $M$ be a  connected  oriented manifold of dimension $q$. 
   Let $\phi \in H^n_{\pi}(B\G_q^+; \mR)$ be a spherical class, for $n > 2q$.
   Then for $\wtalpha =  \int_M    \widetilde{\nu}^*(\phi)  \in H^{n-q}(\cS(M) ; \mR)$ there exists $[f_0] \in \cK_{n-q}(M)$ such that the evaluation
   $\langle \wtalpha, [f_0] \rangle \ne 0$. 
  \end{thm}
  \proof
   Let $\F$ denote the product foliation on $N = \mS^{n-q} \times M$ with leaves $\mS^{n-q} \times \{y\}$ for $y \in M$.
Let $\widetilde{\nu}_{\F} \colon   N \to B\G_q^+$ classify  the foliation $\F$, which is a lifting of the   map  $\nu_{\F} \colon N \to B{\rm SO}(q)$     classifying   the normal bundle of $\F$. We construct a compactly supported perturbation of $\widetilde{\nu}$.

Choose basepoints  $x_0 \in M$   and    $\theta_{n-q} \in \mS^{n-q}$, and set $N = \mS^{n-q} \times M$ with basepoint $(\theta_{n-q}, x_0) \in N$. 
Set $w_0 = \nu_{\F}(\theta_{n-q}, m_0)  \in B{\rm SO}(q)$ and set the basepoint $w = \widetilde{\nu}_M(\theta_{n-q} , m_0) \in B\G_q^+$. 
Let  $\iota_w \colon B\oG_q \to B\G_q^+$ denote the inclusion of the fiber over the basepoint $w_0$.  Then $w \in B\oG_q \subset B\G_q^+$.  

  Let $f' \colon \mS^n \to B\G_q^+$  be a map such that   $\langle \phi , \widehat{f'} \rangle \ne 0$. The group $\pi_n(B{\rm SO}(q)) \otimes \mQ = 0$ as $n > 2q$, so there exists a multiple of $f'$ in the image of the map $\pi_n(B\oG_q) \to \pi_n(B\G_q^+)$. Thus  we may assume that $f' = \iota_w \circ f$ for $f \colon \mS^n \to B\oG_q$ and  that $f$ is a pointed map, with $f(\theta_{n-q}, m_0) = w$.

Choose a   contractible disk neighborhood $\mD^n \cong U \subset N$ with the basepoint $(\theta_{n-q} , m_0)  \in \partial U$, where $\partial U$ is homeomorphic to $\mS^{n-1}$.   Let $N/U$ denote the quotient space of $N$ by the set  $U$, then $N$ and  $N/U$ are homeomorphic. Note that the restriction of $\widetilde{\nu}_{\F} | U \colon U \to B\G_q^+$ is homotopic to the constant map to the basepoint $w \in B\G_q^+$.  Also, the quotient $U/\partial U \cong \mS^n$ so collapsing the boundary $\partial U$ yields  a map $b \colon N \to \mS^n \vee N$. In colloquial terms, this process ``bubbles off'' a copy of $\mS^n$ from  $N$.

Define  $\widehat{f} = (\widehat{f}_1 \vee \widehat{f}_2) \circ b \colon N \to B\G_q^+$ as the join of two maps, where     $\widehat{f}_1 = \iota_w \circ f \colon \mS^n \to B\oG_q \to B\G_q^+$ and $\widehat{f}_2 = \widetilde{\nu}_{\F} \circ \pi  \colon N \to N/U \cong N \to  B\G_q^+$.
 Calculate the pairing of the cohomology class $\phi$ with the  homology class of $\widehat{f}$:
\begin{equation}
\langle \phi, [\widehat{f}] \rangle = \langle \phi, [\widehat{f}_1] \rangle + \langle \phi , [\widehat{f}_2] \rangle = \langle \phi, [\widehat{f}_1] \rangle  = \langle \iota_w^*(\phi), [f] \rangle \ne 0 \ ,
\end{equation}
where we use that $\widehat{f}_2 \colon \colon N \to B\G_q^+$ is homotopic to a map into the $q$-skeleton, so the evaluation of the class $\phi$ vanishes as its degree $n > q$.
 Thus the pull-back class $\widehat{f}^*(\phi) \in H^n(N ; \mR)$ is non-vanishing.
 The map   $\widehat{f}  \colon N = \mS^{n-q} \times M \to B\G_q^+$ is a lift of the map $\nu_M \colon N \to B{\rm SO}(q)$, so induces   maps $f_0 \colon \mS^{n-q} \to \cS(M)$ and $f_1 \colon  \mS^{n-q} \times M \to \cS(M) \times $.
 We then have the commutative diagram in cohomology, where the maps $\int_M$ are integration over the fiber (see Appendix~\ref{subsec-integration}):

\begin{picture}(300,64)(-30,0)  
\put(24,50){$H^n(\mS^{n-q} \times M ; \mR)$}
\put(130,53){\vector(-1,0){20}}
\put(118,60){$f_1^*$}

\put(72,42){\vector(0,-1){17}}
\put(80,30){$\int_M$}
\put(36,10){$H^{n-q}(\mS^{n-q} ; \mR)$}
\put(130,14){\vector(-1,0){20}}
\put(118,20){$f_0^*$}
\put(140,50){$H^n(\cS(M) \times M ; \mR)$}
\put(264,53){\vector(-1,0){30}}
\put(180,42){\vector(0,-1){17}}
\put(186,30){$\int_M$}
 \put(250,60){$\tilde{\nu}^*$}
\put(280,50){$H^n(B\G_q^+ ; \mR)$}
\put(150,10){$H^{n-q}(\cS(M) ; \mR)$}
\end{picture}

 Then $\widehat{f} = \widetilde{\nu} \circ f_1$ so $\widehat{f}^*(\phi) \ne 0$  implies that 
 $\ds  f_0^* \int_M    \widetilde{\nu}^*(\phi)  =  \int_M  f_1^* \circ \widetilde{\nu}^*(\phi) \ne 0$.
   Set $\wtalpha =  \int_M    \widetilde{\nu}^*(\phi)  \in H^{n-q}(\cS(M) ; \mR)$ which then pairs nontrivially with the class $[f_0] \in \pi_{n-q}(\cS(M))$.
  \endproof

 \begin{remark}\label{rmk-Thurston}
 {\rm
 The construction in the above proof is equivalent to a well-known construction of Thurston which  first appeared in      \cite{Thurston1974}. 
 Given a basepoint preserving  map  $f \colon \mS^n \to B\oG_q$, it determines a delooped map $f_0 \colon \mS^{n-q} \to \Omega^q B\oG_q$ and a corresponding map $f_1 \colon \mS^{n-q} \times \mR^q \to B\oG_q$ which maps the complement of a compact subset to the basepoint of $B\oG_q$.
 Choose a basepoint  $x_0 \in M$, and choose  disk neighborhood $x_0 \in U \subset M$ and a diffeomorphism $\psi \colon \mR^q \cong \mD^q \to U \subset  M$. Then $\psi$ induces a map $f_{\psi} \colon \mS^{n-q} \times \mR^q \to \cS(M) \times M$ which is just the map $f_1$  in the proof.
}
\end{remark}

For the next   results, we assume  that   the tangent bundle to $M$ is framed, so the classifying map $\nu_M \colon M \to B{\rm SO}(q)$ is the constant map, and the framing defines a lifting  $\widetilde{\nu}_M \colon M \to B\oG_q$ to the fiber. Then $\cS(M) = \Maps_c (M, B\oG_q)$ is a space of maps with basepoint $\widetilde{\nu}_M$, and there are additional tricks for constructing non-trivial homotopy classes in $\cK_*(M)$.  We first give a construction which requires the additional assumption that $M$ is compact.
 
  \begin{thm}\label{thm-method2}
  Let $M$ be a  compact connected   oriented manifold of dimension $q$, and assume the tangent bundle $TM$ is trivial. 
   Let $\phi \in H^n_{\pi}(B\oG_q; \mR)$ be a spherical class, for $n \geq q+2$.
   For $\wtbeta =  \sigma^* \circ \widetilde{\nu}^*(\phi)  \in H^{n}(\cS(M) ; \mR)$ there exists $[g_0] \in \cK_{n}(M)$ such that the evaluation
   $\langle \wtbeta, [g_0] \rangle \ne 0$. 
  \end{thm}
\proof

Let $\phi \in H^n(B\oG_q; \mR)$ be a spherical class, which pairs non-trivially with the homology class defined by a map  $f \colon \mS^n \to B\oG_q$. Assume that $f$ maps  $\theta_n \in \mS^n$   to  $w \in B\oG_q$.  
Let  $\widehat{f}  \colon N = \mS^n \times M \to B\oG_q$ be the induced map on the product, which is a compact perturbation of $\widetilde{\nu}_M$ as $M$ is compact. 

Then $\widehat{f}^*(\phi) \in H^n(\mS^n \times M ; \mR)$ is non-vanishing, and supported on the first factor of the decomposition $H^*(\mS^n \times M ; \mR) \cong H^*( \mS^n  ; \mR) \otimes H^*(M ; \mR)$.  
The map   $\widehat{f}$ is a lift of the constant map $\nu_{\F} \colon N \to B{\rm SO}(q)$, and as $M$ is compact, $\widehat{f}$  induces   maps $f_0 \colon \mS^{n-q} \to \cS(M)$ and $f_1 \colon  \mS^{n-q} \times M \to \cS(M)$.

Let $\sigma_1 \colon \mS^n \to \mS^n \times \{x_0\} \subset \mS^n \times M$ be the inclusion of the first factor. 
Then   we have a    commutative diagram in cohomology:
 
\begin{picture}(300,64)(-30,0)  
\put(24,50){$H^n(\mS^{n} \times M ; \mR)$}
\put(134,53){\vector(-1,0){30}}
\put(118,60){$f_1^*$}

\put(62,42){\vector(0,-1){17}}
\put(70,30){$\sigma_1^*$}
\put(36,10){$H^{n}(\mS^{n} ; \mR)$}
\put(134,14){\vector(-1,0){30}}
\put(118,20){$f_0^*$}
\put(140,50){$H^n(\cS(M) \times M ; \mR)$}
\put(264,53){\vector(-1,0){30}}
\put(180,42){\vector(0,-1){17}}
\put(186,30){$\sigma_1^*$}
 \put(250,60){$\tilde{\nu}^*$}
\put(280,50){$H^n(B\oG_q ; \mR)$}
\put(150,10){$H^{n}(\cS(M) ; \mR)$}
\end{picture}

  Then $\widehat{f} = \widetilde{\nu} \circ f_1$ so $\widehat{f}^*(\phi) \ne 0$  implies   
 $\ds  f_0^* \circ  \sigma_1^* \circ  \widetilde{\nu}^*(\phi)  =  \sigma_1^* \circ   f_1^* \circ \widetilde{\nu}^*(\phi) \ne 0$.
   Set $\wtbeta =  \sigma_1^* \circ  \widetilde{\nu}^*(\phi)  \in H^{n-q}(\cS(M) ; \mR)$ which then pairs nontrivially with the class $[f_0] \in \pi_{n}(\cS(M))$.
\endproof

 \begin{remark}\label{rmk-section}
 {\rm
 Note that the construction of the classes $[f_0]$ and $\wtbeta$ depends on the existence of the classifying map $\widehat{f}  \colon N = \mS^n \times M \to B\G_q^+$ that has image in the fiber $B\oG_q$. When the tangent bundle of $M$ is not trivial, the map $\widehat{f}$ has to cover the classifying map $\nu_{\F} \colon N \to B{\rm SO}(q)$ which is not trivial. Thus, the fiber over the image of $\nu_{\F}$ need not be constant, and it is no longer obvious how to construct the lifted map $\widehat{f}$. In fact, the twisting of the fiber over $\nu_{\F}$ may result in a twisting of the chains in $B\G_q^+$ which may yield  a non-trivial differential in the spectral sequence for the fibration $B\oG_q \to B\G_q^+ \to B{\rm SO}(q)$  killing an  extension of the  class $\phi \in H^n(B\oG_q ; \mR)$ to $ H^n(B\G_q^+ ; \mR)$.  
 }
\end{remark}

When $X$ is not compact, recall that $\widehat{X}$ denotes the $1$-point compactification of $X$, and otherwise $\widehat{X} = X$. 
  For the next constructions, we  use the   following   notion. 
   \begin{defn}\label{def-cospherical}
 Let $X$ be a connected topological space. A homology class $C \in H_k(\widehat{X}; \mZ)$ is said to be \emph{co-spherical} if there exists a   map $\xi \colon  \widehat{X} \to \mS^k$ such that $\xi_*(C)   \in H_k(\mS^k ; \mZ)$ is non-trivial.
 Let $H_*^{\pi}(\widehat{X} ; \mZ)$ denote the subgroup of co-spherical classes.
  \end{defn}
The fundamental class of $M$ is co-spherical, as seen by the mapping $\xi \colon M \to \mS^q$ of $M$ to its  quotient   by its $(q-1)$-skeleton, and if $M$ is non-compact, composed with the collapse of a neighborhood of infinity to the basepoint, which is then homotopy equivalent to $\mS^q$. This map was the key to the construction in the proof of Theorem~\ref{thm-method1}.

The next construction uses a combination of the ideas behind the proofs of Theorems~\ref{thm-method1} and \ref{thm-method2}.  
We assume there is given a  co-spherical class    $[C] \in H^{\pi}_k(\widehat{M}; \mZ)$        for some $0 < k < q$, and use this to make a compact perturbation of the base map $\widetilde{\nu}_M$.   However, the perturbation is no longer assumed to be on a contractible disc in $N = \mS^k \times M$, so an additional assumption is needed on the classifying map for the tangent bundle over  the perturbation. The simplest assumption is that $M$ has trivialized tangent bundle, but it suffices to assume that the tangent bundle is trivial over the compact support of the perturbation. 
   For example, for $M$ compact, the homology  $H_1(M ; \mQ)$ is generated by co-spherical classes, as 
  $\ds H^1(M ; \mZ) \cong [M, K(\mZ, 1)] \cong  [M, \mS^1]$, and a tubular neighborhood of each generator $\mS^1 \subset M$ will have trivial tangent bundle. We formulate the next result with the assumption that $TM$ is trivial, and leave extensions of the method to more general  cases as exercises. Theorem~\ref{thm-g=2}  illustrates the application of these remarks.
     \begin{thm}\label{thm-method3}
  Let $M$ be a connected   oriented manifold of dimension $q$, and assume the tangent bundle $TM$ is trivial. 
   Let $\phi \in H^n_{\pi}(B\oG_q; \mR)$ be a spherical class, for $n \geq   q+2$.
   Given a co-spherical class $C \in H_k^{\pi}(\widehat{M}; \mZ)$, there is a class
   $\wtgamma_{C} =  \int_{C} \widetilde{\nu}^*(\phi)  \in H^{n-k}(\cS(M) ; \mR)$ and   $[f_{C}] \in \cK_{n-k}(M)$ such that the evaluation
   $\langle \wtgamma_{C}, [f_{C}] \rangle \ne 0$. 
  \end{thm}
\proof
Let $f \colon \mS^n \to B\oG_q$ whose homology class pairs non-trivially with   $\phi \in H^n(B\oG_q; \mR)$.

Let $\xi \colon  \widehat{M} \to \mS^k$ be such that $\xi_*(C)   \in H^k(\mS^k ; \mZ)$ is non-trivial.  Note that   $\xi$ maps the complement of some compact set $K \subset M$ to the basepoint $\theta_k \in \mS^k$. The classifying map $\nu_M \colon M \to B{\rm SO}(q)$ is homotopic to a constant, and the fiber $B\oG_q$ is $(q+1)$-connected, so we can assume that $\nu_M$ is constant on some open neighborhood of $K$.

Let $N = \mS^{n-k} \times M$.
 Let $\wedge \colon \mS^{n-k} \times \mS^k \to \mS^n$ denote the smash product, then let
\begin{equation}\label{eq-fhat}
\widehat{f}_{C}   \colon N = \mS^{n-k} \times M \stackrel{id \times \xi}{\longrightarrow}  \mS^{n-k} \times \mS^k \stackrel{\wedge}{\longrightarrow}  \mS^n \stackrel{f}{\longrightarrow}  B\oG_q \ .
\end{equation}
 Then $\widehat{f}_{C}^*(\phi) \in H_c^n(\mS^{n-k} \times M ; \mR) \cong H^{n-k}(\mS^{n-k} \times M , {x_0} ; \mR)$ is non-vanishing in the cohomology with compact supports,  
 and pairs non-trivially with the cycle defined by 
$(id \times \xi) \colon     \mS^{n-k} \times C \to \mS^{n-k} \times \widehat{M}$. Let $f_C \colon \mS^{n-k} \to \cS(M)$ denote the adjoint to $\widehat{f}_C$.

 Let $\int_{C} \colon H^n(\mS^{n-k} \times M ; \mR) \to  H^{n-k}(\mS^{n-k} ; \mR)$ denote the integration over the cycle $C$,  which is just the slant product with the homology class $C$ followed by projection to $H^{n-k}(\mS^{n-k} ; \mR)$. 
Then as before, we have the    commutative diagram in cohomology:
 
\begin{picture}(300,64)(-30,0)  
\put(24,50){$H_c^n(\mS^{n-k} \times M ; \mR)$}
\put(148,53){\vector(-1,0){40}}
\put(112,60){$(f_{C} \times id)^*$}

\put(72,42){\vector(0,-1){17}}
\put(78,30){$\int_{C}$}
\put(36,10){$H^{n-k}(\mS^{n-k} ; \mR)$}
\put(148,14){\vector(-1,0){40}}
\put(130,20){$f_{C}^*$}
\put(160,50){$H_c^n(\cS(M) \times M ; \mR)$}
\put(288,53){\vector(-1,0){40}}
\put(200,42){\vector(0,-1){17}}
\put(206,30){$\int_{C}$}
 \put(266,60){$\tilde{\nu}^*$}
\put(300,50){$H^n(B\oG_q ; \mR)$}
\put(168,10){$H^{n-k}(\cS(M) ; \mR)$}
\end{picture}

  As $\widehat{f}^*_C = (f_{C} \times id)^* \circ \widetilde{\nu}^*$, we then have   $\ds  f_C^* \circ  \int_{C}    \widetilde{\nu}^*(\phi)  =  \int_{C}   \widehat{f}^*_C \ne 0$.

   Set $\wtgamma_{C} =  \int_{C}  \widetilde{\nu}^*(\phi)  \in H^{n-k}(\cS(M) ; \mR)$ which then pairs nontrivially with  $[f_{C}] \in \pi_{n-k}(\cS(M))$.
 \endproof

 \section{Secondary characteristic classes}\label{sec-secondary}
 
The techniques in the previous section  for constructing classes in $\cK_*(M)$ assume  there is   a spherical class $\phi \in H^n_{\pi}(B\G_q^+; \mR)$. 
The   secondary classes for smooth foliations provides the canonical (and essentially only)  approach to showing the existence of spherical classes for $B\G_q^+$ for degrees  $n > 2q$. We give a brief introduction to their construction.   Lawson    in his survey \cite{Lawson1977} gives a more extensive discussion of the secondary invariants and their properties. 

The normal bundle $Q$ to a smooth-foliation $\F$, when restricted to a leaf $L_x$ of $\F$, has  a natural flat connection $\nabla^{L_x}$ defined by the leafwise parallel transport on $Q$ restricted to $L_x$. An    \emph{adapted connection}  $\nabla^{\F}$ on $Q \to M$ is any connection whose restrictions to leaves equals this natural flat connection. An adapted connection  need not be flat over $M$.  The connection data provided by $\nabla^{\F}$ can be thought of as a  ``linearization'' of the normal structure to $\F$ along the leaves.  Thus, $\nabla^{\F}$  captures aspects of the data provided by the Haefliger groupoid $\GF^r$ of $\F$ --  it is a ``partial linearization'' of the highly nonlinear data which defines the homotopy type of $B\GF$.  In this section, we discuss the applications of this partial linearization  to the study of the space $B\G_q$. 

The seminal observation was made by Bott around 1970. The cohomology ring $H^*(BO(q); \mR) \cong \mR[p_1 , \ldots , p_k]$ where $2k \leq q$, and $p_{j}$ has graded degree $4 j$.

 \begin{thm}[Bott Vanishing \cite{Bott1970}] \label{BVT} Let $\F$ be a codimension-$q$, $C^2$-foliation.   Let  $\nu_{Q}   \colon M \to BO(q)$ be the classifying map for the normal bundle $Q$. Then $\nu_Q^* \colon H^{\ell}(BO(q); \mR) \to H^{\ell}(M; \mR)$ is the trivial map for $\ell > 2q$.
\end{thm}

 Morita observed in \cite[Section~3]{Morita1977} that  there exists a codimension-$2$ smooth foliation  on a closed 4-manifold $M$  for which 
$p_1(Q) = \nu_Q^*(p_1) \in H^4(M; \mR)$ is non-zero.   By taking products of this example, one observes that   the classes $\nu^*(p_1^k) \in H^{4k}(B\G_q; \mR)$  are non-vanishing for $4k \leq 2q$. Moskowitz gave in \cite{Moskowitz1985} an alternate construction  of a compact 4-manifold with a trivial plane field that has $p_1(Q) \ne 0$ for its normal bundle. 
The existence of the codimension-2 foliation in the examples of both  Morita and Moskowitz uses  \cite[Theorem~1]{Thurston1974b} of Thurston to construct the foliation on the manifold. 
A simplified version of Thurston's proof for the case of 2-dimensional foliations is given by Mitsumatsu and Vogt in their work \cite[Theorem~1.1]{MitsumatsuVogt2017}, which is more ``elementary'' in its methods.  In any case,  it seems there is   no explicit construction    of a $C^2$-foliation for which 
$\nu_Q^* \colon H^{\ell}(BO(q); \mR) \to H^{\ell}(M; \mR)$ is non-trivial in the range   $q < \ell \leq 2q$.

It is a remarkable observation that Theorem~\ref{BVT} is false for integral coefficients:
\begin{thm}[Bott-Heitsch \cite{BottHeitsch1972}] \label{BH} The universal normal bundle map, 
\begin{equation}
\nu^* \colon H^{\ell}(BO(q); \mZ) \to H^{\ell}(B\G_q; \mZ)
\end{equation}
 is injective for all $\ell \geq 0$.
\end{thm}
Theorem~\ref{BH} implies that for $q \geq 2$ and $r \geq 2$, the space $B\G_q$ does not have the homotopy type of a finite type CW complex. More is true, that for all $\ell > q/2$,   a CW model for $B\G_q$ must have infinitely many cells in     dimension  $4\ell -1$.

We next recall the construction  of the secondary characteristic classes for $C^2$-foliations.  
 The key observation is that the Bott Vanishing Theorem holds at the level of the differential forms representing the characteristic classes of the normal bundle, and not just for their cohomology classes.

Denote by $I(\mathfrak{gl}(q, \mR))$ the graded ring of adjoint-invariant polynomials on the Lie algebra $\mathfrak{gl}(q, \mR)$ of the real general linear group $GL(q, \mR)$.  As a ring, $I(\mathfrak{gl}(q, \mR)) \cong \mR[c_1, c_2, \ldots , c_q]$  is a polynomial algebra on $q$ generators, where the $i^{th}$-Chern polynomial $c_i$ has  graded degree $2i$.  
Let $I(\mathfrak{gl}(q, \mR))^{(q+1)}$ denote the ideal of polynomials of degree greater than $q$, and introduce    the quotient ring,
$$I_q \equiv I(\mathfrak{gl}(q, \mR))_{q} = I(\mathfrak{gl}(q, \mR))/I(\mathfrak{gl}(q, \mR))^{(q+1)} \cong \mR[c_1, c_2, \ldots , c_q]_{2q}$$ 
which is  isomorphic to a   polynomial ring truncated in graded degrees larger than $2q$.  Associate to each generator $c_i$ the closed $2i$-form $c_i(\Omega({\nabla^{\F}})) \in \Omega_{deR}^{2i}(M)$, so that one obtains the Chern homomorphism 
$\Delta^{\F} \colon \mR[c_1, c_2, \ldots , c_q] \to \Omega_{deR}^{ev}(M)$.
 \begin{thm}[Strong Bott Vanishing \cite{BottHaefliger1972}] \label{SBVT} Let $\F$ be a codimension-$q$, $C^2$-foliation, and $\nabla^{\F}$ an adapted connected on the normal bundle $Q$.   Then for any polynomial $c_J \in \mR[c_1, c_2, \ldots , c_q]$ of graded degree $deg(c_J) > 2q$, the Chern form $c_J(\Omega({\nabla^{\F}}))$ is identically zero. Thus, there is an induced map of differential graded algebras:  
 \begin{equation}\label{eq-SBVT}
\Delta^{\F} \colon   \mR[c_1, c_2, \ldots , c_q]_{2q} \longrightarrow \Omega_{deR}^{ev}(M)  \ .
\end{equation}
\end{thm}
 Associated to the DGA map $\Delta^{\F}$ is a map of   minimal models, whose homotopy class is itself an invariant under concordance for foliations, and determines additional invariants for foliations  \cite{Hurder1980,Hurder1981a}.

Now assume that the normal bundle to the foliation $\F$ on $M$ is trivial. The choice of a framing of $Q$, denoted by $s$, induces an isomorphism $Q \cong M \times {\mathbb R}^q$. Let $\nabla^s$ denoted the connection for the trivial bundle, then the curvature forms of $\nabla^s$ vanish identically. It follows that the transgression form $y_i =Tc_i(\nabla^{\F}, \nabla^s) \in \Omega_{deR}^{2i-1}(M)$   of $c_i$   satisfies the equation 
$dy_i = c_i(\Omega({\nabla^{\F}}))$.
Consider the DGA complex 
\begin{equation}\label{eq-defWq}
W_q  = \Lambda(y_1, y_2, \ldots, y_{q}) \otimes \mR[c_1, c_2, \ldots , c_q]_{2q} \ ,
\end{equation}
where the differential is defined by $d(h_i\otimes 1) = 1 \otimes c_i$.

The data $(\F, s,\nabla^{\F})$  determine a map of differential algebras  $ \Delta^{\F,s} \colon W_q \to \Omega^*(M)$.
The induced map in cohomology, $\ds \Delta^{\F,s}  \colon H^*(W_q) \to H^*(M)$,   depends only on the homotopy class of the framing $s$ and the framed concordance class of $\F$, and moreover, this construction is functorial.

\begin{thm}    There is a    well-defined   universal characteristic map
\begin{equation}
\Delta   \colon    H^*(W_q)    \to    H^*(B\overline{\Gamma}_q ; \mR) 
\end{equation}
Given a   codimension-$q$ foliation $\F$ with framing $s$, the classifying map $\ds \Delta^{\F,s} \colon H^*(W_q) \to H^*(M; \mR)$ factors through the universal map:

\begin{picture}(100,45)\label{comm_diag}
\put(170,40){$H^*(B\overline{\Gamma}_q ; \mR) $}
\put(190,30){\vector(0,-1){17}}
\put(130,13){\vector(3,2){30}}
\put(200,20){$h_{\F}^*$}
\put(130,25){$\Delta$}
\put(150,8){$\Delta^{\F, s}$}
\put(100,0){$H^*(W_q)$}
\put(145,2){$\vector(1,0){22}$}
\put(176,0){$H^*(M; \mR)$}
\end{picture}
\end{thm}

The \emph{Vey basis} for $H^{*}(W_q)$ was given in \cite{Godbillon1974}:
\begin{prop}\label{prop-veybasisWq}
The following set of monomials in $W_q$ form a basis for $H^{*}(W_q)$:
\begin{eqnarray} \label{eq-VeybasisWq}
y_I c_J  & \ {\rm such \  that} \ &  I = (i_1, \ldots , i_s) \ {\rm with} \ 1 \leq i_1 < \cdots < i_s \leq q      \\
  &   &   J = (j_1, \ldots , j_{\ell}) \ {\rm with} \  j_1 \leq \cdots \leq j_{\ell} \ , \  j_1 + \cdots j_{\ell} \leq q \nonumber\\
    &   &    i_1 + j_1 + \cdots j_{\ell} \geq q+1 \ , \  i_1 \leq    j_1 \ .  \nonumber
\end{eqnarray}
\end{prop}
    
   Now consider a smooth foliation $\F$ of a manifold $M$, with no assumption that its normal bundle is trivial. In place of the flat connection $\nabla^s$ associated to a framing of $Q$, let $\nabla^g$ be the connection on $Q$ associated to a Riemannian metric on $Q$. Then the closed forms $c_{2i}(\Omega(\nabla^g)) \in \Omega^{4i}_{deR}(M)$ need not vanish, and in fact their cohomology classes define the Pontrjagin forms for the normal bundle.  On the other hand,  the forms $c_{2i+1}(\Omega(\nabla^g)) \in \Omega^{4i+2}_{deR}(M)$ vanish by the skew symmetry of the curvature matrix $\Omega(\nabla^g)$.
  Thus, we can repeat the above construction of a characteristic map for the subcomplex $W_q$ of $W_q$ defined by
    \begin{equation}
WO_q  = \Lambda(y_1, y_3, \ldots, y_{q'}) \otimes \mR[c_1, c_2, \ldots , c_q]_{2q} \ , 
\end{equation}
where $q' \leq q$ is the largest odd integer with $q' \leq q$.

  The data $(\F, \nabla^{g})$  determine a map of differential algebras  $ \Delta^{\F} \colon WO_q \to \Omega^*(M)$.
The induced map in cohomology, $\ds \Delta^{\F}  \colon H^*(WO_q) \to H^*(M)$,   depends only on the   framed concordance class of $\F$,  but not on the choice of the Riemannian metric on $Q$. It follows that there is a   well-defined   universal characteristic map
$\Delta   \colon    H^*(WO_q)    \to    H^*(B\Gamma_q ; \mR)$. 

The \emph{Vey basis} for $H^{*}(WO_q)$ was given in \cite{Godbillon1974}:
\begin{prop}\label{prop-veybasisWOq}
The following set of monomials in $WO_q$ form a basis for $H^{*}(WO_q)$:
\begin{eqnarray}\label{eq-VeybasisWOq}
y_I c_J  & \ {\rm such \  that} \ &  I = (i_1, \ldots , i_s) \ {\rm with} \ 1 \leq i_1 < \cdots < i_s \leq q \ , \ {\rm each} \ i_k \ {\rm is \ odd}  \\
  &   &   J = (j_1, \ldots , j_{\ell}) \ {\rm with} \  j_1 \leq \cdots \leq j_{\ell} \ , \  j_1 + \cdots j_{\ell} \leq q  \nonumber \\
    &   &    i_1 + j_1 + \cdots j_{\ell} \geq q+1 \ , \  i_1 \leq \ {\rm any \ odd} \  j_k \ .   \nonumber
\end{eqnarray}
\end{prop}
 Note that   $H^{2q+1}(WO_q) \cong H^{2q+1}(W_q)$ when $q$ is even.

The   \emph{secondary characteristic classes} of   foliations are  defined as  those in the image of the maps  $\Delta$ for degree at least $2q+1$, in either the framed or unframed cases.  
 There are special subclasses of the secondary invariants corresponding to special properties, which we mention:

  \begin{defn}\label{def-specialcases} Let  $h_I \wedge c_J$ be an element of the Vey basis for $H^*(WO_q)$ or $H^*(W_q)$.
  \begin{itemize}
\item  $h_1 \wedge c_J$ is said to be a \emph{generalized Godbillon-Vey class}.
\item  $h_I \wedge c_J$ is said to be \emph{residual} if $j_1 + \cdots j_{\ell} = q$.
\item $h_I \wedge c_J$ is said to be \emph{rigid} if $i_1 + j_1 + \cdots j_{\ell} \geq q+2$.
\end{itemize}
  \end{defn}

    The study of the image of the universal maps $\Delta$   have been the primary source of information about the (non-trivial) homotopy type  of   $B\Gamma_q$ and $B\overline{\Gamma}_q$. The outstanding problem remains:
 \begin{prob} \label{conj-injective}
Show that   the   map  $\Delta   \colon    H^*(W_q)    \to    H^*(B\overline{\Gamma}_q ; \mR)$  is injective for     $q \geq 1$.
 \end{prob}

 \section{Foliated manifolds with non-trivial secondary classes}\label{sec-NVsecondary}
 
    We    survey the literature on  examples of foliated manifolds which are used to show the non-triviality of some of the universal  classes in $H^*(B\G_q ; \mR)$ and $H^*(B\oG_q ; \mR)$.

 First, we recall the case for $q=1$ where there $H^{3}(WO_1) \cong \mR$ is generated by the class $y_1c_1$, and for all other degrees, the cohomology vanishes. Let $\F$ be an oriented codimension-one foliation of a manifold $M$, then 
    $\Delta^{\F}(h_1) = \frac{1}{2\pi} \eta  \in \Omega_{deR}^1(M)$ is the Reeb class, and so 
    $\Delta^{\F}(h_1  c_1) = \frac{1}{4\pi^2} \eta \wedge d\eta \in \Omega_{deR}^3(M)$ represents the Godbillon-Vey class     \cite{GodbillonVey1971}.

\begin{thm} [Godbillon-Vey \cite{GodbillonVey1971}]   Let $\F$ be a codimension-one, $C^2$ foliation on $M$ with trivial normal bundle. Then the  3-form $\Delta^{\F}(h_1  c_1)$ is closed, and its cohomology class  
 $GV(\F)   \in H^{3}(M; \mR)$ is independent of all choices. Moreover, $GV(\F)$ depends only on the concordance class of $\F$.
 \end{thm}
 
The   paper \cite{GodbillonVey1971} also included an example by Roussarie to show this class was non-zero. Soon afterwards, Thurston gave  a construction in \cite{Thurston1972} of families of examples of foliations on the 3-sphere $\mS^3$, such that the $GV(\F) \in H^3(\mS^3; \mZ) \cong \mR$  assumes a continuous range of real values. As a consequence:

\begin{thm} [Thurston \cite{Thurston1972}]   \label{thm-thurstonq=1}
The Godbillon-Vey class is a surjection   
 $GV \colon   \pi_3(B\Gamma_1) \to \mR$. 
\end{thm}

The appendix by Brooks to   \cite{Bott1978} gives a  detailed explanation of the examples Thurston constructed.   The Seminaire Bourbaki article by Ghys \cite{Ghys1989} is a basic reference for the properties of the Godbillon-Vey class; the author's survey  \cite{Hurder2002a} is a more recent update.

 For   a foliation $\F$  of  codimension $q > 1$ with normal bundle framing $s$, the    secondary classes of $\F$ are spanned by the   images $\ds \Delta^{\F,s}(y_I    c_J)$, where $y_I   c_J$ satisfies \eqref{eq-VeybasisWq}.   Without the assumption that the normal bundle is framed, we   restrict to monomials $y_I   c_J$ as in \eqref{eq-VeybasisWOq}.

The original example of Roussarie \cite{GodbillonVey1971}, and its extensions to   codimension $q > 1$, start with a semi-simple Lie group $G$.  Choose   closed subgroups $K \subset H \subset G$, with $K$ compact.  Then $G/K$ is foliated by the left cosets of $H/K$.
Choose a cocompact, torsion-free lattice $\G \subset G$, then the  foliation of $G/K$ descends to   a foliation $\F$ on the compact manifold  $M = \G\backslash G/K$ which is a    locally homogeneous space. The calculation of the secondary invariants for such foliations  then follows from explicit calculations in Lie algebra cohomology, using    Cartan's approach to the cohomology of homogeneous spaces.  Examples of this type are studied in 
\cite{Baker1978,Fuks1976,KT1974d, KT1975a, KT1978a,KT1979, Pelletier1983,Pittie1976a, Pittie1979, Tabachnikov1984, Tabachnikov1985, Tabachnikov1986, Yamato1975}. 
 For example,  Baker  shows in  \cite{Baker1978}:
 \begin{thm} \label{thm-injectivityframed}
 Let $q = 2m > 4$. Then the set of residual classes 
$$\{h_1h_2h_{i_1}\cdots h_{i_k}   c_1^{q} \ ; \  h_1h_2h_{i_1}\cdots h_{i_k}   c_1^{q-2} c_2   \mid 2<i_1<\cdots<i_k\leq m\}$$
in $H^*(W_{q})$ map under $\Delta$ to linearly independent classes in $H^*(B\overline{\Gamma}_q; \mR)$.
 \end{thm}
 The non-vanishing results of Kamber and Tondeur follow a similar format, but are more extensive, as given in Theorem~7.95 of \cite{KT1975a} for example. The conclusion of all these approaches is to show that the    universal map $\Delta$ is  injective on various subspaces of  $H^*(W_q)$ and $H^*(WO_q)$.
   
  For the case of foliations with codimension $q =2$,   Rasmussen \cite{Rasmussen1980} modified  the construction by Thurston of codimension-one foliations with varying Godbillon-Vey class. Thurston's construction used the weak-stable foliation of the geodesic flow on a compact Riemann surface (with boundary) of constant negative curvature. Rasmussen extended these ideas to the case of   compact hyperbolic 3-manifolds.  He showed there exists families $\{\F_{\lambda} \mid \lambda \in \mR\}$  of smooth foliations in codimension-$2$ on oriented compact $5$-manifolds $M_1$ and $M_2$ for which the secondary    classes $\{ h_1   c_1^2, h_1   c_2 \}$ vary continuously and independently when evaluated on the fundamental classes of $M_1$ and $M_2$.

 For the case of foliations with codimension $q \geq 3$,   Heitsch developed a residue theory for the residual secondary classes of smooth foliations \cite{Heitsch1977, Heitsch1978}, and used this to construct extensive families of foliations whose secondary classes varied continuously \cite[Theorems~6.1, 6.2, 6.3]{Heitsch1978}.

  Together, these examples yield that for a fixed collection of classes in the image of $\Delta$, there are continuous families of cycles $h_{\F_{\lambda}} \colon M \to B\G^+_q$ or $h_{\F_{\lambda},s} \colon M \to B\overline{\Gamma}_q$ such that the evaluation of these fixed secondary classes on the cycles defined by the $\F_{\lambda}$ vary continuously.  
 Thus, the groups $H_*(B\Gamma^+_q ; \mZ)$ and $H_*(B\overline{\Gamma}_q ; \mZ)$ must be a truly enormous.

  There are two very remarkable results of a different nature, one  for  foliations of codimension $1$ by Tsuboi,   and the other   for foliations of codimension $2$ by Boullay, on the existence of divisible subgroups in the homology of $B\Gamma_q^+$. 
  
  \begin{thm}  \cite[Theorem~1.2]{Tsuboi1984}   \label{prob-divisable1}
   $H_3(B\G_1^+ ; \mZ)$  contains a subgroup isomorphic to $\mR$.
 \end{thm}

 \begin{thm} \cite[Corollary~4]{Boullay1996}  \label{prob-divisable2}
   $H_5(B\G_2^+ ; \mZ)$  contains a subgroup isomorphic to $\mR^2$.
 \end{thm}
The work of Boullay uses the results on divisibility in algebraic K-groups to construct his examples.  Nariman discusses this connection in   \cite{Nariman2022}.

As mentioned above, Morita constructed in     \cite{Morita1977} a  compact 4-manifold $M$ with a  codimension 2 foliation  that has non-trivial first Pontrjagin class $p_1 \in H^4(M; \mR)$.  By taking products  of the examples constructed by Tsuboi or Boullay, depending on whether the codimension $q$ is odd of even, with an appropriate number of copies of this manifold, one obtains cycles in arbitrary codimension $q$ that produce divisible subgroups in the homology of degree $2q+1$.  Using the same calculations as in    \cite{Morita1977}, one obtains:
 \begin{cor}
 For all $q \geq 1$,  $H_{2q+1}(B\G_q^+ ; \mZ)$  contains a subgroup isomorphic to $\mR$.
 \end{cor}
   Actually, much more can be shown using    homotopy techniques, as discussed in \cite{Hurder2023b}.
   
 The notion of \emph{discontinuous invariants} for foliations was developed by Morita in his works \cite{Morita1984a,Morita1985} (and see also \cite[Chapter~3]{Morita2001})   give an approach for the   study of the variable foliation invariants, and is also analogous  to the divisibility property for the integer homology groups. The universal Godbillon-Vey class $GV \in H^3(B\G_1^+ ; \mR)$ defines a map $\overline{GV} \colon B\G_1^+ \to K(\mR^{\delta}, 3)$ where $\mR^{\delta}$ denotes the group of reals with the discrete topology, and $K(\mR^{\delta},3)$ is the Eilenberg-MacLane space in dimension 3 for $\mR^{\delta}$.
 While the homotopy groups of $K(\mR^{\delta}, 3)$ are prescribed, the cohomology of this space is absolutely enormous. The image of 
 $\overline{GV}^* \colon H^*(K(\mR^{\delta},3) ; \mZ) \to H^*(B\G_1^+ ; \mZ)$ are the \emph{discontinuous invariants} for codimension one foliations. 
 \begin{prob}[Morita \cite{Morita1985}] \label{prob-Morita}
 Determine the image of $\overline{GV}^*$ in $H^*(B\G_1^+ ; \mZ)$. 
  \end{prob}
  There are an analogous constructions of discontinuous invariants in each codimension $q > 1$ as well, using  the variable secondary classes to define maps to products of Eilenberg-MacLane spaces. The non-triviality of products of the discontinuous invariants is completely unknown for smooth foliations, though the work of Tsuboi in \cite{Tsuboi1998} shows there are non-trivial products for foliations which are transversally Lipshitz.

    What is striking, looking back at the roughly 15 years between 1971 and 1986 during which the study of secondary invariants for foliations  was most actively researched, is   the many approaches and expositions of the subject, especially  for example   
\cite{Baker1978,BernshteinRozenfeld1972,BernshteinRozenfeld1973,Bott1972a,Bott1975c,Bott1978,BottHaefliger1972,Fuks1973,Fuks1976,GF1968,GF1969b,KT1974c,KT1974c,KT1975a,KT1978a,Lawson1977,Morita2001,Pittie1976a,Pittie1979}.
The  notation  in these papers is not  consistent.

It is also striking    how limited are the types of examples presented in these works.  
 All of the ``explicit constructions'' of foliations in the literature with non-trivial secondary characteristic classes are either locally homogeneous, or deformations of locally homogeneous actions. Essentially, they  are all   generalizations and/or modifications of the original constructions of Roussarie \cite{GodbillonVey1971} and Thurston  \cite{Thurston1972} for the Godbillon-Vey class in codimension-one. Nothing is known about \emph{explicit} constructions of examples with non-vanishing secondary classes beyond these ``almost homogeneous'' examples.

 \section{Spherical secondary classes}\label{sec-spherical}

In this section, we give a review the known results about  the spherical secondary classes for the spaces $B\G_q^+$ and $B\oG_q$. 
  A key tool is the Rational Hurewicz  Theorem in \ref{subsec-hurewicz}.   We  also use  the \emph{brace product} in Section~\ref{subsec-brace}   to extend the non-triviality results for $B\overline{\Gamma}_q$ from degree $2q+1$ to    higher degrees.
  First, we give the prototypical construction of spherical cohomology classes.

\begin{thm}\label{thm-variablehomotopy}
Let $\phi \in H^{2q+1}(B\overline{\Gamma}_{q} ; \mR)$ be a cohomology class such that there exists a continuous family of framed of codimension-q foliations $\{(\F_{\lambda}, s_{\lambda}) \mid \lambda \in \Lambda\}$    on a     manifold $Y$, such that the pull-back $\chi_{\F_{\lambda}}^*(\phi) \in H^{2q+1}(Y ; \mR)$ is non-trivial, and varies continuously with $\lambda$ if $\Lambda$ is not an isolated set. Then for each $\lambda \in \Lambda$,  there exists $[f_{\lambda}] \in \pi_{2q+1}(B\overline{\Gamma}_{q})$  so that the pairing 
$\langle \phi , [f_{\lambda}] \rangle$ is non-vanishing and varies continuously with $\lambda$. That is, $\phi$ induces a map
 $\phi  \colon \pi_{2q+1}(B\overline{\Gamma}_{q}) \to \mR$ which is surjective.
\end{thm}
\proof
First note that as the evaluation map $\phi \colon H_{2q+1}(B\overline{\Gamma}_{q} ; \mZ) \to \mR$ is a homomorphism, if it is onto an open set, then it is onto all of $\mR$.

Now let   $\{F_{\lambda} \mid \lambda \in \Lambda\}$ be a continuous   family of codimension-q foliations on $Y$. Suppose that the normal bundles to the foliations admit  trivializations $s_{\lambda}$ that depend continuously on $\lambda$. Then there exists a continuous family of classifying maps for these foliations, denoted by
 $\chi_{\F_{\lambda}} \colon Y \to B\overline{\Gamma}_q$, where we suppress the dependence on the framing $s_{\lambda}$. Let $\phi \in H^{2q+1}(B\overline{\Gamma}_q ; \mR)$ be a cohomology class  such that the pull-back $\chi_{\F_{\lambda}}^*(\phi) \in H^{\ell}(Y; \mR)$ is non-trivial for $\lambda \in \Lambda$.
 
 Since $Y$ has the homotopy type of a   CW complex, there is a   CW complex $Y_{(q+1)}$ with all cells in dimension at least $q+2$, and a   map $\xi \colon Y \to Y_{(q+1)}$ such that the map $H_{2q+1}(Y ; \mZ) \to H_{2q+1}(Y_{(q+1)} ; \mZ)$ is injective. The space $Y_{(q+1)}$ is obtained by collapsing the cells of dimension less than $q+2$ in a CW decomposition of  $Y$. As a consequence of Theorem~\ref{thm-MT} and basic obstruction theory, the classifying maps $\chi_{\F_{\lambda}}$ descend to a   family of maps 
 $\overline{\chi}_{\F_{\lambda}} \colon Y_{(q+1)} \to B\overline{\Gamma}_q$ such that $\overline{\chi}_{\F_{\lambda}}^*(\phi) \in H^{2q+1}(Y_{(q+1)}; \mR)$ is non-trivial for $\lambda \in \Lambda$.

By the Rational Hurewicz Theorem~\ref{thm-RHIT} as applied to  the space $Y_{(q+1)}$  for $r = q+1$, the map
  $$h \colon \pi_{2q+1}(Y_{(q+1)}) \otimes \mQ  \longrightarrow   H_{2q+1}(Y_{(q+1)} ; \mQ) \cong   H_{2q+1}(Y ; \mQ)$$
   is an isomorphism. The class $\overline{\chi}_{\F_{\lambda}}^*(\phi) \in H^{2q+1}(Y_{(q+1)}; \mR)$ is non-vanishing, so there exists a cycle $[Z] \in H_{2q+1}(Y_{(q+1)}, \mZ)$ such that the pairing $\langle \overline{\chi}_{\F_{\lambda}}^*(\phi), [Z] \rangle \ne 0$. Moreover, if $\Lambda$ does not consists of isolated points, the evaluation pairing varies continuously with the parameter $\lambda \in \Lambda$, so 
 $\langle \overline{\chi}_{\F_{\lambda'}}^*(\phi), [Z] \rangle \ne 0$ for $\lambda' \in \Lambda$ sufficiently close to $\lambda$.
Choose a map $f_{[Z]} \colon \mS^{2q+1} \to Y_{(q+1)}$ such that the image $h([f_{[Z]}]) \in H_{2q+1}(Y_{(q+1)} ; \mZ)$ is a non-zero rational multiple of $[Z]$. Then the composition 
\begin{equation}
\overline{\chi}_{\F_{\lambda}} \circ f_{[Z]} \colon \mS^{2q+1} \to B\overline{\Gamma}_q
\end{equation}
defines a family of spherical (2q+1)-cycles in $B\overline{\Gamma}_q$ which pair      non-trivially with the class $\phi$.  
 \endproof

Next, consider the homotopy properties for the space $B\overline{\Gamma}_q$  derived from the fibration sequences   \eqref{eq-commutingdiagram}. The  first observation is that $B\overline{\Gamma}_q$ is $(q+1)$-connected implies that the map $\nu \colon B\Gamma_q^+  \to  B{\rm SO}(q) $ admits a section over the $(q+2)$-skeleton   $B^{(q+2)} \subset B{\rm SO}(q)$ for a CW decomposition of $B{\rm SO}(q)$.  For $q=4\ell -1$ or $q=4\ell -2$, there is a class $\xi_{\ell} \colon \mS^{4\ell} \to B{\rm SO}(q)$ such that $\xi_{\ell}^*(p_{\ell}) \in H^{\ell}(\mS^{4\ell}, \mZ)$ is non-vanishing, and $\xi_{\ell}$ lifts to a map $\widehat{\xi}_{\ell} \to B\Gamma_q^+$ such that $\widehat{\xi}_{\ell}^* \circ \nu^*(p_{\ell}) \ne 0$, hence $\nu^*(p_{\ell}) \ne 0$. 

By the Bott Vanishing Theorem,  $\nu^*(p_{\ell}^2) = 0$, so the map $\widehat{\xi}_{\ell} \vee \widehat{\xi}_{\ell} \colon \mS^{4\ell} \vee \mS^{4\ell} \to B\Gamma_q^+$ cannot be extended to the product $\mS^{4\ell} \times \mS^{4\ell}$. The obstruction to this extension is the Whitehead product $[\widehat{\xi}_{\ell}, \widehat{\xi}_{\ell}] \in \pi_{8\ell -1}(B\Gamma_q^+)$ which must be non-vanishing in rational homotopy theory. Schweitzer and Whitman observed   that   this homotopy class is non-zero in the note  \cite{SchweitzerWhitman1978}, using a theory of residues for the Pontrjagin classes.

Haefliger observed in the note \cite{Haefliger1978b} the following more general formulation. Suppose there is given a collection of forms $\{\phi_i \mid i \in \cI\} \subset \Omega^*_{deR}(B\Gamma_q^+)$ in the de~Rham complex of $B\Gamma_q^+$ which  satisfy $\phi_i \wedge \phi_j = 0$ for $i \ne j$. Moreover, assume  the linear functionals they define   $[\phi_i] \colon \pi_*(B\Gamma_q^+) \to \mR$ are rationally independent. Choose a collection of maps $\{f_i \colon \mS^{n_i} \to B\Gamma_q^+ \mid i \in \cI, n_i = \deg(\phi^i)\}$ for which their classes are independent when restricted to the homology classes they determine. Then for the induced map on their wedge product,
\begin{equation}
\widehat{f} \colon X = \bigvee_{i \in \cI} \ \mS^{n_i} \to B\Gamma_q^+ \ ,
\end{equation}
the induced map in rational homotopy theory,  $\widehat{f}_{\#} \colon \pi_*(X) \otimes \mQ \to \pi_*(B\Gamma_q^+) \otimes \mQ$, is injective. The rational homotopy groups $\pi_*(X) \otimes \mQ$ have the structure of a free graded Lie algebra, with multiplication given by the Whitehead product, and this space is infinite dimensional if the index set $\cI$ is not a singleton. 
Note that the  rational Hurewicz homomorphism $h_* \colon \pi_*(B\G_q^+) \to H_*(B\G_q^+ ; \mQ)$ maps all of these Whitehead product classes to zero.

Both of the observations, by Schweitzer and Whitman and by Haefliger, are subsumed by the  author's thesis work \cite{Hurder1980,Hurder1981a} using Sullivan's approach to rational homotopy theory using differential forms to construct higher order invariants, the dual homotopy classes.  There are many expositions of rational homotopy theory, but the original papers  \cite{DGMS1975,GM1981,Sullivan1976,Sullivan1978} explain what is used below. 
We briefly recall some aspects of this approach. 

Recall that a model for a differential graded algebra $\cA$ (over $\mQ$ or $\mR$) is a DGA map $\psi \colon \fM \to \cA$
where $\fM$ is a free DGA and the map $\psi$ induces an isomorphism on cohomology groups.  Then $\fM$ is a minimal model for $\cA$ if it is the minimal such DGA.  Minimal models exist under some basic restrictions.  Let $\fM(\cA) \to \cA$ denote a minimal model for $\cA$. 

For a semi-simplicial space $X$, let $\Omega^*_{deR}(X)$ denote the DGA of semi-simplicial forms on $X$, and let $\fM(X)$ denote a minimal model for $\Omega^*_{deR}(X)$. Let $\pi^*(\fM(X))$ denote the space of indecomposable elements in $\fM(X)$.
A basic result of minimal model theory is that there is a natural isomorphism $\pi^*(\fM(X))  \cong {\rm Hom}(\pi_*(X) , \mR)$ so that the space $\pi^*(\fM(X))$ consists of ``dual homotopy'' classes.

The map \eqref{eq-SBVT} defined by the Strong Bott Vanishing Theorem~\ref{SBVT} is natural, so induces a DGA map 
  \begin{equation}
\Delta \colon I_q  \longrightarrow \Omega_{deR}^{*}(B\G^+_q)  \ ,
\end{equation}
where recall that we defined $I_q \equiv   \mR[c_1, c_2, \ldots , c_q]_{2q}$ to be the truncated polynomial algebra. This determines a map on minimal models $\fM(\Delta) \colon \fM(I_q) \to \fM(B\G^+_q)$ which induces a map on dual homotopy, denoted by 
 $\Delta^{\#} \colon \pi^*(I_q) \to {\rm Hom}(\pi_*(B\G^+_q) , \mR)$.
The classes in the image of $\Delta^{\#}$ are called the \emph{dual homotopy invariants} of foliations in \cite{Hurder1981a}. For example,  corresponding to the indecomposable element of $\fM(I_q)$ defined by the Pontrjagin class $p_{\ell} = c_{2\ell} \in I_q$, there is an indecomposable element $y_{8\ell-1} \in \fM(I_q)$, and the image $\Delta^{\#}(y_{8\ell-1}) \in {\rm Hom}(\pi_{8\ell-1}(B\G^+_q) , \mR)$ detects the Whitehead product constructed by Schweitzer and Whitman. Moreover, it was shown that the class $\Delta^{\#}(y_{8\ell-1})$ is equal to the evaluation of the secondary class $\Delta(y_{2\ell} c_{2\ell}) \in H^{8\ell-1}(B\overline{\Gamma}_q ; \mR)$. Similar results hold for other elements of $\pi^*(I_q)$ under the image of $\Delta^{\#}$, so the study of the dual homotopy invariants for foliations gives a useful method of constructing non-trivial spherical cohomology classes in 
$H^{*}(B\overline{\Gamma}_q ; \mR)$. Many more such examples are given in the works \cite{Hurder1980,Hurder1981a,Hurder1981b,Hurder1985a}. It should be noted that the complete structure of the model $\fM(I_q)$ is unknown, as its calculation becomes exponentially more complicated in higher degrees as the codimension   $q$ increases. It is known that $\pi^*(I_q)$ contains graded Lie algebras of increasingly large size, which are the analogues of the examples studied by Haefliger in \cite{Haefliger1978b}.

Next, we construct additional spherical classes in $H^*(B\overline{\Gamma}_q ; \mR)$ using  the \emph{brace product} in Section~\ref{subsec-brace}.   
Recall that   $\iota \colon B^{(q+2)} \subset B{\rm SO}(q)$ denotes the    $(q+2)$-skeleton for a CW decomposition of $B{\rm SO}(q)$, and 
let 
$\sigma_q \colon  B^{(q+2)} \to B\Gamma_q^+$ be a section to the fibration $\nu \colon B\Gamma_q^+ \to  B{\rm SO}(q)$.
Introduce the pull-back $E_q$ of this fibration over the inclusion  $B_{(q+2)} \subset B{\rm SO}(q)$, so we obtain a fibration
 \begin{equation}\label{eq-truncatedfibration}
B\overline{\Gamma}_q \stackrel{\iota}{\longrightarrow} E_q \stackrel{\nu}{\longrightarrow}  B^{(q+2)}
\end{equation}
 which admits a section   denoted by $\sigma  \colon B^{(q+2)} \to E_q$. 
 Given classes $\alpha \in \pi_i(B^{(q+2)})$ and $\beta \in \pi_j(B\overline{\Gamma}_q)$ the brace product 
 $\{ \alpha  ,  \beta \} \in \pi_{i+j-1}(B\oG_q)$. We apply this construction to spherical secondary classes, to construct additional spherical secondary classes of higher degree.

 Let $0 < 4\ell \leq q+1$ then there exists a map $\tau_{\ell} \colon \mS^{4\ell} \to B^{(q+2)}$ such that the Pontrjagin class $p_{\ell} \in H^{4\ell}(B{\rm SO}(q) ; \mZ)$ pairs non-trivially with the image of the composition $\iota \circ \tau_{\ell} \to B{\rm SO}(q)$. Then the composition $\sigma \circ \tau_{\ell} \colon \mS^{4\ell} \to E_q$ is a lift of this map, and the brace product defines a mapping $\{\tau_{\ell}, \cdot \} \colon \pi_j(B\overline{\Gamma}_q) \to \pi_{4\ell + j -1}(B\overline{\Gamma}_q)$.   We use this mapping to show:

 \begin{thm}\label{thm-bracing}
 Assume $ q \geq 3$. 
 Let $y_{i_1} c_J \in WO_q$ be an element of the Vey basis of degree $2q+1$, and suppose there exists a map $f \colon \mS^{2q+1} \to B\oG_q$ such that $\Delta(y_{i_1} c_j)$ pairs non-trivially on the homology class $[f] \in H_{2q+1}(B\G_q^+ ; \mR)$.
 That is, $\Delta(y_{i_1} c_J) \in H_{\pi}^{2q+1}(B\G_q^+ ; \mR)$. Assume that $I' = (i_1', \ldots , i_r')$ satisfies
 \begin{equation}\label{eq-bracedclasses}
 i_1 <  i_1'  <  i_2' <  \cdots <  i_r'   \ , \ 2 i_r' \leq q+1 \ , \ {\rm and } \  i_{\ell}' \ {\rm is \ even \ for} \ 1 \leq \ell \leq r  \ .
\end{equation}
 Then the secondary class $\Delta(y_I y_{I'} c_J) \in H_{\pi}^{2q+k}(B\oG_q ; \mR)$.    That is, it is also non-trivial and spherical.
 \end{thm}
 \proof
 We construct a mapping $g \colon \mS^k \to B\oG_q$ such that the pairing $\langle \Delta(y_{i_1} y_{I'} c_J) , [g] \rangle  \ne 0$.
 We have  $2i_1' \leq q+1$, so there exists a lifting $\whtau_{2 i_1'} \colon \mS^{2 i_1'} \to B\G_q^+$  that pairs non-trivially with $\nu^*(c_{2 i_1'}) \in H^{2i_1'}(B\G_q^+ ; \mR)$.
 For  the   map $f \colon \mS^{2q+1} \to B\G_q^+$, form the brace product   $\{[\whtau_{2 i_1'}], [f]\}_{\sigma} \in \pi_{2q+2i_1'}(B\oG_q)$.
 
 The class $[\whtau_{2 i_1'}]$ pairs non-trivially with $\nu^*(c_{2 i_1'})$ and $[f]$ pairs non-trivially with $\Delta(y_{i_1} c_J)$, so if the brace product vanishes in $\pi_{2q+2i_1'}(B\oG_q) \otimes \mQ$, then the product 
 $\Delta(y_{i_1} c_{2i_1'} c_J) \in H^{2q+ 2 i_1' +1}(B\G_q^+ ; \mR)$ must be non-vanishing. But the  product $c_{2i_1'} c_J = 0$ in $WO_q$ so this is a contradiction.
 
 Now proceed recursively on the index set $I' = (i_1', \ldots , i_r')$. Set $f_0 = f \colon \mS^{2q+1} \to B\oG_q$ and 
 $f_1 \colon \mS^{2q+2i_1'} \to B\oG_q$ which represents the brace product $\{[\whtau_{2 i_1'}], [f]\}_{\sigma}$.
 Now assume for $\ell < r$ that maps $f_{\ell} \colon \mS^{2q+k} \to B\oG_q$ have been constructed which pair non-trivially with the class
 $\Delta(y_{i_1} y_{i_1'} \cdots y_{i_{\ell}'} c_J) \in H^{2q+k}(B\oG_q ; \mR)$. 
 Let $f_{\ell+1} \colon \mS^{2q+k+2i_{\ell+1}' -1} \to B\oG_q$ represent the brace product 
  $\{[\whtau_{2 i_{\ell +1}'}], [f_{\ell}]\}_{\sigma}$. 
  We   show this class pairs non-trivially with the cohomology class   
\begin{equation}\label{eq-indectivesc}
\Delta(y_{i_1} y_{i_1'} \cdots y_{i_{\ell+1}'} c_J) \in H^{2q+k+2i_{\ell+1}' -1}(B\oG_q ; \mR) \ . 
\end{equation}
Use the map $\whtau_{2 i_{\ell+1}'}$ to form a pull-back fibration
$\ds  B\oG_q \to E_{\ell+1} \to \mS^{2i'_{\ell +1}}$. 

The space $E_{\ell+1}$ is $(2i'_{\ell +1} - 1)$-connected, so the Pontrjagin classes $\Delta(c_{i_{j}'})$ vanish when restricted to $E_{\ell+1}$, for    $1 \leq j \leq \ell$. Thus the class $\Delta(y_{i_1} y_{i_1'} \cdots y_{i_{\ell}'} c_J)$ is well-defined on $E_{\ell+1}$. It follows that the brace product $\{[\whtau_{2 i_{\ell +1}'}], [f_{\ell}]\}_{\sigma}$ is non-vanishing in 
$\pi_{2q+k+2i_{\ell+1}' -1}(B\oG_q)$ by the same argument as before, and pairs non-trivially with the secondary class 
in \eqref{eq-indectivesc}.

We thus obtain the map $g = f_r \colon \mS^{2q+k} \to B\oG_q$ which pairs non-trivially with $\Delta(y_I y_{I'} c_J)$.
  \endproof
 
 \begin{remark}
 {\rm
 There is a more formal proof of Theorem~\ref{thm-bracing} using the calculation of the minimal model for the space $E_q$ in \eqref{eq-truncatedfibration}. The steps are the same, using what amounts to an iterated Hirsch construction (that is,  a Postnikov tower) as is typical for building minimal models. The above description has a more geometric flavor. 
 }
 \end{remark}
 
We conclude by recalling a   theorem  from \cite{Hurder1985a},  which summarizes some of the  known results on the  non-triviality of the higher degree homotopy groups $\pi_k(B\overline{\Gamma}_{q})$ that follows using the dual homotopy invariants. The conclusion is that the homotopy groups of $B\oG_q$ are non-trivial in an infinite range      and uncountably generated,  as the degree $\ell$ tends to infinity. 
 The integers  $\{v_{q,\ell}\}$  in the statement below are defined in Section~2.8 of \cite{Hurder1985a}, and are defined as the ranks of various free graded Lie algebras in the minimal model for the truncated polynomial algebra $ \mR[c_1, c_2, \ldots , c_q]_{2q}$.

\begin{thm}  \cite[Theorem~1]{Hurder1985a} \label{thm-lotsaclasses}
For  $q > 1$ and $\ell > 2q$, there is an epimorphism of abelian groups 
\begin{equation}
h^* \colon \pi_{\ell}(B\overline{\Gamma}_{2q+1}) \to \mR^{v_{q, \ell}} \ ,
\end{equation}
where the sequences of non-negative integers $\{v_{q,\ell}\}$  satisfy:
 \begin{enumerate}
\item For $q=2$, $\ds \lim_{k\to \infty} v_{2, 4k+1} = \infty$, with $v_{2, 4k+1} > 0$ for all $k > 0$;
\item For $q=3$, $\ds \lim_{k\to \infty} v_{3, 3k+1} = \infty$, with $v_{3, 3k+1} > 0$ for all $k > 0$;
\item For $q>3$, $\ds \lim_{k \to \infty} v_{q, k} = \infty$.
\end{enumerate}
 \end{thm}

This has as an immediate consequence.
 \begin{cor}\label{thm-lotsahomology}
 For each $q \geq 1$, the homology groups $H_*(B\overline{\Diff_c(\mR^q)} ; \mQ)$ contains a free exterior algebra whose basis elements  in degree $\ell$ correspond to a Hamel basis for $\mR^{v_{q,\ell}}$.
 \end{cor}

\section{Flat Diff(M)-bundles with non-trivial characteristic classes}\label{sec-flatbundles}

In this section, we   describe the explicit construction of non-trivial classes in $H_*(B\overline{\Diff_c(M)} ; \mQ)$ 
starting with the cycles in $H_*(\cS(M); \mQ)$  constructed in the previous sections. We then consider 
  the problem of determining the image  of these classes  in $ H_*(B\Diff_c^+(M)^{\delta} ; \mQ)$.

   \begin{prop}\label{prop-homology}
   Let $M$ be a compact, connected oriented manifold of dimension $q$.
  Given  $\phi \in H^n_{\pi}(B\G_q^+; \mR)$  and a map 
   $f \colon \mS^n \to B\oG_q$ whose homology class pairs non-trivially with   $\phi$, and also  assume  there is given  a co-spherical   class       $C \in H_k^{\pi}(\widehat{M}; \mZ)$. Then there exists $\psi_C \in H^{n-k}(B\overline{\Diff_c(M)} ; \mR)$, 
   a compact oriented manifold $Z_C$ and map $g_C \colon Z_C \to  B\overline{\Diff_c(M)}$, and   rational number $r \ne 0$ such that 
   \begin{equation}\label{eq-pairing}
\langle \psi_C , [g_C(Z_C)] \rangle \ = r \cdot \langle \phi , [f(\mS^n)] \rangle \ ,
\end{equation}
   \end{prop}
   \proof
      Recall that $\widetilde{\nu}_M   \colon M\Diff_c^+(M)^{\delta}   \to B\G_q^+$ classifies the foliation on the total space over  $B\Diff_c^+(M)^{\delta}$, and by restriction also gives a map $\widetilde{\nu}_M   \colon B\overline{\Diff_c(M)} \times M   \to B\G_q^+$.

The proof of Theorem~\ref{thm-method3} constructs a  class     $\wtgamma_{C}    \in H_{\pi}^{n-k}(\cS(M) ; \mR)$  and     defines as  in \eqref{eq-fhat} a map
$\widehat{f}_{C}   \colon N = \mS^{n-k} \times M \to  B\oG_q$ such that   $\widehat{f}_{C}^*(\phi) \in H^n(\mS^{n-k} \times M ; \mR)$ is non-vanishing.

Let $f_C \colon \mS^{n-k} \to \cS(M)$ be the adjoint to $\widehat{f}_C$, then   $\langle \wtgamma_{C}, [f_C] \rangle \ne 0$.

The Mather-Thurston Theorem~\ref{thm-MTHT} says that for the adjoint map $\tau \colon B\overline{\Diff_c(M)} \to \cS(M)$, the   induced map $\tau_* \colon H_*(B\overline{\Diff_c(M)} ; \mZ) \to H_*(\cS(M) ; \mZ)$ is an isomorphism. We make this isomorphism explicit for the spherical class determined by $f_C \colon \mS^{n-k} \to \cS(M)$.

The  map $\widehat{f}_C \colon N = \mS^{n-k} \times M  \to B\G_q^+$   defines 
 a Haefliger structure on $\mS^{n-k} \times M$ whose normal bundle is the lift of the tangent bundle to $M$.
By  the Mather-Thurston Theorem~\ref{thm-MTHT}, there is 
\begin{itemize}
\item a cycle $g_C \colon Z_C \to B\overline{\Diff_c(M)}$; 
\item an oriented   simplicial complex $W_C$ of dimension $(n-k+1)$; 
\item a map $b_C \colon W_C \to \cS(M)$ with boundary components $\partial W_C = \partial W_C^+ \cup \partial W_C^-$, such that $b_C(W_C^-) = f(\mS^{n-k})$ and $b_C(W_C^+) = g_C(Z_C)$
\end{itemize}
 That is, $b_C(W_C)$ is a boundary between the cycles $\tau \circ g_C(Z_C)$ and $f_C(\mS^{n-k})$ in $\cS(M)$.
 
 \begin{remark}\label{rmk-complexity}
 {\rm 
 Since $H_*(B\oG_q ; \mQ)$ and the oriented bordism group $\Omega_*^{SO}(B\oG_q) \otimes \mQ$  are isomorphic,  we can assume that both $Z_C$ and $W_C$ are oriented compact manifolds if desired,   replacing $f_C$ with a positive multiple of itself, if necessary.  For $q \geq 2$, the proof of Mather-Thurston by Meigniez \cite{Meigniez2021} then constructs the simplicial complex  $W_C$ via a sequence of foliated surgeries on $\mS^{n-k} \times M$. The cycle  $g_C \colon Z_C \to B\overline{\Diff_c(M)}$ then arises as the boundary of the   space $W_C$ so constructed. 
 The other proofs of the Mather-Thurston Theorem in \cite{Mather1975,Mather1979,McDuff1979a,Thurston1974} prove that  $\tau_*$ is an isomorphism by more abstract methods, so reveal even less information about the space $Z_C$ and the action of its fundamental group on $M$.
}
\end{remark}

Let $\widehat{b}_C \colon W_C \times M \to B\G_q^+$ denote the adjoint of $b_C$, and set   $\widehat{g}_C = \widetilde{\nu}_M \circ (g_C \times id) \colon Z_C \times M \to B\G_q^+$. 
Then   $\widehat{b}_C(W_C \times M)$ is a foliated concordance $\widetilde{\nu}_{W_C} \colon W_C \times M \to B\G_q^+$ between $\widehat{f}_C(\mS^{k} \times M)$ and $\widehat{g}_C (Z_f \times M)$ in $B\G_q^+$.
That is, we have the  commutative diagram:

\begin{picture}(300,100)(-30,15)  
\put(128,100){$B\overline{\Diff_c(M)} \times M$}
\put(288,100){$B\G_q^+$}
\put(160,72){\vector(0,1){20}}
\put(124,78){$g_C \times id$}
\put(296,72){\vector(0,1){20}}

\put(46,20){$B\overline{\Diff_c(M)}$}
\put(124,90){\vector(-1,-1){50}}
\put(120,30){$g_C$}
\put(140,24){\vector(-1,0){40}}

\put(142,60){$Z_C \times M$}
\put(208,104){\vector(1,0){50}}
\put(248,80){\vector(2,1){30}}

\put(186,63){\vector(3,1){20}}
\put(264,63){\vector(-2,1){16}}
\put(208,70){$W_C \times M$}
\put(224,110){$\widetilde{\nu}_M$}
\put(238,88){$\widehat{b}_C$}

\put(156,54){\vector(0,-1){20}}
\put(162,34){\vector(0,1){20}}
\put(166,42){$\sigma$}
\put(146,42){$\pi$}
\put(156,20){$Z_C$}
 \put(274,60){$\mS^{n-k} \times M$}
\put(296,54){\vector(0,-1){20}}
\put(300,78){$\widehat{f}_C$}
\put(300,42){$\pi$}
 \put(286,20){$\mS^{n-k}$}
\end{picture}
 
 If the tangent bundle $TM$ is trivial, then we can substitute $B\oG_q$ for $B\G_q^+$ in the above diagram.

Set $\widehat{\phi}_C = \widehat{f}_C^*(\phi) \in H^{n}(\mS^{n-k} \times M; \mR)$ and $\widehat{\psi}_C = (g_C^* \circ \widetilde{\nu}_M^*(\phi) \in H^{n}(Z_C \times M ; \mR)$. 
 By the Kunneth formula,   we have
  \begin{equation}\label{eq-kunneth}
H^{*}(\mS^{n-k} \times M ; \mR)  \  \cong  \  H^{*}(\mS^{n-k} ; \mR) \otimes H^{*}(M ; \mR) 
\  , \   H^{*}(Z_C \times M ; \mR)  \  \cong  \  H^{*}(Z_C ; \mR) \otimes H^{*}(M ; \mR)  \ .  
\end{equation}
The slant product with the class $C \in H_k(M; \mZ)$ defines $\backslash C \colon H^{*}(M ; \mR) \to H^{*-k}(M ; \mR)$, and we compose this with projection $\Pi_1$ to the first factor in \eqref{eq-kunneth} to define:
\begin{equation}
\phi_C = \Pi_1( \widehat{\phi}_C \backslash C) \quad , \quad  \psi_C = \Pi_1( \widehat{\psi}_C \backslash C)  \ .
\end{equation}

Then as the coboundary $\widehat{b}_C$ is induced on the product on     factor of $M$, we have:
\begin{equation}
\langle \psi_C , [Z_C] \rangle   \ = \ \langle \widehat{\psi}_C , [Z_C] \otimes C \rangle   
  =    \langle \widehat{\phi}_C , [\mS^{n-k}] \otimes C \rangle   \ = \ \langle \phi_C , [\mS^{n-k}] \rangle  =   \langle \phi , [f(\mS^n)] \rangle
\end{equation}
This completes the proof of the claim in \eqref{eq-pairing}.
\endproof

 We   give a simple extension of Proposition~\ref{prop-homology}, which is used  in the following sections 
 to construct uncountable families of linearly independent homology classes in $H_*(B\overline{\Diff(M)} ; \mQ)$.

 \begin{prop}\label{prop-topclasses}
   Let $M$ be a compact, connected oriented manifold of dimension $q$, and assume there is given a co-spherical   class       $C \in H_k^{\pi}(\widehat{M}; \mZ)$.
Let $\{\phi_1, \ldots , \phi_m\} \subset H_{\pi}^n(B\G_q^+ \; \mR)$ for $n > 2q$, and let $\{f_i \colon \mS^n \to B\G_q^+ \mid i \in \cI\}$  be a collection of maps  such that the  homology classes they define, $\{[f_i]   \mid i \in \cI\} \subset H_n(B\G_q^+ \; \mQ)$, are rationally linearly independent for the vector mapping $\vec{\phi} = (\phi_1, \ldots , \phi_m) \colon H_n(B\G_q^+ ; \mQ) \to \mR^m$.
Then there exists $\{\psi_1, \ldots , \psi_m\} \in H^{n-k}(B\overline{\Diff_c(M)} ; \mR)$ and 
for each $i \in \cI$, there exists a compact oriented $(n-k)$-manifold $N_i$ and map $g_i \colon N_i \to B\overline{\Diff_c(M)}$ such that the classes 
$\{[g_i(N_i)] \mid i \in \cI\} \subset H_{n-k}(B\overline{\Diff_c(M)} ; \mQ)$ are rationally linearly independent in the image of the vector mapping $\vec{\psi}  \colon H_{n-k}(B\overline{\Diff_c(M)} ; \mQ) \to \mR^m$.
 \end{prop}
 \proof
 For     $\widetilde{\nu}_M \colon B\overline{\Diff_c(M)} \times M \to B\G_q^+$ and  $1 \leq \ell \leq m$, set $\widehat{\psi}_{\ell} = \widetilde{\nu}_M^*(\psi_{\ell}) \in H^n(B\overline{\Diff_c(M)} \times M; \mR)$  and    $\psi_{\ell} = \widehat{\psi}_{\ell} \backslash C \in H^{n-k}(B\overline{\Diff_c(M)}; \mR)$.

Observe that $\pi_n(B{\rm SO}(q)) \otimes \mQ$ is trivial for $n > 2q$. Thus for each $i \in \cI$, there exists an integer $m_i > 0$ such that the class $m_i \cdot [f_i] \in \pi_n(B\G_q^+)$ is in the image of $\pi_n(B\oG_q) \to \pi_n(B\G_q^+)$. So without loss of generality we can assume that each map $f_i$ factors through $B\oG_q$.

For each $i \in \cI$, apply the method of proof of Proposition~\ref{prop-homology} to the map $f_i$ and cycle $C$, to construct a manifold $N_i$ of dimension $n-k$, and map $g_i \colon N_i \to B\overline{\Diff_c(M)}$.  

Suppose that the    class
$\ds A = \sum_{i\in \cI} \ r_i \cdot [g_i(N_i)] \in H_{n-k}(B\overline{\Diff_c(M)} ; \mQ)$ is zero, where $r_i \in \mQ$ and only finitely many are non-zero. Then the class $\widetilde{\nu}_M(A \cup C) \in H_n(B\G_q^+ ; \mQ)$ is also zero.
 Then by \eqref{eq-pairing}, 
$$ \sum_{i\in \cI} \ r_i' \cdot  f_i(\mS^n) \  =  \  \widetilde{\nu}_M(A \cup C) \in H_n(B\G_q^+ ; \mQ) \ , $$
where $r_i'$ is a non-zero rational multiple of $r_i$ hence is again rational.
But this contradicts our assumption  that these classes are rationally linearly independent.
 \endproof
 
   Next,   consider the problem whether the classes constructed in Proposition~\ref{prop-homology} are non-trivial in the image of the map $(\iota_0)_* \colon H_*(B\overline{\Diff_c(M)} ; \mQ) \to  H_*(B\Diff_c^+(M)^{\delta} ; \mQ)$. 
We have the commutative diagram:

\begin{picture}(300,70)(-20,0)  
 
\put(86,50){$B\overline{\Diff_c(M)} \times M$}

\put(200,50){$M\Diff_c^+(M)^{\delta}$}
\put(160,53){\vector(1,0){30}}
\put(260,53){\vector(1,0){30}}
\put(120,40){\vector(0,-1){17}}
\put(224,40){\vector(0,-1){17}}
\put(96,10){$B\overline{\Diff_c(M)}$}
\put(230,32){$\pi$}
\put(110,32){$\pi$}

\put(270,60){$\tilde{\nu}_M$}
\put(306,50){$B\Gamma_q^+$}
\put(172,20){$\iota_{0}$}
\put(162,60){$\iota_0 \times id$}
\put(272,20){$\iota_{\delta}$}
\put(200,10){$B\Diff_c^+(M)^{\delta}$}
\put(160,14){\vector(1,0){30}}
\put(260,14){\vector(1,0){30}}
\put(296,10){$B\Diff_c^+(M)$}
\end{picture}
 
 Given a cycle $g_C \colon Z_C \to B\overline{\Diff_c(M)}$ we obtain a cycle 
 $\iota_0 \circ g_C \colon Z_C \to B\Diff_c^+(M)^{\delta}$, and so the question is whether the homology class it defines in 
 $H_{n-k}(B\Diff_c^+(M)^{\delta} ; \mQ)$ is non-trivial.
 
 For the case $C = M$, then the class $[g_M(Z_M)] \in H_{n-q}(B\Diff_c^+(M)^{\delta} ; \mQ)$ pairs non-trivially with the cohomology class $\psi_M = \int_M \  (\iota_0 \times id)^* \circ \widetilde{\nu}_M^*(\phi) \in H^{n-q}(B\overline{\Diff_c(M)} ; \mR)$. The integration over the fiber map is well-defined for the space $M\Diff_c^+(M)^{\delta}$ as well, so there is a class
 $\widetilde{\psi}_M = \int_M \ \widetilde{\nu}_M^*(\phi) \in H^{n-q}(B\Diff(M)_c^{\delta} ; \mR)$ such that $\iota_0^*(\widetilde{\psi}_M) = \psi_M$. Thus, the image $\iota_0(Z_M) \in H_{n-q}(B\Diff(M)_c^{\delta} ; \mQ)$ is non-trivial.

For the homology classes in $H_{n-k}(B\overline{\Diff_c(M)} ; \mQ)$ constructed using the methods of Theorems~\ref{thm-method2} and \ref{thm-method3}, for $0 \leq k < q$, the relation between the   groups $H_{n}(M\Diff_c^+(M)^{\delta} ; \mQ)$ and $H_{n-k}(B\Diff_c^+(M)^{\delta} ; \mQ)$ is   more complicated, and depends on the homology spectral sequence associated to the fibration \eqref{eq-basicfibration}. Even for the case when $M= \mS^1$, the above argument fails for the case $k=0$, as there is no class in degree $3$ in $H_*(B\Diff^+(\mS^1)^{\delta}; \mQ)$; see  Theorem~\ref{thm-MoritaEuler3} below.

 Next, recall from    Remark~\ref{rmk-Thurston} the alternative interpretation of the construction in Theorem~\ref{thm-method1}, which expands on a construction of Thurston in \cite{Thurston1974}.
  
  Choose $x \in M$ and let $V \subset M$ be an open neighborhood of $x$.   Choose   distinct points $\{x_1, \ldots , x_m\}$ in  $V$, and for each $1\leq \ell \leq m$ choose an open disk neighborhood $x_{\ell} \in U_{\ell} \subset V$, such that their closures are disjoint, $\overline{U}_i \cap \overline{U}_j = \emptyset$ for $i \ne j$. Moreover, assume that each closure $\overline{U}_{\ell}$ is diffeomorphic to the closed unit disk $\mD^q \subset \mR^q$, and choose a diffeomorphism $h_{\ell} \colon \mD^q \to \overline{U}_{\ell}$. Then the restriction $h_{\ell} \colon \mB^q \to U_{\ell}$ is a diffeomorphism of the open unit ball to $U_{\ell}$.
     
   Fix a diffeomorphism $\lambda \colon \mR^q \to \mB^q$, then define the inclusion mapping $\rho_{\ell} \colon  \overline{\Diff_c(\mR^q)} \to \overline{\Diff(M)}$, where a map $g \in \overline{\Diff_c(\mR^q)}$ is conjugated by the map $h_{\ell} \circ \lambda$ to a diffeomorphism of $U_{\ell}$ with compact support in  $U_{\ell}$, then extended as the identity on the rest of $M$. Note that the diffeomorphisms in the images of $\rho_i$ and $\rho_j$ commute if $i \ne j$, and thus we obtain a homomorphism
   \begin{equation}\label{eq-Oproduct}
\Upsilon_m = \rho_1 \times \cdots \times \rho_m \colon \overline{\Diff_c(\mR^q)} \times \cdots \times \overline{\Diff_c(\mR^q)} \to \overline{\Diff(M)} \ .
\end{equation}
Then the classifying map for $R_m$ yields a map on homology
     \begin{equation}\label{eq-Hproduct}
  (B\Upsilon_m)_*  \colon H_*(\overline{\Diff_c(\mR^q)} ; \mQ) \otimes \cdots \otimes H_*(\overline{\Diff_c(\mR^q)} ;\mQ) \to H_*(\overline{\Diff(M)} ; \mQ) \ .
\end{equation}
 By  Theorem~\ref{thm-MTT},  the adjoint induces an isomorphism $\tau_* \colon H_*(B\overline{\Diff_c(\mR^q)} ; \mQ) \to H_*(\Omega^q B\oG_q ; \mQ)$, where $\Omega^q$ denotes the $q$-fold based loops.
 
 Recall that if $X$ is a simply connected topological space, then the loop space $\Omega X$ is an H-space, and using the Pontrjagin product, its homology $H_*(\Omega X ; \mQ)$ is a free polynomial algebra over $\pi_*(\Omega X) \otimes \mQ$. 
Thus, $H_*(B\overline{\Diff_c(\mR^q)} ;\mQ) $ is a free polynomial algebra over $\pi_{*+q}(B\oG_q) \otimes \mQ$.
As we showed in Section~\ref{sec-spherical}, and Theorem~\ref{thm-lotsaclasses} in particular, the spaces $\pi_*(B\oG_q)$ are uncountably large, and so the same follows even more so for the groups $H_*(B\overline{\Diff_c(\mR^q)} ;\mQ)$.

Note that  the  representations $\rho_{\ell}$ and $\rho_{\ell'}$ are conjugate in $\overline{\Diff(M)}$, so 
 the product \eqref{eq-Hproduct} also induces a multiplication for the homology classes in the image of each map $\rho_{\ell}$.  
 Thus, as observed by Morita in \cite[Section~4]{Morita1984a}, the Pontrjagin product on $H_*(B\overline{\Diff_c(\mR^q)} ;\mQ)$ coincides with the ``spacial product'' induced from the representations $\rho_{\ell}$.
 
 The proof of Proposition~\ref{prop-homology} and the arguments above shows that each induced   map 
 $$(B\rho_{\ell})_* \colon H_{q+1}B\overline{\Diff_c(\mR^q)} ; \mQ)  \to H_{q+1}(B\Diff^+(M)^{\delta} ; \mQ)$$
  is injective. For the case when $M = \mS^1$, Morita shows in \cite{Morita1984a} that the   image of the product map \eqref{eq-Hproduct}   map is injective in $H_*(B\Diff^+(\mS^1)^{\delta} ; \mQ)$ for all $m >1$ as well. For the case when $M = \Sigma_g$ is a closed surface of sufficiently high genus $g \gg 2$,   Nariman shows in \cite[Theorem~0.17]{Nariman2017b} that again, the image of the product map is injective for all $m >1$.  
The image of the product map \eqref{eq-Hproduct} is unknown  for $m > 1$ in  all other cases.

 \section{Characteristic classes for dimension 1}\label{sec-q=1}

The results for    the groups  $H_*(B\overline{\Diff(\mR)}; \mZ)$, $H_*(B\overline{\Diff(\mS^1)}; \mZ)$  and $H_*(B\Diff^+(\mS^1)^{\delta}; \mZ)$ are all the result of  the construction of examples by Thurston, which  show that the Godbillon-Very class is spherical.    Morita gave a thorough analysis of  these groups   in his works  \cite{Morita1984a,Morita1985,Morita2001}.   
 In this section we recall these results, and the relation to the techniques developed in previous sections.
 
Recall that  $B\overline{\Gamma}_1 = B\G_1^+$  so there is a map
$\widehat{\nu} \colon B\overline{\Diff_c(\mR)} \times \mR \to B\G_1^+$ where $\widehat{\nu}$ maps   infinity to the basepoint given by the product foliation.  Mather proved in \cite{Mather1971,Mather1973} that the adjunct map $\tau \colon B\overline{\Diff_c(\mR)} \to \Omega B\G_1^+$ is a homology isomorphism, which is the first case of Theorem~\ref{thm-MTT}.  As $\overline{\Diff_c(\mR)}$ is a simple group \cite{Thurston1974}, the group $H_1(\Omega B\G_1^+ ; \mZ)$ vanishes, and it follows that   $\pi_1(\Omega B\G_1^+)$ vanishes, so that $\pi_3(B\G_1^+)$ is the first non-vanishing homotopy group of $B\G_1^+$, hence there is an isomorphism  
$H_3(B\G_1^+ ; \mZ) \cong \pi_3(B\G_1^+)$ by the Hurewicz Theorem. 

By Thurston \cite{Thurston1972} the Godbillon-Vey class is spherical, hence one has:

\begin{prop}\label{prop-uncountable1-q=1}
There exists a   collection of maps $\{f_i \colon \mS^3 \to B\G_1^+ \mid i \in \cI\}$ where $\cI$ is an uncountable set, and their classes  $\{[f_i] \mid i \in \cI\} \subset \pi_3(B\G_1^+)$ are linearly independent over $\mQ$. 
\end{prop}
\proof
The Godbillon-Vey class defines a surjection   
 $GV \colon \pi_3(B\Gamma_1^+) \cong H^3(B\Diff_c(\mR)^{\delta} ; \mZ) \to \mR$. Choose a Hamel basis $\{a_i \mid i \in \cI\} \subset \mR$ of $\mR$ over $\mQ$.
 Then for each $i \in \cI$ choose a map $f_i \colon \mS^3 \to B\G_1^+$ such that $\langle GV , [f_i]\rangle = a_i$.
\endproof

The next result follows using   standard properties of homology and  H-spaces.
\begin{thm}\label{thm-uncountable2-q=1}
The   cohomology ring  $H^*(\Omega B\Gamma_1^+ ; \mQ)$ contains a   free polynomial algebra with basis an uncountable collection $\{\phi_i \mid i \in \cI\} \subset H^2(\Omega B\Gamma_1^+ ; \mQ)$.
\end{thm}
\proof
We have $\pi_2(\Omega B\Gamma_1^+) \cong H_2(\Omega B\Gamma_1^+ ; \mZ)$ so $H^2(\Omega B\Gamma_1^+ ; \mQ) \cong \Hom(\pi_2(\Omega B\Gamma_1^+) , \mQ)$. For each $i \in \cI$, let $\phi_i \in H^2(\Omega B\Gamma_1^+ ; \mQ)$ be such that $\phi_i([f_j]) = \delta^j_i$ where $f_j \colon \mS^2 \to \Omega B\Gamma_1^+$ is the map chosen in Proposition~\ref{prop-uncountable1-q=1} and $\delta_i^j$ is the   indicatrix function.

 Let $\vec{n_k}$ denote an integer valued function on the $k$-tuples $I = (i_1, \ldots , i_k)$ with $i_{\ell} \in \cI$, which is zero except on a finite number of $I$. For each such $I$ set $\phi_I = \phi_{i_1} \wedge \cdots \wedge \phi_{i_k}$.
 
 Let $\ds \Phi(\vec{n_k}) = \sum_I ~ n_I \phi_I \in H^{2k}(\Omega B\Gamma_1^+ ; \mQ)$ be the finite sum.
 Suppose that $\Phi(\vec{n_k}) =0$, and let $I$ be a $k$-tuple with $n_I \ne 0$.
 Let $f_I \colon \mS^2 \times \cdots \times \mS^2 \to \Omega B\Gamma_1^+$ be the Pontrjagin product of the classes $\{f_{i_1}, \ldots , f_{i_k}\}$, then $0 = \Phi([f_I]) = n_I a_{i_1} \cdots a_{i_k}$ implies that $n_I = 0$, which is a contradiction.
\endproof

  \begin{cor}\label{cor-uncountable-q=1}
The groups $H^{2k}(B\Diff_c^+(\mR)^{\delta}; \mQ)$ are uncountable for all $k > 0$.
\end{cor}
\proof
By Mather's Theorem we have 
$$H^{2k}(B\Diff_c^+(\mR)^{\delta}; \mQ) \cong \Hom(H_{2k}(B\Diff_c^+(\mR)^{\delta}; \mZ) ; \mQ) \cong \Hom(\pi_{2k}(\Omega B\G_1^+) , \mQ) \cong H^{2k}(\Omega B\Gamma_1^+ ; \mQ) \ .$$

  For each non-zero polynomial $\Phi(\vec{n_k})$ as above, let $\overline{\Phi(\vec{n_k})} \in H^{2k}(B\Diff_c^+(\mR)^{\delta}; \mQ)$ be the corresponding class. 
Then for each $k > 0$, the collection of classes $\{\overline{\Phi(\vec{n_k})}\} \subset H^{2k}(B\Diff_c^+(\mR)^{\delta}; \mQ)$ span  a      vector space with uncountable dimension.
\endproof

Let $X$ be a space with basepoint $x_0 \in X$, let $\Omega X$ be the space of loops based at $x_0$. The reduced suspension $\Sigma \Omega X$ is the quotient space of $\Omega X \times [0,1]$ modulo   $(\Omega X \times \{0\} \cup \Omega X \times \{1\} \cup \{\widetilde{x_0}\} \times [0,1])$. The fundamental map $k \colon \Sigma \Omega X \to X$ is then given by $k(\omega, t) = \omega(t)$.  
The  homology suspension map 
 $E_* \colon  H_*(\Omega X ; \cA) \to H_{*+1}(\Sigma \Omega X ; \cA)$  is induced by the geometric suspension of singular simplices. The composition 
 \begin{equation}
\sigma = - k_*  \circ E \colon H_n(\Omega X ; \cA) \to H_{n+1}(X ; \cA)
\end{equation}
 is   induced  by the transgression map  of the path space fibration $\cP X \to X$. 
 Barcus and Meyer show:
 \begin{thm}\cite[Proposition~2.4]{BarcusMeyer1958} \label{thm-BM58}
 Let $X$ be    $r$-connected, then $\sigma$ is isomorphism for  $0 < n < 2r$.
 \end{thm}
 Apply this result to the $2$-connected space $X = B\G_1^+$   to obtain:
 \begin{cor}
 The   suspension $\sigma \colon H_2(\Omega B\G_1^+ ; \mZ) \to H_{3}(B\G_1^+ ; \mZ)$ is an isomorphism.
 \end{cor}
 This is just a restatement of the isomorphisms  used to construct the classes in Proposition~\ref{thm-uncountable2-q=1}. 

 \begin{prob}\label{prob-productsq=1}
 Show that $\sigma \colon H_{2n}(\Omega B\G_1^+ ; \mZ) \to H_{2n+1}(B\G_1^+ ; \mZ)$  is non-trivial for  $n > 1$.
 \end{prob}
A solution to Problem~\ref{prob-productsq=1} is related to a  solution of Problem~\ref{prob-Morita}. Morita observed   that the non-triviality of the map 
  $\overline{GV}^* \colon H^6(K(\mR^{\delta},3) ; \mZ) \to H^6(B\G_1^+ ; \mZ)$ follows from the vanishing of the rational Whitehead products of the classes chosen in the proof of  Proposition~\ref{prop-uncountable1-q=1}; see   \cite[Proposition~6.4]{Morita1985}).  By the Samelson result \cite[Theorem]{Samelson1953} (or \cite[Theorem~6.2]{Morita1985}) the product vanishes if the images of the suspension map $\sigma$ annihilates their Pontrjagin products.  

Next, we have the commutative diagram: 
 
\begin{picture}(300,100)(-30,0)  
\put(60,90){$\mS^1$}
\put(64,80){\vector(0,-1){17}}
\put(100,93){\vector(1,0){20}}
\put(30,50){$B\overline{\Diff(\mS^1)} \times \mS^1$}
\put(100,53){\vector(1,0){20}}
\put(108,60){$\iota_{1}$}

\put(64,40){\vector(0,-1){17}}
\put(60,23){\vector(0,1){17}}
\put(70,32){$\pi$}
\put(48,32){$\sigma$}
\put(36,10){$B\overline{\Diff(\mS^1)}$}
\put(100,14){\vector(1,0){20}}
\put(108,20){$\iota_{0}$}

\put(160,90){$\mS^1$}
\put(162,80){\vector(0,-1){17}}

\put(140,50){$M\Diff^+(\mS^1)^{\delta}$}
\put(200,53){\vector(1,0){30}}
\put(162,40){\vector(0,-1){17}}
\put(170,32){$\pi$}
 \put(212,60){$\tilde{\nu}$}
\put(240,50){$B\Gamma_1^+ $}
\put(212,20){$\iota_{\delta}$}
\put(140,10){$B\Diff^+(\mS^1)^{\delta}$}
\put(200,14){\vector(1,0){30}}
\put(235,10){$B\Diff^+(\mS^1) \cong B{\rm SO}(2)$}
\end{picture}

The composition $\tilde{\nu} \circ \iota_1$ induces the adjunction map  
$$\tau \colon B\overline{\Diff(\mS^1)} \to \cS(\mS^1) \cong \Maps(\mS^1,   B\G_1^+) = \Lambda B\G_1^+ \ ,$$
 where $\Lambda B\G_1^+$ denotes the free loop space, and 
 the map $\tau$ induces isomorphisms in homology.
  
 Recall there is a fibration $\Omega\G_1^+ \to \Lambda B\G_1^+ \to B\G_1^+$, which admits a section $s \colon B\G_1^+ \to \Lambda B\G_1^+$ where $s(x)$ is the constant loop at  $x$. It follows that 
 $\pi_*(\Lambda B\G_1^+) = \pi_*(\Omega B\G_1^+) \oplus \pi_*(B\G_1^+)$.
 Since $\pi_1(B\G_1^+) = \pi_2(B\G_1^+) = \{0\}$, by \cite{SVP1976} the rational minimal model for $\Lambda B\G_1^+$ contains as a DGA summand the free algebra $\Lambda [\pi_3(B\G_1^+)  \oplus  \pi_2(\Omega B\G_1^+)]$ with trivial differential.
  
 Let   $gv = \widetilde{\nu}^*(GV) \in H^2(M\Diff_c^+(\mS^1)^{\delta}, \mR)$, then   $\iota_1^*(gv)$  decomposes into classes $\alpha, \beta$ in the expression
\begin{equation}\label{eq-splittingq=1}
 \beta + \alpha \otimes [d\theta] \in H^3(B\overline{\Diff(\mS^1)} \times \mS^1 ; \mR) \cong H^3(B\overline{\Diff(\mS^1)} ; \mR) \oplus H^2(B\overline{\Diff(\mS^1)} ; \mR) \otimes H^1(\mS^1 ; \mR) \ ,
\end{equation}
where $[d\theta] \in H^1(\mS^1 ; \mR)$ is the generator, and  $\alpha$ is obtained by taking the cap product with the fundamental class $[\mS^1]$. That is, $\alpha = \int_{\mS^1} \ \iota_1^*(GV) \in H^2(B\overline{\Diff(\mS^1)} ; \mR)$ is the class constructed in the proof of Theorem~\ref{thm-method1}.
Morita gave the following expression for $\alpha$ in 
 \cite[Corollary~3.2]{Morita1984a} 
\begin{equation}\label{eq-MoritaInt1}
 \alpha = \int_{\mS^1} \  \iota_1^*(gv) =   \iota_0^* \int_{\mS^1} \  gv   \in H^2(B\overline{\Diff(\mS^1)};  \mR) \ .
\end{equation}

 Furthermore,   $\beta  = \sigma^* \circ \  \iota_1^*(gv) \in H^3(B\overline{\Diff(\mS^1)} ; \mR)$ is the class constructed  in the proof of Theorem~\ref{thm-method2}.

We recall two further results by Morita to complete the discussion of the dimension $1$ case.

  \begin{thm}  \cite[Theorem~4.3]{Morita1984a}\label{thm-MoritaEuler3}  
  For all $n \geq 1$,  the classes $\{\alpha^{n} , \alpha^{n-1} \beta\} \in  H^{*}(B\overline{\Diff(\mS^1)} ; \mR)$ are linearly independent, and  give independently variable maps on the     groups $H^{*}(B\overline{\Diff(\mS^1)} ; \mZ)$.
 \end{thm}
 \begin{cor}
For $n \geq 1$, there are surjective maps $\alpha^{n} \colon H_{2n}(B\overline{\Diff(\mS^1)} ; \mZ) \to \mR$ and 
$\alpha^{n-1} \beta \colon H_{2n+1}(B\overline{\Diff(\mS^1)} ; \mZ) \to \mR$.
 \end{cor}
   
 Finally, Morita showed that the image of the Euler class $\overline{\chi} \in H^2(B{\rm SO}(2) ; \mZ)$ and its powers are mapped injectively, 
 $\chi^n = \iota_{\delta}^*(\overline{\chi}^n) \in H^{2n}(B\Diff^+(\mS^1)^{\delta}  ; \mZ)$.  
Even more is true, 
 \begin{thm} \cite[Theorem~1.1]{Morita1984a}\label{thm-MoritaEuler2}  
  For all $n \geq 1$,  the classes $\{\alpha^n , \chi^n \} \in  H^{2n}(B\Diff^+(\mS^1)^{\delta} ; \mR)$ are linearly independent, and  the evaluation $\alpha^n \colon  H_{2n}(B\Diff^+(\mS^1)^{\delta} ; \mZ) \to \mR$ is a surjection.
  \end{thm}
The proof  in \cite[Section~7]{Morita1984a} uses that the space $\mS^1$ admits   finite   coverings  of arbitrary degree.  Nariman uses an extension of this idea in his work \cite{Nariman2023b}.

\section{Characteristic classes for dimension 2}\label{sec-q=2}

The analysis of    the homology and cohomology groups of the spaces   $B\overline{\Diff_c(M)}$  and $B\Diff_c^+(M)^{\delta}$ for an oriented   surface $M$   has several major differences from the $1$-dimensional case, and the study of $\cS(M)$ gets more interesting. We describe these results for $q=2$ in some detail, to compare and contrast them from the results in Section~\ref{sec-q=1}.

The first difference is that the spaces of secondary invariants, $H^*(WO_2)$ and $H^*(W_2)$ as introduced in Section~\ref{sec-secondary}, now consists of more classes than just the Godbillon-Vey invariant:
\begin{eqnarray*}
H^{*> 4}(WO_2) & = & \bigwedge (y_1 c_1^2 , \ y_1 c_2  ) \  , \ |y_1 c_1^2| = |y_1 c_2| = 5 \\
H^*(W_2)  & = & \bigwedge (y_1 c_1^2 , \ y_1 c_2 , \ y_2 c_2 , \ y_1 y_2 c_2)  \ ,  \ |y_1 c_1^2| = |y_1 c_2| = 5 , \ |y_2 c_2| = 7 ,  \ | y_1 y_2 c_2|  = 8 \ .
\end{eqnarray*}
 The author showed in \cite[Theorem~2.5]{Hurder1985a}  that both $\{\Delta(y_1c_1^2), \Delta(y_1c_2)\} \subset H^5(B\G_2^+ ; \mR)$ are spherical classes, and so give rise to  generalized Godbillon-Vey classes $\{gv_1, gv_2\} \subset H^5(M\Diff_c^+(M)^{\delta} ; \mR)$ defined as in \eqref{eq-gvq=2} below, and  the    classes $\{\alpha_1, \alpha_2\} \subset H^3(B\Diff_c^+(M) ; \mR)$ defined by integration over the fiber in \eqref{eq-alpha=2}. In addition, there are also non-trivial classes in  $H^4(B\overline{\Diff_c(M)} ; \mR)$ defined using co-spherical classes in $H_1(M; \mZ)$.
 
 The class $\Delta(y_2c_2) \in H^7(B\oG_2 ; \mR)$  is a ``rigid'' class, and it is unknown if it  is non-zero, though Pittie showed in \cite{Pittie1979} that it must vanish on locally homogeneous foliations.  The class $\Delta(y_1 y_2 c_1^2) \in H^8(B\overline{\G}_2 ; \mR)$
  is non-trivial, as shown by Baker \cite{Baker1978} for example. However, it is unknown if this class is non-vanishing on the image of the Hurewicz map $\pi_8(B\overline{\G}_2) \to H_8(B\overline{\G}_2 ; \mZ)$.
  
 The second difference is that for $M = \mS^2$ and $M = \Sigma_g$ a surface of genus $g \geq 2$, the tangent bundle is not trivial, so the method of Theorem~\ref{thm-method2} does not apply,  and the method of Theorem~\ref{thm-method3}  must be adjusted for its application. On the other hand, for $M = \mT^2$, all of the methods of Section~\ref{sec-constructions} can be applied, yielding many      additional non-trivial homology classes in $H_*(B\overline{\Diff(\mT^2)}; \mQ)$.

We    recall the known facts about  the spherical cohomology groups  $H_{\pi}^*(B\G_2^+ ; \mR)$ and $H_{\pi}^*(B\oG_2 ; \mR)$. 

Let $\widehat{\nu} \colon B\overline{\Diff(\mR^2)}\times \mR^2 \to B\overline{\Gamma}_2$ be the classifying map, where $\widehat{\nu}$ maps   infinity to the basepoint $w$ given by the point foliation.  The adjunct map $\tau \colon B\overline{\Diff(\mR^2)}  \to \Omega^2 B\overline{\Gamma}_2$ is a homology isomorphism by Thurston. As $\overline{\Diff(\mR^2)}$ is a simple group \cite{Thurston1974}, the group $H_1(B\overline{\Gamma}_2; \mZ)$ vanishes, and it follows that   $\pi_1(\Omega^2 B\overline{\Gamma}_2)$ vanishes, so that $\pi_4(B\overline{\Gamma}_2)$ is the first possibly non-vanishing homotopy group of $B\overline{\Gamma}_2$.

By the Hurewicz Theorem,  there is an isomorphism   $\pi_4(B\overline{\Gamma}_2) \cong H_4(B\overline{\Gamma}_2 ; \mZ)$. 
In addition,   the Rational Hurewicz Theorem~\ref{thm-RHIT} yields:
 
 \begin{prop}\label{prop-RHIT2}
 The Hurewicz map induces   isomorphisms
\begin{equation}
h \colon \pi_{i}(B\overline{\Gamma}_2) \otimes \mQ  \longrightarrow   H_{i}(B\overline{\Gamma}_2); \mQ)
\end{equation}
for $1 \leq i < 7$ and a surjection for $i=7$.
 \end{prop}
 
  Rasmussen constructed in \cite{Rasmussen1980} two different families of foliations, $\F_{1, \lambda}$ on $Y_1$ and  $\F_{2, \lambda}$ on $Y_2$, where $Y_1$ and $Y_2$ are oriented compact $5$-manifolds without boundary such that  evaluation of the secondary classes $\{y_1c_1^2, y_1 c_2\}$ on these foliations  vary independently, and thus evaluation of the universal secondary classes defines a map
\begin{equation}\label{eq-secondaryNVq=2}
\{ \Delta (h_1   c_1^2),  \Delta (h_1   c_2) \} \colon H_5(B\Gamma^+_2 ; \mZ) \to \mR^2 \ .
\end{equation}
The image contains an open subset of $\mR^2$, hence must be all of $\mR^2$. Note that there is no claim in \cite{Rasmussen1980} that the normal bundles to these foliations are trivial. Then by Theorem~\ref{thm-variablehomotopy} we have: 
 \begin{thm}\cite[Theorem~2.5]{Hurder1985a} \label{thm-nontrivialq=2}
The   classes $\{\Delta (h_1   c_1^2),  \Delta (h_1   c_2)\}$ yield  a surjection
\begin{equation}
\{ \Delta (h_1   c_1^2) , \Delta (h_1   c_2) \} \colon \pi_5(B\overline{\Gamma}_{2}) \to \pi_5(B\G^+_{2}) \to \mR^2  \ .
\end{equation}
\end{thm}
As in the proof of Proposition~\ref{prop-uncountable1-q=1},  we then have:
 
\begin{prop}\label{prop-uncountable1-q=2}
Choose a Hamel basis $\{a_i \mid i \in \cI\} \subset \mR$ of $\mR$ over $\mQ$.  Then 
there exists a   collection of maps $\{f_{i,j} \colon \mS^5 \to B\oG_2 \mid i,j \in \cI\}$  such that 
$$ \langle \Delta (h_1   c_1^2) , [f_{i,j}] \rangle = a_i  \ , \  \langle \Delta (h_1   c_2) , [f_{i,j}] \rangle = a_j   \  . $$
In particular, the classes
  $\{[f_{i,j}]   \mid i,j \in \cI\} \subset \pi_5(B\oG_2)$ are linearly independent over $\mQ$. 
\end{prop}

As in the proof of Theorem~\ref{thm-uncountable2-q=1}, we obtain:

 \begin{thm}\label{thm-uncountable2-q=2}
The   cohomology ring  $H^*(\Omega^2 B\overline{\Gamma}_2 ; \mQ)$ contains a   free exterior algebra with basis an uncountable collection $\{\phi_{i,j} \mid i,j \in \cI\} \subset H^3(\Omega^2 B\overline{\Gamma}_2 ; \mQ)$.
\end{thm}

 \begin{cor}\label{cor-uncountable-q=2}
The groups $H^{3\ell}(B\overline{\Diff_c(\mR^2)}; \mQ)$ are uncountable for all $\ell > 0$.
\end{cor}
We pose the analogous version of Problem~\ref{prob-productsq=1} for the suspension map of codimension $2$ foliations:
\begin{prob}\label{prob-productsq=2}
 Show that $\sigma_2 \colon H_{3\ell}(\Omega^2 B\overline{\G}_2 ; \mQ) \to H_{3\ell +2}(B\overline{\G}_2 ; \mQ)$  is non-trivial for  $\ell > 1$.
 \end{prob}

Recall that $\widetilde{\nu}_M   \colon M\Diff_c^+(M)^{\delta}   \to B\oG_2^+$ classifies the foliation on the total space over  $B\Diff_c^+(M)^{\delta}$.
Define the generalized Godbillon-Vey classes 
\begin{equation}\label{eq-gvq=2}
gv_1 =   \widetilde{\nu}_M^*(\Delta (h_1   c_1^2)) \ , \  gv_2 = \widetilde{\nu}_M^*(\Delta (h_1   c_2)) \ \in \ H^{5}(M\Diff_c^+(M)^{\delta} ; \mR) \ .
\end{equation}
   Define   classes $\{\alpha_1 , \alpha_2 \} \subset  H^3(B\Diff_c^+(M)^{\delta}; \mR)$ by integration over the fiber:
\begin{equation}\label{eq-alpha=2}
\alpha_1 = \int_M \   gv_1  \ , \     \alpha_2 = \int_M \   gv_2 \ .
\end{equation}

Next, let $\{f_{i,j} \colon \mS^5 \to B\oG_2 \mid i,j \in \cI\}$ be the collection of maps constructed in Proposition~\ref{prop-uncountable1-q=2}.   

By Proposition~\ref{prop-topclasses}, for each $i,j \in \cI$, there exists a closed compact oriented $3$-manifold $N_{i,j}$ and map $g_{i,j} \colon N_{i,j} \to B\overline{\Diff_c(M)}$ such that the classes 
$\{[g_{i,j}(N_{i,j})] \mid i,j \in \cI\} \subset H_{3}(B\overline{\Diff_c(M)} ; \mQ)$ are rationally linearly independent in the image of the   map $\vec{gv} =(gv_1 , gv_2)  \colon H_{3}(B\overline{\Diff_c(M)} ; \mQ) \to \mR^2$. 
By the comments at the end of Section~\ref{sec-flatbundles}, these classes inject into $H_3(B\Diff_c^+(M)^{\delta} ; \mQ)$, so we have:

  \begin{thm}\label{thm-topclassesq=2}
For a   connected oriented surface $M$, the classes 
\begin{equation}
\{[g_{i.j}(N_{i,j})] \mid i,j \in \cI\} \subset H_{3}(B\Diff_c^+(M)^{\delta} ; \mQ) 
\end{equation}
are linearly independent,  and detected by the classes $\{gv_1 , gv_2 \} \subset H^3(B\Diff_c^+(M)^{\delta} ; \mR)$.
 \end{thm}

We next detail the constructions of classes for surfaces, on a case by case analysis  by their genus.

For $M = \mS^2$,   the genus $g=0$ case, the classifying map $\nu_{\mS^2} \colon \mS^2 \to B{\rm SO}(2)$ is non-trivial, and $H_1(\mS^2 ; \mQ)$ is trivial, so the only results are those in Theorem~\ref{thm-topclassesq=2}:
\begin{prop}
$H_3(B\Diff^+(\mS^2)^{\delta} ; \mQ)$ contains an uncountably generated rational subgroup which maps onto $\mR^q$ by evaluation with the   Godbillon-Vey classes $\{gv_1, gv_2\} \in H^3(B\Diff^+(\mS^2)^{\delta} ; \mR)$.
\end{prop}

For $M = \mT^2$,   the genus $g=1$ case, there are many more constructions of classes in $H_*(\overline{\Diff(\mT^2)} ; \mQ)$ but the classes in $H_*(B\Diff^+(\mT^2)^{\delta} ; \mQ)$ are the same as for $\mS^2$.
 First,   Theorems~\ref{thm-method2} and \ref{thm-nontrivialq=2} yield:
  \begin{prop}\label{prop-bottomclasses2}
There is an injection
 \begin{equation}\label{eq-sectionmap2}
 \sigma^* \circ \  \widetilde{\nu}^* \  \colon    H^5_{\pi}(B\oG_q ;  \mR) \to H^{5}(B\overline{\Diff(\mT^2)} ; \mR) \ .
\end{equation}
In particular, for 
\begin{equation}\label{eq-gamma2}
\beta_1 =  \sigma^* \circ \iota_1^* \circ \widetilde{\nu}^*(\Delta (h_1   c_1^2))  \ , \  \beta_2 =  \sigma^*  \circ \iota_1^*  \circ \widetilde{\nu}^*(\Delta (h_1   c_2))
\end{equation}
the classes $\{\beta_1 , \beta_2 \} \subset H^{5}(B\overline{\Diff(\mT^2)} ; \mR)$ are linearly independent. 
 \end{prop}

Then by  Theorem~\ref{thm-nontrivialq=2} and Proposition~\ref{prop-uncountable1-q=2} we have: 
 
  \begin{thm}\label{thm-topclasses+q=2}
For each $i,j \in \cI$, there exists a compact oriented $5$-manifold $N'_{i,j}$ and map $\widehat{f_{i,j}'} \colon  N'_{i,j} \to B\overline{\Diff(\mT^2)}$ such that the homology classes 
$\{[\widehat{f_{i,j}'}(N'_{i,j})] \mid i,j \in \cI\} \subset H_{5}( B\overline{\Diff(\mT^2)} ; \mQ)$ are linearly independent, and detected by the  cohomology classes $\{ \beta_1 , \beta_2 \} \subset H^{5}(B\overline{\Diff(\mT^2)} ; \mR)$.
 \end{thm}

Next, we use  Theorem~\ref{thm-method3} and Proposition~\ref{prop-uncountable1-q=2} to construct linearly independent  classes in $H^4(B\overline{\Diff(\mT^2)} ;\mR)$.
Chose a factorization  $\mT^2 = \mS^1 \times \mS^1$, which defines the projections 
  $$p_1 \colon \mT^2 \to \mS^1 \times \{\theta_1\} = \mS^1  \ , \ 
  p_2 \colon \mT^2 \to   \{\theta_1\} \times \mS^1 \to \mS^1 \  $$
  and sections
  $$\eta_1 \colon  \mS^1  \cong  \mS^1 \times \{\theta_1\} \to  \mT^2  \ , \  
  \eta_2 \colon \mS^1  \cong  \{\theta_1\} \times \mS^1  \to  \mT^2  \  .$$
The sections define  co-spherical classes, again denoted by  $\eta_1, \eta_2 \in H_1(\mT^2 ;  \mZ)$.
 Then by the proof  of Theorem~\ref{thm-method3} define the classes
  \begin{equation}
\gamma^1_1 = \int_{\eta_1} \ gv_1  \ , \  \gamma^1_2 = \int_{\eta_1} \ gv_2   \ , \  \gamma^2_1 = \int_{\eta_2} \ gv_1   \ , \  \gamma^2_2 = \int_{\eta_2} \ gv_2 \ .
\end{equation}

 Next, combine Proposition~\ref{prop-topclasses} using $C= \mS^1 \times \{\theta_1\}$ or $C= \{\theta_1\} \times \mS^1$, with  Proposition~\ref{prop-uncountable1-q=2} and  Theorem~\ref{thm-method3} to conclude that the  set
$\{\gamma_1^1 , \gamma_2^1 , \gamma_1^2 , \gamma_2^2\} \subset  H^4(B\overline{\Diff(\mT^2)} ;\mR)$ is linearly independent, and there are four  families of uncountably many linearly independent classes in $H_4(B\overline{\Diff(\mT^2)} ;\mQ)$ which pair non-trivially with the $\gamma^i_j$-classes.  
  
 Thus, for $M = \mT^2$ we have defined the following linearly independent classes, for $1 \leq i,  j \leq 2$, which pair non-trivially with corresponding uncountable generated subgroups of $H_*(B\overline{\Diff(\mT^2)} ; \mQ)$: 
\begin{itemize}
\item  $gv_i \in \  H^{5}(B\overline{\Diff(\mT^2)} \times \mT^2 ; \mR)$. 
\item  $\alpha_i = \int_{\mT^2} \  gv_i    \  \in \  H^{3}(B\overline{\Diff(\mT^2)} ; \mR)$.  
\item $\gamma^i_j = \int_{\eta_i} \ gv_j  \ \in \ H^{4}(B\overline{\Diff(\mT^2)}   ; \mR)$.
\item $\beta_i = \sigma^* gv_i  \ \in \ H^{5}(B\overline{\Diff(\mT^2)}  ; \mR) $. 
\end{itemize}

The basic fibration sequence \eqref{eq-basicfibration} in this case becomes
$$ B\overline{\Diff(\mT^2)} \stackrel{\iota}{\longrightarrow} B\Diff^+(\mT^2)^{\delta} \stackrel{\iota_{\delta}}{\longrightarrow} B\Diff^+(\mT^2) \cong B{\rm SL}(2,\mZ)  \  . $$
 The  fundamental group ${\rm SL}(2,\mZ) = \pi_1(B\Diff^+(\mT^2) )$ acts on the homology group $H_1(\overline{\Diff(\mT^2} ; \mQ)$ and the only  invariant subspace is the trivial one. Thus,  the classes $\{\gamma_j^i\}  \subset  \ H^{4}(B\overline{\Diff(\mT^2)}   ; \mR)$ defined above by slant product with  the homology 1-cycles $\eta_i$  are killed in the homology spectral sequence for the fibration \eqref{eq-basicfibration}. That is, their images in $H_1(B\Diff^+(\mT^2)^{\delta} ; \mQ)$ are trivial.

For $M = \Sigma_g$,   the genus $g\geq 2$ case, there are many more constructions of classes in $H_*(\overline{\Diff(\Sigma_g)} ; \mQ)$ and the complexity of the study of $\cS(M)$ in higher dimensions $q \geq 2$ emerges.

The   classes  $\{ \alpha_1, \alpha_2 \} \subset H^3(B\Diff_c^+(\Sigma)^{\delta}; \mR)$  defined by \eqref{eq-alpha=2} are non-trivial by  Theorem~\ref{thm-nontrivialq=2}, and by Theorem~\ref{thm-topclassesq=2} they define a surjection of  $H_{3}(B\Diff^+(\Sigma_g)^{\delta} ; \mQ) $ onto $\mR^2$.

 The classifying map $\nu_{\Sigma_g} \colon \Sigma_g \to B{\rm SO}(2)$ is not trivial, as the surfaces have non-zero Euler class. However, by deleting a point, set  $\Sigma_g^* = \Sigma_g \setminus \{m_0\}$, and now the tangent bundle of $\Sigma_g^*$  is trivial.
 Also, the homology group $H_1(\Sigma_g^* ; \mZ) \cong \mZ^{2g}$ is generated by $2g$ simple closed curves $\{C_i \mid 1 \leq i \leq 2g\}$ based at a point $z_0 \ne m_0$ which give generators of $H_1(\Sigma_g^* ; \mZ)$.
 We then use Theorem~\ref{thm-method3} to construct classes, for $1 \leq i \leq 2g$ and $j=1,2$, 
 \begin{equation}
\gamma_j^i = \int_{C_i} \ gv_j \in H^4(B\overline{\Diff_c(\Sigma_g^*)} ; \mR)
\end{equation}

Correspondingly, there are $4g$  families of   linearly independent classes in $H_4(B\overline{\Diff_c(\Sigma_g^*)} ;\mQ)$ such that the classes $\{\gamma^i_j \mid 1 \leq i \leq 2g, j=1,2\}$ evaluate on them to give a surjection onto $\mR^{4g}$.  
Observe that there are inclusions 
\begin{equation}
H^4(B\overline{\Diff_c(\Sigma_g^*)} ;\mR) \subset H^4(B\overline{\Diff(\Sigma_g)} ;\mR) \ , \ H_4(B\overline{\Diff_c(\Sigma_g^*)} ;\mQ) \subset H_4(B\overline{\Diff(\Sigma_g)} ;\mQ) \ , 
\end{equation}
so we have:
\begin{thm}\label{thm-g=2}
For $g \geq 2$, there exists linearly independent classes 
\begin{equation}
\{\gamma^i_j \mid 1 \leq i \leq 2g, j=1,2\} \subset H^4(B\overline{\Diff(\Sigma_g)} ;\mR)
\end{equation}
which evaluate on $H_4(B\overline{\Diff(\Sigma_g)} ;\mQ)$ to yield a surjection onto $\mR^{4g}$.
\end{thm}
In addition, the curves $\{C_i \mid 1 \leq i \leq 2g\}$ can be chosen to be   intersect only in the basepoint $z_0$. It follows that the classes they define in  $H_4(B\overline{\Diff(\Sigma_g)} ;\mQ)$ have disjoint support, so we obtain a map analogous to the product map in \eqref{eq-Hproduct}. This yields non-trivial homology classes in $H_{4k}(B\overline{\Diff(\Sigma_g)} ;\mQ)$ for   $k$ arbitrarily large.
 
 The connected component $\Diff_0(\Sigma_g)$ is contractible, and the quotient ${\cM\cC\cG}_g = \Diff^+(\Sigma_g)/\Diff_0(\Sigma_g)$ is the \emph{mapping class group} for genus $g$ surfaces.   The group $\Diff_0(\Sigma_g)$ acts as the identity  on $H_1(\Sigma_g; \mZ)$ so the action of  $\Diff^+(\Sigma_g)$ on $H_1(\Sigma_g; \mQ)$ factors through an action of $\cM\cC\cG_g$. This actions are very well understood, and in particular there are no invariant subspaces for the action. Thus, the homology classes in Theorem~\ref{thm-g=2} map to trivial classes in $H_1(B\Diff^+(\Sigma_g)^{\delta}; \mQ)$

  Consider the map $\iota_{\delta}^* \colon H^*(B\Diff^+(\Sigma_g) ; \mR) \to H^*(B\Diff^+(\Sigma_g)^{\delta} ; \mR)$.
We have $B\Diff^+(\Sigma_g)  \to B{\cM\cC\cG}_g$ is a homotopy equivalence, 
and the cohomology groups $H^*(B{\cM\cC\cG}_g ; \mQ)$ have been extensively studied; see    Margalit \cite[Section~9]{Margalit2019} and Morita~\cite{Morita2006}. 
In a celebrated work \cite{Morita1984b,Morita1987}, Morita constructed many families of classes in $H^*(\Diff^+(\Sigma_g) ; \mQ)$ where $\Sigma_g$ is a compact surface of genus $g \geq 2$, and showed  the map $\iota_{\delta}^* \colon H^*(\Diff^+(\Sigma_g) ; \mQ) \to H^*(\Diff^+(\Sigma_g)^{\delta} ; \mQ)$  is not injective.

  We briefly recall the construction of the Morita-Mumford-Miller classes in $H^*(B\Diff^+(\Sigma_g) ; \mQ)$. The tangent bundle to the fibers of the universal map $\Sigma_g \to M\Diff^+(\Sigma_g)  \to B\Diff^+(\Sigma_g)$ are $2$-dimensional, and let $e_g \in H^2(M\Diff^+(\Sigma_g)) ; \mQ)$ denote its Euler class. Then integrate the powers of $e_g$ over the fibers to obtain the Morita-Mumford-Miller \emph{tautological classes}:
  \begin{equation}\label{eq-MMMclasses}
\kappa_{\ell} = \int_{\Sigma_g} \ e_g^{\ell+1} \in H^{2\ell}(B\Diff^+(\Sigma_g) ; \mQ) \cong H^{2\ell}({B\cM\cC\cG}_g ; \mQ) \ , \ \ \ell \geq 1 \ . 
\end{equation}
Miller \cite{MillerE1986} and Morita \cite{Morita1987} proved the following seminal result:
\begin{thm}\cite[Theorem~6.1]{Morita1987}
For $n \geq 1$ there exists a genus $g(n) \geq 2$, where $g(n) \leq 6n$,  and canonical relations on the graded polynomial ring $\mQ[\kappa_1, \ldots , \kappa_{g-2}]$ such that the map 
\begin{equation}
\mQ[\kappa_1, \ldots , \kappa_{g-2}]/{\rm relations} \to H^{*}(B\Diff^+(\Sigma_g) ; \mQ) 
\end{equation}
is injective up to degree $2n$ for all $g \geq g(n)$.
\end{thm}
A celebrated paper by  Madsen and Weiss proved the remarkable result:
\begin{thm}\cite{MadsenWeiss2007}\label{thm-MW}
For $n \geq 1$ and $g$ sufficiently large, there is an isomorphism
\begin{equation}
\{\mQ[\kappa_1, \ldots , \kappa_{g-2}]/{\rm relations}\}^{2n} \cong H^{2n}(B\Diff^+(\Sigma_g)  ; \mQ) 
\end{equation}
\end{thm}
This result was given alternative, ``more geometric'' proofs by   
 Eliashberg, Galatius and Mishachev \cite{EGM2011}
  and   Hatcher \cite{Hatcher2011}.  
  In contrast to the conclusion of Theorem~\ref{thm-MW}, Morita showed in the celebrated work \cite{Morita1987,Morita2001} that the higher MMM-classes always vanish in   $H^*(B\Diff^+(\Sigma_g)^{\delta} ; \mQ)$.
  \begin{thm}\label{thm-Morita-MMM1}
  The class $\kappa_{\ell} \in H^{6\ell}(B\Diff^+(\Sigma_g)^{\delta} ; \mQ)$ is trivial for $\ell \geq 3$.
     \end{thm}
     The proof that $\kappa_{\ell} = 0$ follows from the Bott Vanishing Theorem. Observe that  $e_g^2 = p_i$ is the first Pontrjagin class of the tangent bundle to the fibers, and so $e_g^{\ell+1} \in H^{6\ell +2}(B\Diff^+(\Sigma_g)^{\delta} ; \mQ)$ vanishes for $\ell \geq 3$. It is unknown whether $\kappa_2 = 0$, as it is one of the intriguing open problems whether the cube of the Euler  class $e_g^{3} \in H^{6}(B\G_2^+ ; \mQ)$ must vanish.  

In contrast to this, Kotschick and Morita   proved in \cite{KotschickMorita2005} the following result, which was given an alternate proof by  Nariman \cite[Theorem~0.12]{Nariman2017b}:
      \begin{thm}\label{thm-Morita-MMM2}
  There is an inclusion
  $\mQ[\kappa_1] \subset H^{2\ell}(B\Diff^+(\Sigma_g)^{\delta} ; \mQ)$ stably. That is, given $\ell \geq 3$ there exists $g_{\ell}$ so that for $g \geq g_{\ell}$ the class   $\kappa_1^{\ell} \ne 0$ in $H^{2\ell}(B\Diff^+(\Sigma_g)^{\delta} ; \mQ)$.
       \end{thm}

 As discussed by Margalit \cite[Section~9]{Margalit2019}, the cohomology above the range of Theorem~\ref{thm-MW} contains many classes not of the form of the MMM-classes, and in fact the constructions of these classes show they grow exponentially as $g$ increases.
   The general problem is thus:
  \begin{prob}\label{prob-MMM1}
  For fixed $g \geq 2$, determine which classes in $H^*(B{\cM\cC\cG}_g ; \mQ)$ have non-zero image under the map
  $\iota_{\delta}^* \colon H^*(B{\cM\cC\cG}_g ; \mQ) \cong H^*(B\Diff^+(\Sigma_g)  ; \mQ) \to H^*(B\Diff^+(\Sigma_g)^{\delta} ; \mQ)$.
    \end{prob}
There are further results related to Problem~\ref{prob-MMM1} in the works by Kotschick and Morita in \cite{KotschickMorita2005}, Bowden in \cite{Bowden2012} and Nariman in \cite{Nariman2017a,Nariman2017b}.

\section{Characteristic classes for dimension 3}\label{sec-q=3}

Let $M$ be a connected oriented $3$-manifold. The tangent bundle to $M$ is   trivial, so   the methods in  Section~\ref{sec-constructions}  can be applied to construct families of classes in $H_*(B\overline{\Diff_c(M)}; \mQ)$. Moreover,  for $M$ of dimension $q \geq 3$, there are additional families of spherical classes in $H^*(B\G_q^+ ; \mQ)$ beyond    the Godbillon-Vey type classes in $H^{2q+1}(B\oG_3; \mR)$. 
For the case when $q=3$, these   new spherical classes appear in degree $H^{10}(B\oG_3 ; \mR)$ which give rise to even more families of  homology classes in $H_*(B\overline{\Diff_c(M)}; \mQ)$.
First, we recall the Vey basis \cite{Godbillon1974}: 
\begin{eqnarray*}
H^{*>6}(WO_3) & = & \bigwedge (y_1 c_1^3 , \ y_1 c_1 c_2 , \ y_1 c_3,  \ y_1 y_3 c_1^3 , \ y_1 y_3 c_1 c_2 ,  \ y_3 c_3 )  \\
H^*(W_3)  & = & \bigwedge (y_1 c_1^3 , \ y_1 c_1 c_2 , \  y_1 c_3, \  y_2 c_2 , \ y_1 y_2 c_1^3 , \ y_1 y_2 y_3 c_1 c_2 , \ y_1 y_2 y_3c_1^3 , \ y_1 y_2 c_1 c_2 ,  \  y_3c_3)    \ .
\end{eqnarray*}

 Theorem~\ref{thm-variablehomotopy} implies the classes $\{\Delta(y_1c_1^3), \Delta(y_1 c_1 c_2) , \Delta(y_1 c_3) \} \in H^7(B\G_3^+ ; \mR)$ are spherical classes, as was previously  observed in    \cite{Hurder1980,Hurder1981a}. The linear independence of the classes in the image of $\Delta \colon H^{7}(W_3) \to H^7(B\oG_3 ; \mR)$ follows from the works cited for the proof of Theorem~\ref{thm-injectivityframed}.   
  
In addition, for dimension $q=3$, Theorem~\ref{thm-bracing} yields additional families of spherical classes as follows.  The first Pontrjagin class $p_1 \in H^4(B{\rm SO}(3) ; \mR)$ is spherical, detected by a map $\xi_4 \colon \mS^4 \to B{\rm SO}(3)$.  As $B\oG_3$ is $4$-connected, there is a lift $\widehat{\xi}_1 \colon \mS^4 \to B\G_3^+$ with $\nu \circ \widehat{\xi}_1 = \xi_1$.
The brace product with the class $[\widehat{\xi}_1] \in \pi_4(B\G_3^+)$ as defined in \eqref{eq-brace} gives a mapping
\begin{equation}\label{eq-braceq=3}
\{ [\widehat{\xi}_1], \cdot  \}  \colon   \pi_j(B\overline{\Gamma}_3) \to \pi_{j+3}(B\overline{\Gamma}_3) \ .
\end{equation}
The classes $\{\Delta(y_1 c_1^3), \Delta(y_1 c_1 c_2), \Delta(y_1 c_3)\} \subset H_7(B\oG_3 ; \mR)$ are spherical, so using the mapping \eqref{eq-braceq=3}  the classes 
$\{\Delta(y_1 y_2 c_1^3), \Delta(y_1 y_2 c_1 c_2), \Delta(y_1 y_2 c_3)\} \subset H^{10}(B\oG_3 ; \mR)$
are also spherical.  
Using the procedures of the previous sections, we then obtain:
   
 \begin{thm}  \label{thm-nontrivialq=3}
The   classes $\{\Delta(y_1c_1^3), \Delta(y_1 c_1 c_2) , \Delta(y_1 c_3) \} \in H^7(B\G_3^+ ; \mR)$    yield  a surjection
\begin{equation}
\{ \Delta (y_1c_1^3) , \Delta (y_1 c_1 c_2) , \Delta (y_1 c_3) \} \colon \pi_7(B\overline{\Gamma}_3) \to \pi_7(B\Gamma^+_3) \to \mR^3  \ .
\end{equation}
\end{thm}

\smallskip

\begin{prop}\label{prop-uncountable1-q=3}
Choose a Hamel basis $\{a_i \mid i \in \cI\} \subset \mR$ of $\mR$ over $\mQ$.  Then 
there exists a   collection of maps $\{f_{i,j,k} \colon \mS^7 \to B\oG_3 \mid i,j,k \in \cI\}$  such that 
$$ \langle \Delta (y_1   c_1^3) , [f_{i,j,k}] \rangle = a_i  \ , \  \langle \Delta (y_1 c_1  c_2) , [f_{i,j,k}] \rangle = a_j 
 \ , \  \langle \Delta (y_1 c_3) , [f_{i,j,k}] \rangle = a_k    \  . $$
In particular, the classes
  $\{[f_{i,j,k}]   \mid i,j,k \in \cI\} \subset \pi_7(B\oG_3)$ are linearly independent over $\mQ$. 
\end{prop}

 \smallskip

 \begin{thm}\label{thm-uncountable2-q=3}
The   cohomology ring  $H^*(\Omega^3 B\overline{\Gamma}_3 ; \mQ)$ contains a   free exterior algebra with basis an uncountable collection $\{\phi_{i,j,k} \mid i,j,k \in \cI\} \subset H^4(\Omega^3 B\overline{\Gamma}_3 ; \mQ)$.
\end{thm}

 \begin{cor}\label{cor-uncountable-q=3}
The groups $H^{4\ell}(B\overline{\Diff_c(\mR^3)}; \mQ)$ are uncountable for all $\ell > 0$.
\end{cor}

 \medskip

\begin{prob}\label{prob-productsq=3}
 Show that $\sigma_3 \colon H_{4\ell }(\Omega^3 B\overline{\G}_3 ; \mQ) \to H_{4\ell+3}(B\overline{\G}_3 ; \mQ)$  is non-trivial for  $\ell > 1$.
 \end{prob}

Recall that $\widetilde{\nu}_M   \colon M\Diff_c^+(M)^{\delta}   \to B\oG_3^+$ classifies the foliation on the total space over  $B\Diff_c^+(M)^{\delta}$.
Define the generalized Godbillon-Vey classes in $H^7(M\Diff_c^+(M)^{\delta} ; \mR)$:
\begin{equation}\label{eq-gvq=3}
gv_1 =   \widetilde{\nu}_M^*(\Delta (h_1   c_1^3))  \ , \ \  gv_2 = \widetilde{\nu}_M^*(\Delta (h_1 c_1  c_2))    \  , \ \  gv_3 = \widetilde{\nu}_M^*(\Delta (h_1 c_3))  \ .
\end{equation}
   Define   classes $\{\alpha_1 , \alpha_2 , \alpha_3 \} \subset  H^4(B\Diff_c^+(M)^{\delta}; \mR)$ by integration over the fiber:
\begin{equation}\label{eq-alpha=3}
\alpha_1 = \int_M \   gv_1  \ , \    \alpha_2 = \int_M \   gv_2  \ , \   \alpha_3 = \int_M \   gv_3 \ .
\end{equation}

 Next, $\{f_{i,j,k} \colon \mS^7 \to B\oG_3 \mid i,j,k \in \cI\}$   be the collection of maps constructed in Proposition~\ref{prop-uncountable1-q=3}.   
 By Proposition~\ref{prop-topclasses}, for each $i,j,k \in \cI$, there exists a compact oriented $4$-manifold $N_{i,j,k}$ and map $g_{i,j,k} \colon N_{i,j,k} \to B\overline{\Diff_c(M)}$ such that the classes 
$\{[g_{i,j,k}(N_{i,j,k})] \mid i,j,k \in \cI\} \subset H_{4}(B\overline{\Diff_c(M)} ; \mQ)$ are rationally linearly independent in the image of  $\vec{\alpha} =(\alpha_1 , \alpha_2, \alpha_3)  \colon H_{4}(B\overline{\Diff_c(M)} ; \mQ) \to \mR^3$. 
By the comments at the end of Section~\ref{sec-flatbundles}, the image of these classes   inject into $H_4(B\Diff_c^+(M)^{\delta} ; \mQ)$:

  \begin{thm}\label{thm-topclassesq=3}
For a   connected oriented $3$-manifold $M$, the cycles 
\begin{equation}
\{[g_{i.j,k}(N_{i,j,k})] \mid i,j,k \in \cI\} \subset H_{4}(B\Diff_c^+(M)^{\delta} ; \mQ) 
\end{equation}
are linearly independent,  and evaluation of  the classes $\{\alpha_1 , \alpha_2, \alpha_3 \} \subset H^4(B\Diff_c^+(M)^{\delta} ; \mR)$ on these cycles yields a surjection onto $\mR^3$.
 \end{thm}

Next, as $M$ has trivial tangent bundle,  Theorem~\ref{thm-method2}  yields additional classes in $H^*(B\oG_q^+ ; \mR)$.   
For the section $\sigma \colon B\overline{\Diff_c(M)} \to B\overline{\Diff_c(M)} \times M$, define   classes in $H^7(B\overline{\Diff_c(M)} ; \mR)$:
\begin{equation}
\beta_1 = \sigma^* (gv_1)  \ , \   \beta_2 = \sigma^* (gv_2)  \ , \ \beta_3 = \sigma^* (gv_3) \ .
\end{equation}

Then  modify the maps 
   $\{f_{i,j,k} \colon \mS^7 \to B\oG_3 \mid i,j,k \in \cI\}$   constructed in Proposition~\ref{prop-uncountable1-q=3},  using the Mather-Thurston Theorem,  to  obtain:

  \begin{thm}\label{thm-topclasses+q=3}
For   $i,j,k \in \cI$, there exists a compact oriented $7$-manifold $N'_{i,j,k}$ and map $\widehat{f}_{i,j,k}'  \colon N'_{i,j,k} \to B\overline{\Diff_c(M)}$ such that the  cycles 
$\{[\widehat{f}_{i,j,k}'(N'_{i,j,k})] \mid i,j,k \in \cI\} \subset H_{7}( B\overline{\Diff_c(M)}  ; \mQ)$ are linearly independent,   and evaluation of  the classes $\{ \beta_1 , \beta_2 , \beta_3\} \subset H^{7}(B\overline{\Diff_c(M)} ; \mR)$ on these cycles yields a surjection onto $\mR^3$.
 \end{thm}

  Using the brace product and proceeding as before, we obtain:

 \begin{thm}  \label{thm-nontrivialq=3-2}
 $\{\Delta(y_1 y_2 c_1^3), \Delta(y_1 y_2 c_1 c_2) , \Delta(y_1 y_2 c_3) \} \in H^{10}(B\oG_3 ; \mR)$    yield  a surjection
\begin{equation}
\{ \{ \Delta(y_1 y_2 c_1^3), \Delta(y_1 y_2 c_1 c_2), \Delta(y_1 y_2 c_3)\} \colon \pi_{10}(B\overline{\Gamma}_3) \to \pi_{10}(B\Gamma^+_3) \to \mR^3  \ .
\end{equation}
\end{thm}
\smallskip

  \begin{thm}\label{thm-topclasses+q=3-2}
For   $i,j,k \in \cI$, there exists a compact oriented $10$-manifold $N''_{i,j,k}$ and map $\widehat{f}_{i,j,k}''  \colon N''_{i,j,k} \to B\overline{\Diff_c(M)}$ such that the   cycles 
$\{[\widehat{f}_{i,j,k}''(N''_{i,j,k})] \mid i,j,k \in \cI\} \subset H_{10}( B\overline{\Diff_c(M)}  ; \mQ)$ are linearly independent, and evaluation of  the classes
$$\{ (\iota_0 \circ \sigma)^*(\Delta(y_1 y_2 c_1^3)) ,  (\iota_0 \circ \sigma)^*(\Delta(y_1 y_2 c_1 c_2)) ,  (\iota_0 \circ \sigma)^*(\Delta(y_1 y_2 c_3))\} \subset H^{10}(B\overline{\Diff_c(M)} ; \mR) $$
on these cycles yields a surjection onto $\mR^3$.
 \end{thm}

 Using the methods of Theorem~\ref{thm-method3}, when  $H^1(M ; \mQ)$ is non-trivial, additional families of classes can be constructed in $H_{6}( B\overline{\Diff_c(M)}  ; \mQ)$ and $H_{9}( B\overline{\Diff_c(M)}  ; \mQ)$,  for each class $0 \ne Z \in  H^1(M ; \mZ)$ defined by a mapping $M \to \mS^1$. Details are left to the reader.

Note that the cycles constructed in the proofs of Theorems~\ref{thm-topclasses+q=3} and \ref{thm-topclasses+q=3-2},and using the methods of Theorem~\ref{thm-method3} are all supported in the interior of $M$. A $3$-manifold admits a decomposition into its irreducible components
\begin{equation}
M = M_1 \ \#_{\mS^2} \ M_2 \ \#_{\mS^2} \  \cdots \ \#_{\mS^2} \  M_{k} \ .
\end{equation}
Then we can apply the above results to construct families of cycles for each irreducible component $M_{\ell}$ and they inject into $H_*(B\overline{\Diff_c(M)} ; \mQ)$. One can thus apply these procedures in a ``mix-and-match'' construction. 

The classes obtained using Theorem~\ref{thm-topclassesq=3} are non-trivial in $H_4(B\Diff_c^+(M)^{\delta} ; \mQ)$, but for the other classes constructed above in $H_*(B\overline{\Diff_c(M)} ; \mQ)$, it seems to be completely unknown whether any of them survive for the map to $H_*(B\Diff_c^+(M)^{\delta} ; \mQ)$.

Finally, assume that $M$ is compact, and consider  the basic fibration  sequence \eqref{eq-basicfibration} and its 
 associated Leray-Serre spectral sequence, with $E_2$-term 
 \begin{equation}\label{eq-ss3}
E_2^{r,s}(B\Diff^+(M)^{\delta} ; \mQ) \cong H^r(B\Diff^+(M) ; H^s(B\overline{\Diff(M)}) ; \mQ) \ . 
\end{equation}
The topological type of $\Diff^+(M)$ has been determined for many compact $3$-manifolds   \cite{Hatcher2012},  so an analysis  of \eqref{eq-ss3} may be feasible for such manifolds. We consider some select cases.

Let $M = \mS^3$, then   the inclusion ${\rm SO}(4) \subset \Diff^+(\mS^3)$ is a homotopy equivalence by Hatcher \cite{Hatcher1983}. Thus $H^*(B\Diff^+(\mS^3) ; \mQ) \cong \mQ[e_4,p_1]$ is a polynomial ring with generators the first Pontrjagin class $p_1$  and the Euler class $e_4$, both of degree $4$. Let $\widetilde{p}_1 = \iota_{\delta}^*(p_1)$ and $\widetilde{e}_4 = \iota_{\delta}^*(e_4)$ in  $H^4(B\Diff^+(\mS^3)^{\delta} ; \mQ)$. 
Nariman has proved the following remarkable result:
\begin{thm}\cite[Theorem~1.2]{Nariman2023b}
There is an inclusion $\mQ[\widetilde{e}_4 , \widetilde{p}_1] \subset  H^4(B\Diff^+(\mS^3)^{\delta} ; \mQ)$.
\end{thm}

Let $M = \mS^1 \times \mS^2$,  Hatcher showed in  \cite{Hatcher1981} that 
 the group $\Diff(\mS^1 \times \mS^2)$ has the homotopy type of ${\rm O}(2) \times {\rm O}(3) \times \Omega {\rm O}(3)$, though the H-space structure is not preserved.
 \begin{prob}
 Determine the image of   $\iota_{\delta}^* \colon H^*(B\Diff(\mS^1 \times \mS^2) ; \mQ) \to H^*(B\Diff(\mS^1 \times \mS^2)^{\delta} ; \mQ)$.
 \end{prob}

Let $M = \mT^3 = \mS^1 \times \mS^1 \times \mS^1$,  Hamstrom showed in  \cite{Hamstrom1965} that the homotopy groups $\pi_{\ell}(\Diff(\mT^3))$ are trivial for $\ell > 1$. On the other hand, there are inclusions $\mT^3 \subset \Diff_0(\mT^3)$ and ${\rm SL}(3,\mZ) \to \Diff^+(\mT^3)$, and their images generate a subgroup isomorphic to the semi-direct product $G = \mT^3 \rtimes {\rm SL}(3,\mZ)$. The following problem seems approachable.
\begin{prob}
 Determine the image  of  $\iota_{\delta}^* \colon H^*(B\Diff(\mT^3) ; \mQ) \to H^*(B\Diff(\mT^3)^{\delta} ; \mQ)$.
 \end{prob}

Finally, let $M_{\kappa}$ denote a compact oriented $3$-manifold without boundary, and constant curvature $\kappa = \pm1$. Then Bamler and Kleiner determine the homotopy type of $\Diff^+(M_{\kappa})$ in  \cite{BamlerKleiner2023}; see also the monograph  \cite{HKMR2012}.  Again, it seems that nothing is known about the image of the map $\iota_{\delta}^*$ in this case.

\section{Characteristic classes for higher dimensions}\label{sec-q>3}

Let $M$ be an oriented connected manifold  of dimension  $q > 3$. The constructions of classes in $H^*(B\Diff_c^+(M); \mR)$ and $H_*(B\overline{\Diff_c(M)} ; \mQ)$ proceeds using a combination of the methods described in Sections~\ref{sec-NVsecondary},  \ref{sec-spherical} and \ref{sec-flatbundles}, but with added complexity due to the increasing number of linearly independent spherical classes in $H^{2q+1}(B\G_q^+ ; \mR)$ and $H^{2q+1}(B\oG_q ; \mR)$. In addition, there are   multiple brace product operations on $\pi_*(B\oG_q)$ that create additional spherical classes. We give a broad description of the constructions that are possible, but more precise statements require fixing the dimension and topological type of the manifold $M$.

For manifolds of dimension $q \geq 5$ there   also arise a   variety of non-trivial classes in $H^*(B\Diff_c^+(M) ; \mQ)$ constructed using pseudo-isotopy theory. The relationship of these with $H^*(B\Diff_c^+(M)^{\delta} ; \mQ)$ is completely unknown.
 It would be especially interesting to understand the images under $\iota_{\delta}^*$ of the exotic classes described in the works  \cite{GRW2014b,KRW2000,Tillmann2012}.

First, we describe the general form of the results obtained in previous sections.
Let $\cV_q$ be the collection of elements $y_i c_J \in WO_q$ of the Vey basis for $H^{2q+1}(WO_q)$ which are independently variable as   in Theorems~6.1 and 6.2, \cite{Heitsch1978}. Let $v_q$ denote the cardinality of the set $\cV_q$. 

Theorem~\ref{thm-variablehomotopy}  implies that  $\Delta(\cV_q)$ is a linearly independent set of spherical classes in $H^{2q+1}(B\G_q^+ ; \mR)$,  whose evaluation on 
$H_{2q+1}(B\G_q^+ ; \mQ)$ define surjections onto $\mR^{v_q}$.

Define
$\alpha_{i,J} = \int_M \ \widetilde{\nu}^*(\Delta(y_i c_J)) \in H^{q+1}(B\Diff_c^+(M)^{\delta} ; \mR)$ for $y_i c_J \in \cV_q$.

By Proposition~\ref{prop-topclasses}, there are $v_q$ families of uncountably generated subgroups of $H_{q+1}(B\overline{\Diff_c(M)} ; \mQ)$, and the evaluation of 
$\vec{\alpha}$ on the images of these classes in $H_{q+1}(B\overline{\Diff_c(M)} ; \mQ)$ yields a surjective map onto $\mR^{v_q}$.
In particular, these classes inject into $H_{q+1}(B\Diff_c^+(M)^{\delta} ; \mQ)$.

Next, consider the application  of the brace product    as in Theorem~\ref{thm-bracing}.  Let $\kappa_q$ be the greatest integer such that $4\kappa_q \leq q+1$. The first solution $\kappa_3 = 1$ appears when $q \geq 3$, and a second solution $\kappa_7 = 2$ appears when $q \geq 7$, and so forth.   Then for each $1\leq \ell \leq \kappa_q$ there exists $\xi_{\ell} \colon \mS^{4\ell} \to B\G_q^+$ such that the pull-back   $\xi_{\ell}^* \circ \nu^*(p_{\ell}) \ne 0$, where  $p_{\ell} \in H^{4\ell}(B{\rm SO}(q) ; \mZ)$ is the $\ell$-th Pontrjagin class.
Then the brace product defines a mapping $\{\xi_{\ell}, \cdot \} \colon \pi_j(B\overline{\Gamma}_q) \to \pi_{j + 4\ell  -1}(B\overline{\Gamma}_q)$.  

Define the extended subset of the Vey basis for $H^*(WO_q)$:
\begin{equation}
\widehat{\cV}_q = \{y_I y_{I'} c_J \mid  y_I c_J \in \cV_q \ , \   I' = (i_1' < \cdots <  i_r') \ {\rm each} \ i_{\ell} \ {\rm even \ with} \ 2i_r' \leq \kappa_q\}
\end{equation}
Let $\widehat{\cV}_{q,k} \subset \widehat{\cV}_q$ denote the terms of degree $k$, and let $\widehat{v}_{q,k}$ denote the cardinality of the set $\widehat{\cV}_{q, k}$. Then as in previous sections, we have the  non-triviality results as follows.

 \begin{thm}  \label{thm-nontrivialq>3}
For $k \geq 2q+1$, the   classes $\Delta(\widehat{\cV}_{q,k})   \in H^k(B\oG_q ; \mR)$    yield  a surjection
\begin{equation}
\Delta(\widehat{\cV}_{q,k})  \colon \pi_k(B\oG_q) \to  \mR^{\widehat{v}_{q,k}}  \ .
\end{equation}
\end{thm}

\medskip

\begin{prop}\label{prop-uncountable1-q>3}
Choose a Hamel basis $\{a_i \mid i \in \cI\} \subset \mR$ of $\mR$ over $\mQ$.  Then for each $k \geq 2q+1$ with $\widehat{v}_{q,k} > 0$,  there exists a   collection of maps $\{f_{I,I',J,i} \colon \mS^k \to B\oG_q \mid y_i y_{I'} c_J \in  \widehat{\cV}_{q,k} \ , \ i\in \cI\}$  such that 
$  \langle \Delta (y_i y_{I'} c_J ) , [f_{I,I',J,i}] \rangle = a_i$, and evaluates to zero on these spherical cycles   otherwise.
In particular, the classes
  $\{[f_{I,I',J,i}] \} \subset \pi_k(B\oG_q)$ are linearly independent over $\mQ$. 
\end{prop}

As before, this implies that the     cohomology ring  $H^*(\Omega^q B\oG_q ; \mQ)$ contains a   free exterior algebra with basis indexed by the set $\widehat{\cV}_{q,k} \times \cI$, and correspondingly the groups $H^{*}(B\overline{\Diff_c(\mR^q)}; \mQ)$ are uncountable.

 In this manner, corresponding to  the variable spherical classes in $H^{2q+1}(W_q)$ we obtain families of infinite dimensional subspaces of $\pi_{k}(B\oG_q)$ for $k=2q+4,2q+8,2q+11,2q+12, \ldots$, and thus  subspaces of $H_{k-q}(B\overline{\Diff_c(M)} ; \mQ)$ with an uncountable number of rationally independent generators. This just accounts for the classes obtained by integration over the fiber method  in Theorem~\ref{thm-method1}. There are additional homology classes corresponding to the methods in Theorems~\ref{thm-method2} and\ref{thm-method3}.

\begin{remark}\label{rmk-wtf}
{\rm
The results on the homology groups $H_{*}(B\overline{\Diff_c(M)} ; \mQ)$ and $H_{*}(B\Diff_c^+(M)^{\delta} ; \mQ)$
show that they contain many non-trivial cycles, and each of these cycles is represented by a map $g \colon N \to B\overline{\Diff_c(M)}$ which is detected by  a secondary class in $H^*(B\oG_q; \mR)$ or in $H^*(B\G_q^+ ; \mR)$. Each such cycle induces  a representation  $g_{\#} \colon \pi_1(N, x_0) \to \Diff_c(M)$. It is unknown what is the ``geometric meaning'' of the non-triviality of the homology class of $g(N)$, or of its pairing with a secondary class used in its construction? The author's works \cite{Hurder1986,HK1987} related the vanishing of the secondary invariants with ergodic properties of the action of $\pi_1(N, x_0)$ on $M$ induced by  $g_{\#}$. However, the construction of the cycles $N$ using foliated surgery, as in Meigniez's proof of the Mather-Thurston Theorem in \cite{Meigniez2021}, suggests that the groups $\pi_1(N, x_0)$ are absolutely enormous, and in particular not amenable, so the question remains essentially open, as to what the constructions in this paper mean in a geometric or dynamical context.
}
\end{remark}

 We conclude with a discussion of another direction of research, which is very much in the spirit of the recent works by Galatius, Kupers, Nariman  and Randal-Williams. Let $M$ be compact, and 
 Consider  the basic fibration  sequence  
$B\overline{\Diff_c(M)} \stackrel{i}{\longrightarrow} B\Diff_c^+(M)^{\delta} \stackrel{\iota}{\longrightarrow} B\Diff_c^+(M)$ and its 
 associated Leray-Serre spectral sequence, with $E_2$-term 
 \begin{equation}\label{eq-ss4}
E_2^{r,s}(B\Diff_c^+(M)^{\delta} ; \mQ) \cong H^r(B\Diff_c^+(M) ; H^s(B\overline{\Diff_c(M)}) ; \mQ) \ . 
\end{equation}

The first step in analyzing \eqref{eq-ss4} is to understand the groups $H^*(B\overline{\Diff_c(M)}) ; \mQ)$, but this appears to be an almost completely unknown for dimensions $q \geq 4$. For example, according to Kupers  \cite{Kupers2019a}, apparently nothing is known about $H^*(B\overline{\Diff(\mS^4)}) ; \mQ)$. Even the simplest case when $M = \mD^q$  for $q \geq 4$ is wide open.   Kupers  \cite{Kupers2019a} writes:
\begin{quotation}
Question 1.1.3. What is the homotopy type of $\Diff_{c}(\mD^q)$? 
For example, what are its path components? What are its homology groups? What are its homotopy groups? Is it homotopy equivalent to a Lie group? Are its homotopy groups finitely generated? Can we relate it to other objects in algebraic topology? Can we relate it to other objects in manifold theory?
\end{quotation}

 One of the original motivations for the work in \cite{Hurder1985b} was a     result of Farrell and Hsiang  \cite{FarrellHsiang1978,Hsiang1979}:

\begin{thm}  \label{thm-FH} For $\ell \ll q$ and $\ell  \equiv 3$ mod 4 and $q$ odd:
\begin{equation}
\pi_{\ell}(\Diff_0(\mT^q)) \otimes \mQ = \mQ \ .
\end{equation}
\end{thm}
 \begin{prob}\label{prob-FH}
Show that the Farrell-Hsiang classes  in Theorem~\ref{thm-FH} vanish in the image of $\iota_{\delta}^* \colon H^k(\Diff_0(\mT^q) ; \mQ) \to H^k(\Diff_0(\mT^q)^{\delta} ; \mQ)$. 
 \end{prob}
The vanishing of classes in the image of $\iota_{\delta}^*$   implies the existence of non-trivial classes in some groups $H_{\ell}(B\overline{\Diff(\mT^q)}; \mQ)$ for $ 1< \ell \ll q$, and hence by the Mather-Thurston Theorem to the  existence of non-trivial classes in $H_k(B\oG_q ; \mQ)$ in the ``unknown range''  $q+1 < k \leq q$, so is of great interest.

The works of 
 Galatius and  Randal-Williams  greatly extend this seminal result  of Farrell and Hsiang in their works  \cite{GRW2014a,GRW2014b,GRW2018,GRW2017,GRW2020}, and the techniques used their have intriguing aspects analogous to ideas in foliation theory.  The (stabilized) homotopy groups of $\Diff_0(M)$ can be calculated  in terms of surgery groups and the Waldhausen K-theory of $M$, as described for example by Hatcher   in \cite{Hatcher1978,Hatcher2012} and Loday \cite{Loday1979}.
The standard method of calculating the Waldhausen K-theory groups introduces the notion of  \emph{block automorphisms},  denoted $\widetilde{\cA}(M)$, that were introduced by Burghelea in \cite{Burghelea1973}. A discussion of the relation between block automorphisms and the study of diffeomorphism groups is given in Section~2 of Hatcher \cite{Hatcher1978}. (See also the works of Burghelea \cite{BurgheleaLashoff1974,Burghelea1978,Burghelea1979}, Cohen \cite{Cohen1987} and Weiss and Williams \cite{WW2001}.) At an intuitive level, the notion of a block automorphism is analogous to the fracturing technique in the proof of the Mather-Thurston Theorem, as described say in Mather's paper \cite{Mather1979}. It seems a very intriguing question whether this analogy has a formal mathematical explanation. This suggests a far-reaching and speculative research program:
  
\begin{prob}\label{prob-Waldhausen}
Find geometric or algebraic correspondences between the Waldhausen K-theory groups of a manifold $M$ and the groups $\cK_*(M)$.
\end{prob}

 \appendix
 \section{Basic homotopy techniques}\label{sec-tool}

  We recall some basic results from homotopy theory that are frequently used in our constructions.
  \subsection{Hurewicz Theorems} \label{subsec-hurewicz}
 The Hurewicz Theorem relates  the homotopy and homology groups.
\begin{thm} \label{thm-hurewicz}
Let $X$ be an $r$-connected space for $r \geq 1$. That is, $X$ is a simply connected space with $\pi_k(X)  = 0$ for $0 \leq k \leq r$.   Then the Hurewicz homomorphism $h_* \colon \pi_{k}(X) \to H_{k}(X; \mZ)$ is an isomorphism for $0< k \leq r+1$.
\end{thm}

  The Rational   Hurewicz Theorem relates $\pi_*(X) \otimes \mQ$ with the groups $H_*(X ; \mQ)$.  
 Klaus and Kreck give an elementary proof of this result in  \cite{KlausKreck2004}, and provide   additional interesting discussions on this fundamental theorem.

\begin{thm}\label{thm-RHIT}
Let $X$ be a simply connected space with $\pi_k(X) \otimes \mQ = 0$ for $0 \leq k \leq r$. Then the Hurewicz map induces an isomorphism  
\begin{equation}
h \colon \pi_{k}(X) \otimes \mQ  \longrightarrow   H_{k}(X ; \mQ)
\end{equation}
for $1 \leq k \leq 2r$, and a surjection for $k=2r+1$.
\end{thm}

  \subsection{Integration over the fiber} \label{subsec-integration}
The operation of ``integration over the fiber''   is often used in the context of sphere bundles, but it applies much more generally.
It was first introduced in a paper by Chern   \cite{Chern1953}, and stated more formally by Borel and Hirzebruch in \cite[Section~8]{BorelHirzebruch1958}, where they explain the operation in terms of the differential on the $E_2$-term of the spectral sequence associated to the fibration. This definition was related to the definition using differential forms by Auer \cite{Auer1973}. Bott and Tu discuss   integration over the fiber for vector bundles in \cite[Section~6]{BottTu1982}.
\begin{thm} \label{thm-IoF}
Let $p \colon E \to B$ be a fibration with fiber $X = p^{-1}(b)$, where $B$  is connected and $X$ is a compact oriented manifold of dimension $q$. Then for $n \geq q$,  there exists a natural homomorphism 
\begin{equation}
\int_X \  \colon H^n(E ; \mR) \to H^{n-q}(B ; \mR)  \ ,  
\end{equation}
called  \emph{integration over the fiber}.
For $\alpha \in H^{n}(E ; \mR)$ and $\omega \in H^{\ell}(B ; \mR)$,  $\omega \cdot \int_X(\alpha) = \int_X(p^*\omega \cdot \alpha)$.
 \end{thm}

  \subsection{Brace products} \label{subsec-brace}

The \emph{brace product} was introduced by James in \cite{James1970}.  
 \begin{prop}\label{prop-brace} 
 Let  
$\ds 
F \stackrel{\iota}{\longrightarrow} E \stackrel{p}{\longrightarrow} B
$
 be a fibration, and suppose there exists a section $\sigma \colon B \to E$.  
 Then there is a well-defined product, called the  \emph{brace product},
 \begin{equation}\label{eq-brace}
 \{ \cdot , \cdot , \}_{\sigma} \colon \pi_i(E) \times \pi_j(F) \to \pi_{i+j-1}(F) \ .
 \end{equation}
 where the subscript $\sigma$ indicates that the product depends upon the choice   of the section $\sigma$. 
 \end{prop}
  \proof
 Given classes $\alpha \in \pi_i(E)$ and $\beta \in \pi_j(F)$ the Whitehead product 
 $[\sigma_{\#}(\alpha) , \iota_{\#}(\beta)] \in \pi_{i+j-1}(E_q)$ is the image of a unique class denoted $\{\alpha, \beta\}_{\sigma} \in \pi_{i+j-1}(B\overline{\Gamma}_q)$. 
 \endproof

The   brace product vanishes when the space $E$ is a product fibration, so the non-triviality of a pair $\{\alpha, \beta\}_{\sigma}$ is a measure of the ``homotopy twisting'' of the fibration.    Kallel and  Sjerve   \cite{KS2001} relate  the brace product to the structure of   free loop space $\Lambda E$.


\begin{thebibliography}{10}

 \bibitem{Auer1973}
{J.W.~Auer},
\newblock {\it Fiber integration in smooth bundles},
\newblock {\bf Pacific J. Math.}, 44:33--43, 1973.

 \bibitem{Baker1978}
{D.~Baker},
\newblock {\it On a class of foliations and the evaluation of their characteristic classes},
\newblock {\bf Comment. Math. Helv.}, 53:334--363, 1978.

 \bibitem{BamlerKleiner2023}
{R.~Bamler and B.~Kleiner},
\newblock {\it Ricci flow and diffeomorphism groups of 3-manifolds},
\newblock {\bf J. Amer. Math. Soc.}, 36:563--589, 2023.

 \bibitem{BarcusMeyer1958}
{W.D.~Barcus and J.-P.~Meyer},
\newblock {\it The suspension of a loop space},
\newblock {\bf Amer. J. Math.}, 80:895--920, 1958.
 
 \bibitem{BenderskyGitler1991}
{M.~Bendersky and S.~Gitler},
\newblock {\it The cohomology of certain function spaces},
\newblock {\bf Trans. Amer. Math. Soc.}, 326:423--440, 1991.

\bibitem{BernshteinRozenfeld1972}
{I.~Bernshtein and B.~Rozenfeld},
\newblock {\it On characteristic classes of foliations},
\newblock {\bf Funkts. Anal. Ego Prilozhen},  6:68-69, 1972.
\newblock Translation: {\bf Funct. Anal. Appl.},  6:60-61, 1972.

\bibitem{BernshteinRozenfeld1973}
{I.~Bernshtein and B.~Rozenfeld},
\newblock {\it Homogeneous spaces of infinite-dimensional Lie algebras and characteristic classes of foliations},
\newblock {\bf Uspehi Mat. Nauk.},  28:103-138, 1973.
\newblock Translation: {\bf Russian Math. Surveys},  28:107-142, 1973.

\bibitem{BorelHirzebruch1958}
{A.~Borel and F.~Hirzebruch},
\newblock {\it Characteristic classes and homogeneous spaces. {I}},
\newblock {\bf Amer. J. Math.}, 80:458--538,  1958.

\bibitem{Bott1970}
{R.~Bott},
\newblock {\it On a topological obstruction to integrability},
\newblock  In {\bf Global Analysis (Proc. Sympos. Pure Math., Vol. XVI, Berkeley, Calif., 1968)},
\newblock {Amer. Math. Soc., Providence, R.I.},  1970:127--131.

\bibitem{BottHaefliger1972}
{R.~Bott and A.~Haefliger},
\newblock {\it On characteristic classes of $\Gamma$-foliations},
\newblock {\bf  Bull. Amer. Math. Soc.}, 78:1039--1044, 1972.

\bibitem{Bott1972a}
{R.~Bott},
\newblock {\it Lectures on characteristic classes and foliations},
\newblock  In {\bf Lectures on algebraic and differential topology (Second Latin American School in Math., Mexico City, 1971)},
\newblock {Notes by Lawrence Conlon, with two appendices by J. Stasheff},
\newblock {Lecture Notes in Math., Vol. 279},  {Springer, Berlin}, 1972:1--94.

\bibitem{BottHeitsch1972}
{R.~Bott and J.~Heitsch},
\newblock {\it A remark on the integral cohomology of $B\Gamma_q$},
\newblock {\bf Topology}, 11:141-146, 1972.

\bibitem{Bott1975c}
{R.~Bott},
\newblock  {\it Gel{'}fand-Fuks Cohomology and Foliations},
\newblock in {\bf Proc. Eleventh Annual Holiday Symposium at New Mexico State University, December 27-31, 1973},
\newblock {Department of Math. Sciences Publication, Las Cruces, NM}, 1975:1--220.

\bibitem{BSS76}
R.~Bott, H.~Shulman and J. Stasheff, 
\textit{On the de Rham theory of certain classifying spaces}, 
Advances in Math. \textbf{20} (1976), 43--56.

\bibitem{Bott1978}
{R.~Bott},
\newblock {\it On some  formulas for the characteristic classes of group actions},
\newblock In {\bf Differential topology, foliations and Gelfand-Fuks cohomology (Proc. Sympos., Pontif\'\i cia Univ. Cat\'olica, Rio de Janeiro, 1976)},
\newblock {Lect. Notes in Math. Vol.  652},   {Springer--Verlag, Berlin}, 1978:25--61.

\bibitem{BottTu1982}
{R.~Bott and L.~Tu},
\newblock {\bf Differential forms in algebraic topology},
\newblock {Graduate Texts in Mathematics},  Vol. 82,  {Springer--Verlag, New York-Berlin}, 1982.

\bibitem{Boullay1996}
{P.~Boullay},
\newblock {\it {$H_3({\rm Diff}^+(S^2);{\mathbf Z})$} contains an uncountable {$\mathbf Q$}-vector space},
\newblock {\bf Topology}, 35:509--520, 1996.

 \bibitem{Bowden2012}
{J.~Bowden},
\newblock {\it The homology of surface diffeomorphism groups and a question of {M}orita},
\newblock {\bf Proc. Amer. Math. Soc.}, 140:2543--2549, 2012.


 \bibitem{Burghelea1973}
{D.~Burghelea},
\newblock {\it On the homotopy type of {${\rm diff}(M^{n})$} and connected problems},
\newblock {\bf Ann. Inst. Fourier (Grenoble)}, 23:3--17, 1973.

\bibitem{BurgheleaLashoff1974}
{D.~Burghelea and R.~Lashoff},
\newblock {\it The homotopy type of the space of diffeomorphisms. {I}, {II}},
\newblock {\bf Trans. Amer. Math. Soc.}, {1--36; ibid. 196 (1974), 37--50}.

\bibitem{Burghelea1978}
{D.~Burghelea},
\newblock {\it Automorphisms of manifolds},
\newblock in {\bf Algebraic and geometric topology ({P}roc. {S}ympos. {P}ure
              {M}ath., {S}tanford {U}niv., {S}tanford, {C}alif., 1976),
              {P}art 1}, 347--371, {Amer. Math. Soc., Providence, R.I.}, 1978.
  
\bibitem{Burghelea1979}
{D.~Burghelea},
\newblock {\it The rational homotopy groups of {D}iff {$(M)$} and {H}omeo{$(M^{n})$} in the stability range},
\newblock in {\bf Algebraic topology, {A}arhus 1978 ({P}roc. {S}ympos., {U}niv.
              {A}arhus, {A}arhus, 1978)}, 
              \newblock {Lecture Notes in Math.}, Vol. 763, {Springer, Berlin}, 1979:604--626.
 
\bibitem{Chern1953}
{S.-S.~Chern},
\newblock {\it On the characteristic classes of complex sphere bundles and algebraic varieties},
\newblock {\bf Amer. J. Math.}, 75:565--597,  1953.

\bibitem{CohenTaylor1978}
{F.R.~Cohen and L.R.~Taylor},
\newblock {\it Computations of {G}el'fand-{F}uks cohomology, the cohomology of function spaces, and the cohomology of configuration spaces},
\newblock in {\bf Geometric applications of homotopy theory ({P}roc. {C}onf.,
              {E}vanston, {I}ll., 1977), {I}}, 
              \newblock {Lecture Notes in Math.}, Vol. 657, {Springer, Berlin}, 1978:106--143.

\bibitem{CohenTaylor1983}
{F.R.~Cohen and L.R.~Taylor},
\newblock {\it The homology of function spaces},
\newblock in {\bf Proceedings of the {N}orthwestern {H}omotopy {T}heory
              {C}onference ({E}vanston, {I}ll., 1982)}, 
              \newblock {Contemp. Math.}, Vol. 19, {Amer. Math. Soc., Providence, RI}, 1983:39--50.

\bibitem{Cohen1987}
{R.~Cohen}, 
\newblock {\it Pseudo-isotopies, {$K$}-theory, and homotopy theory}, 
\newblock in {\bf Homotopy theory ({D}urham, 1985)}, 
\newblock {London Math. Soc. Lecture Note Ser.}, Vol. 15, {Cambridge Univ. Press, Cambridge}, 1987, 35--71.

\bibitem{DGMS1975}
{P.~Deligne, P.~Griffiths, J.~Morgan and D.~Sullivan}, 
\newblock {\it Real homotopy theory of {K}\"{a}hler manifolds}, 
\newblock {\bf Invent. Math.}, 29:245--274, 1975.
  
\bibitem{Dupont1978}
{J.L.~Dupont},
\newblock {\bf Curvature and characteristic classes},
\newblock  {Lecture Notes in Math.}, Vol. 640, {Springer, Berlin}, 1978.

\bibitem{EGM2011}
{Y.~Eliashberg, S.~Galatius and N.~Mishchev}, 
\newblock {\it Madsen-{W}eiss for geometrically minded topologists}, 
\newblock {\bf Geom. Topol.}, 15:411--472, 2011.
 
\bibitem{FarrellHsiang1978}
{F.T.~Farrell and W.-C.~Hsiang}, 
\newblock {\it On the rational homotopy groups of the diffeomorphism groups of discs, spheres and aspherical manifolds} 
\newblock in {\bf Algebraic and geometric topology ({P}roc. {S}ympos. {P}ure
              {M}ath., {S}tanford {U}niv., {S}tanford, {C}alif., 1976),
              {P}art 1}, 
\newblock {Proc. Sympos. Pure Math., XXXII},  Amer. Math. Soc., Providence, R.I. 1978:325--337.

\bibitem{FelixThomas2004}
{Y.~Felix and J.-C.~Thomas}, 
\newblock {\it Configuration spaces and {M}assey products}, 
\newblock {\bf J. Int. Math. Res. Not.}, 33:1685--1702, 2004.
 
\bibitem{Freedman2020}
{M.~Freedman}, 
\newblock {\it Controlled {M}ather-{T}hurston {T}heorems}, 
\newblock {\rm arxiv:2006.00374}.
 
\bibitem{Fuks1973}
{D.B.~Fuks},
\newblock {\it Characteristic classes of foliations},
\newblock {\bf Uspehi Mat. Nauk}, 28:3--17, 1973.
\newblock Translation: {\bf Russ. Math. Surveys }, 28:1--16, 1973.

\bibitem{Fuks1976}
{D.B.~Fuks},
\newblock {\it Finite-dimensional {L}ie algebras of formal vector fields and characteristic classes of homogeneous foliations},
\newblock {\bf Izv. Akad. Nauk SSSR Ser. Mat.}, 40:57--64, 1976.
\newblock Translation: {\bf Math. USSR Izv.}, 10:55--62, 1976.
    
\bibitem{GRW2014a}
{S.~Galatius and O.~Randal-Williams}, 
\newblock {\it Detecting and realising characteristic classes of manifold bundles}, 
\newblock in {\bf Algebraic topology: applications and new directions}, 
\newblock  {Contemp. Math.}, Vol. 620, {Amer. Math. Soc., Providence, RI}, 2014:99--110.
      
\bibitem{GRW2014b}
{S.~Galatius and O.~Randal-Williams}, 
\newblock {\it Stable moduli spaces of high-dimensional manifolds}, 
\newblock {\bf Acta Math.}, 212:257--377, 2014.
   
\bibitem{GRW2017}
{S.~Galatius and O.~Randal-Williams}, 
\newblock {\it Homological stability for moduli spaces of high dimensional manifolds. {II}}, 
\newblock {\bf Ann. of Math. (2)}, 186:127--204, 2017.
  
\bibitem{GRW2018}
{S.~Galatius and O.~Randal-Williams}, 
\newblock {\it Homological stability for moduli spaces of high dimensional manifolds. {I}}, 
\newblock {\bf J. Amer. Math. Soc.}, 31:215--264, 2018.
 
\bibitem{GRW2020}
{S.~Galatius and O.~Randal-Williams}, 
\newblock {\it Moduli spaces of manifolds: a user's guide}, 
\newblock in {\bf Handbook of homotopy theory}, 
\newblock  {CRC Press/Chapman Hall Handb. Math. Ser.},   2020:443--485.
      
\bibitem{GF1968}
{I.M.~Gel{'}fand, and D.B.~Fuks},
\newblock {\it Cohomologies of the {L}ie algebra of vector fields on the circle},
\newblock {\bf  Funkcional. Anal. i Prilo\v zen.},  2:92--93, 1968.

\bibitem{GF1969b}
{I.M.~Gel{'}fand, and D.B.~Fuks},
\newblock {\it Cohomologies of the {L}ie algebra of tangent vector fields of a smooth manifold},
\newblock {\bf  Funkcional. Anal. i Prilo\v zen.},  3:32--52, 1969.
\newblock  Translation: {\bf Funct. Anal. Appl.},  3:194-210,  1969.

\bibitem{Ghys1989}
{\'E.~Ghys},
\newblock {\it Sur l'invariant de {Godbillon-Vey}},
\newblock In {\bf S\'eminaire Bourbaki, Vol.\ 1988/89},
\newblock {Asterisque, Vol. 177-178},
\newblock {Soci\'et\'e Math\'ematique de France}, {Exp.\ No.\ 706, 155--181},  1989.

\bibitem{GodbillonVey1971}
C.~Godbillon and J.~Vey.
\newblock {\it Un invariant des feuilletages de codimension 1},
\newblock {\bf  C.R. Acad. Sci. Paris}, 273:92--95, 1971.
 
 \bibitem{Godbillon1974}
{C.~Godbillon},
\newblock {\it Cohomology of Lie algebras of formal vector fields},
\newblock In {\bf S\'eminaire Bourbaki, 25 \`eme ann\'ee   (1972/1973), Exp. No. 421},
\newblock {Lecture Notes in Math., Vol. 383},    {Springer, Berlin}, 1974:69--87.

\bibitem{Grayson2013}
{D.~Grayson},
\newblock {\it Quillen's work in algebraic {$K$}-theory},
\newblock {\bf J. K-Theory}, 11:527--547, 2013.

\bibitem{Greenberg1985}
{P.~Greenberg},
\newblock {\it Models for actions of certain groupoids},
\newblock {\bf Cahiers Topologie G\'{e}om. Diff\'{e}rentielle Cat\'{e}g.}, 26:33--42, 1985.

\bibitem{GM1981}
{P.A.~Griffiths and J.~Morgan}, 
\newblock {\bf Rational homotopy theory and differential forms}, 
\newblock {Progress in Mathematics, Vol. 16},  {Birkh\"{a}user, Boston, Mass.} 1981.
   
\bibitem{Gromov1969}
{M.L.~Gromov},
\newblock {\it  Stable mappings of foliations into manifolds}
\newblock {\bf  Izv. Akad. nauk SSSR Ser. Mat.}, 33:707--734, 1969.
\newblock Translation: {Math. USSR Izv.}, 33:671--694, 1969.
 		
\bibitem{Haefliger1970}
{A.~Haefliger},
\newblock {\it Feuilletages sur les vari\'et\'es ouvertes},
\newblock {\bf Topology}, 9:183--194, 1970.

\bibitem{Haefliger1971}
{A.~Haefliger},
\newblock {\it Homotopy and integrability},
\newblock in {\bf  Manifolds--Amsterdam 1970 (Proc. Nuffic Summer School)},
\newblock {Lect. Notes in Math.} Vol.  197,    Springer--Verlag,  Berlin, 1971:133--163.

\bibitem{Haefliger1978a}
{A.~Haefliger},
\newblock {\it On the {G}elfand-{F}uks cohomology},
\newblock {\bf Enseign. Math. (2)}, 24:143--160, 1978.




\bibitem{Haefliger1978b}
{A.~Haefliger},
\newblock {\it Whitehead products and differential forms},
\newblock In {\bf Differential topology, foliations and Gelfand-Fuks cohomology (Proc. Sympos., Pontif\'\i cia Univ. Cat\'olica, Rio de Janeiro, 1976)},
\newblock {Lect. Notes in Math. Vol.  652},   {Springer--Verlag, Berlin}, 1978:13--24.

\bibitem{Haefliger1982}
{A.~Haefliger},
\newblock {\it Rational homotopy of the space of sections of a nilpotent bundle},
\newblock   {\bf Trans. Amer. Math. Soc.}, 273:609--620, 1982.

\bibitem{Haller1998}
{S.~Haller},
\newblock {\it Perfectness and simplicity of certain groups of diffeomorphisms},
\newblock   {\bf Ph.D. Thesis}, University of Vienna, 2012. Available at {\text  https://www.mat.univie.ac.at/~stefan/files/diss.pdf}

\bibitem{Hamstrom1965}
{M.-E.~Hamstrom},
\newblock {\it The space of homeomorphisms on a torus},
\newblock {\bf Illinois J. Math.}, 9:59--65, 1965.

\bibitem{Hatcher1978}
{A.~Hatcher}, 
\newblock {\it Concordance spaces, higher simple-homotopy theory, and applications} 
\newblock in {\bf Algebraic and geometric topology ({P}roc. {S}ympos. {P}ure
              {M}ath., {S}tanford {U}niv., {S}tanford, {C}alif., 1976),
              {P}art 1}, 
\newblock {Proc. Sympos. Pure Math., XXXII},   {Amer. Math. Soc., Providence, R.I.}, 1976, 3--21.

\bibitem{Hatcher1981}
{A.~Hatcher},
\newblock {\it On the diffeomorphism group of {$S\sp{1}\times S\sp{2}$}},
\newblock {\bf Proc. Amer. Math. Soc.}, 83:427--430, 1981.

\bibitem{Hatcher1983}
{A.~Hatcher},
\newblock {\it A proof of the {S}male conjecture, {${\rm Diff}(S\sp{3})\simeq {\rm O}(4)$}},
\newblock {\bf Ann. of Math. (2)}, 117:553--607, 1983.

\bibitem{Hatcher2011}
{A.~Hatcher}, 
\newblock {\it A Short Exposition of the {M}adsen-{W}eiss {T}heorem} 
\newblock Available at {\rm http://pi.math.cornell.edu/$\sim$hatcher/Papers/MW.pdf}, also  {archiv:1103.5223.pdf}

\bibitem{Hatcher2012}
{A.~Hatcher}, 
\newblock {\it A 50-Year View of Diffeomorphism Groups} 
\newblock   {Slides from a talk at the 50th Cornell Topology Festival}, 
\newblock Available at {\rm http://pi.math.cornell.edu/$\sim$hatcher/Papers/Diff(M)2012.pdf}

\bibitem{Heitsch1977}
{J.~Heitsch},
\newblock {\it Residues and characteristic classes of foliations},
\newblock {\bf Bull. Amer. Math. Soc.}, 83:397--399, 1977.

\bibitem{Heitsch1978}
{J.~Heitsch},
\newblock {\it Independent variation of secondary classes},
\newblock {\bf Ann. of Math. (2)}, 108:421--460, 1978.

 \bibitem{HKMR2012}
{S.~Hong, J.~Kalliongis, D.~McCullough and J.H.~Rubinstein},
\newblock {\bf Diffeomorphisms of elliptic 3-manifolds},
\newblock {Lect. Notes in Math. Vol.  2055},   {Springer--Verlag, Berlin}, 2012.

\bibitem{Hsiang1979}
{W.C.~Hsiang},
\newblock {\it On {$\pi _{i}{\rm Diff}M^{n}$}},
\newblock in {\bf  Geometric topology ({P}roc. {G}eorgia {T}opology {C}onf.,
              {A}thens, {G}a., 1977)}, 
\newblock {Academic Press, New York-London}, 1979, pages 351--365.
  
          
\bibitem{Hurder1980}
{S.~Hurder},
\newblock {\it Dual homotopy invariants for {G}-foliations},
\newblock {\bf  Ph.D. Thesis}, University of Illinois at Urbana-Champaign, 1980.
     
\bibitem{Hurder1981a}
{S.~Hurder},
\newblock {\it Dual homotopy invariants for {G}-foliations},
\newblock {\bf  Topology}, 20:365--387, 1981.

\bibitem{Hurder1981b}
{S.~Hurder},
\newblock {\it On the secondary classes of foliations with trivial normal bundles},
\newblock {\bf  Comment. Math. Helv.}, 56:307--327, 1981.

\bibitem{Hurder1985a}
{S.~Hurder},
\newblock {\it On the classifying space of smooth foliations},
\newblock {\bf Illinois Jour. Math.}, 29:108-133, 1985.


\bibitem{Hurder1985b}
{S.~Hurder},
\newblock {\it Characteristic classes for flat {\rm Diff}(X) foliations},
\newblock {unpublished manuscript},  1985.

\bibitem{Hurder1986}
{S.~Hurder},
\newblock {\it The {Godbillon} measure of amenable foliations},
\newblock {\bf Jour. Differential Geom.}, 23:347--365, 1986.


\bibitem{Hurder1987}
{S.~Hurder},
\newblock {\it The homotopy theory of foliations},
\newblock {\bf Lectures Series},  Freie Universitat, W. Berlin,  June 15--14, 1987.

\bibitem{HK1987}
{S.~Hurder and A.~Katok},
\newblock {\it Ergodic  Theory and {Weil} measures for foliations},
\newblock {\bf Ann. of Math. (2)}, 126:221--275, 1987.

\bibitem{Hurder2002a}
{S.~Hurder},
\newblock {\it Dynamics and the {G}odbillon-{V}ey class: a History and Survey},
\newblock In {\bf Foliations: Geometry and Dynamics (Warsaw, 2000)},
\newblock {World Scientific Publishing Co. Inc., River Edge, N.J.},  2002:29--60.

        
\bibitem{Hurder2008}
{S.~Hurder},
\newblock {\it Classifying foliations},
\newblock in   {\bf Foliations, Topology and Geometry}, Contemp Math. Vol. 498, pages 1-65, Amer. Math. Soc., Providence, RI, 2009.

\bibitem{Hurder2023b}
{S.~Hurder},
\newblock {\it Divisible secondary characteristic classes},
\newblock {in preparation}, 2023.

\bibitem{James1970}
{I.M.~James},
\newblock {\it On the decomposability of fibre spaces},
\newblock in   {\bf The {S}teenrod {A}lgebra and its {A}pplications ({P}roc.
              {C}onf. to {C}elebrate {N}. {E}. {S}teenrod's {S}ixtieth
              {B}irthday, {B}attelle {M}emorial {I}nst., {C}olumbus, {O}hio,
              1970)}, Lecture Notes in Math. Vol. 168, 125--134, Springer, Berlin-New York, 1970.

 \bibitem{Jekel1984}
{S.~Jekel},
\newblock {\it Some weak equivalences for classifying spaces},
\newblock {\bf Proc. Amer. Math. Soc.}, 90:469--476, 1984.
  
 \bibitem{KS2001}
{S.~Kallel and D.~Sjerve},
\newblock {\it On the topology of fibrations with section and free loop spaces},
\newblock {\bf Proc. London Math. Soc. (3)}, 83:419--442, 2001.
  
 \bibitem{KT1974c}
{F.W.~Kamber and Ph.~Tondeur},
\newblock {\it Characteristic invariants of foliated bundles},
\newblock {\bf Manuscripta Math.}, 11:51--89, 1974.

\bibitem{KT1974d}
{F.W.~Kamber and Ph.~Tondeur},
\newblock {\it Non-trivial characteristic invariants of homogeneous foliated bundles},
\newblock {\bf Ann. Sci. \'Ecole Norm. Sup. (4)}, 8:433--486, 1975.

\bibitem{KT1975a}
{F.W.~Kamber and Ph.~Tondeur},
\newblock {\bf Foliated bundles and characteristic classes},
\newblock {Lect. Notes in Math. Vol. 493},   Springer-Verlag, Berlin, 1975.

\bibitem{KT1978a}
{F.W.~Kamber and Ph.~Tondeur},
\newblock {\it {$G$}-foliations and their characteristic classes},
\newblock {\bf Bull. Amer. Math. Soc.}, 84:1086--1124, 1978.

\bibitem{KT1979}
{F.W.~Kamber and Ph.~Tondeur},
\newblock {\it On the linear independence of certain cohomology classes of {$B\Gamma \sb{q}$}},
\newblock In {\bf Studies in algebraic topology},
\newblock {Adv. in Math. Suppl. Stud.}  5:213--263, 1979.

\bibitem{KlausKreck2004}
{S.~Klaus and M.~Kreck},
\newblock {\it A quick proof of the rational {H}urewicz theorem and a computation of the rational homotopy groups of spheres},
\newblock {\bf Math. Proc. Cambridge Philos. Soc.}, 136:617--623, 2004.
  
\bibitem{KotschickMorita2005}
{D.~Kotschick and S.~Morita},
\newblock {\it Signatures of foliated surface bundles and the symplectomorphism groups of surfaces},
\newblock {\bf Topology}, 44:131--144, 2005.
  
\bibitem{Kupers2019a}
{A.~Kupers},
\newblock {\bf Lectures on Diffeomorphism Groups of Manifolds},
\newblock  {preliminary version of text}, February 2019.
  

\bibitem{KRW2000}
{A.~Kupers and Randal-Williams},
\newblock {\it On diffeomorphisms of even-dimensional discs},
\newblock    {\rm arXiv:2007.13884}.

\bibitem{Lawson1977}
{H.B.~Lawson, Jr.},
\newblock {\bf The quantitative theory of foliations},
\newblock  {NSF Regional Conf. Board   Math. Sci.}, Vol. 27, {American Mathematical Society, Providence, R.I.}, 1977.

  \bibitem{Loday1979}
{J.-L.~Loday},
\newblock {\it Homotopie des espaces de concordances [d'apr\`es {F}.~{W}aldhausen]},
\newblock {S\'{e}minaire {B}ourbaki, 30e ann\'{e}e (1977/78)}, {Lect. Notes in Math. Vol. 710},   Springer, Berlin, 1979, 187--205.

\bibitem{MadsenWeiss2007}
{I.~Madsen and M.~Weiss},
\newblock {\it The stable moduli space of {R}iemann surfaces: {M}umford's conjecture},
\newblock {\bf Ann. of Math. (2)}, 165:843--941, 2007.

\bibitem{Margalit2019}
{D.~Margalit},
\newblock {\it Problems, questions, and conjectures about mapping class groups},
\newblock in {\bf Breadth in contemporary topology}, 
\newblock {Proc. Sympos. Pure Math.}, Vol. 102, {Amer. Math. Soc., Providence, RI}, 2019, pages 157--186.

 
  
\bibitem{Mather1971}
{J.N.~Mather},
\newblock {\it On {H}aefliger's classifying space. {I}},
\newblock {\bf Bull. Amer. Math. Soc.}, 77:1111--1115, 1971.

\bibitem{Mather1973}
{J.N.~Mather},
\newblock {\it Integrability in codimension one},
\newblock {\bf Comment. Math. Helv.}, 48:195--233, 1973.

\bibitem{Mather1975}
{J.N.~Mather},
\newblock {\it Loops and foliations},
\newblock In {\bf Manifolds---Tokyo 1973 (Proc. Internat. Conf., Tokyo, 1973)},
\newblock {Univ. Tokyo Press, Tokyo}, 1975, 175--180.

\bibitem{Mather1979}
{J.N.~Mather},
\newblock {\it On the homology of Haefliger's classifying space},
\newblock in  {\bf Differential topology (Varenna, 1976)},
\newblock  {Lectures presented at the Summer Session held by the Centro Internazionale Matematico Estivo (C.I.M.E.), Varenna, August 25--September 4, 1976},
\newblock {Liguori Editore, Naples}, 1979, 71--116.

\bibitem{McDuff1979a}
{D.~McDuff},
\newblock {\it Foliations and monoids of embeddings},
\newblock in {\bf Geometric topology ({P}roc. {G}eorgia {T}opology {C}onf., {A}thens, {G}a., 1977)},
\newblock {Academic Press, New York-London}, 1979, 429--444.

\bibitem{Meigniez2021}
{G.~Meigniez},
\newblock {\it Quasicomplementary foliations and the {M}ather-{T}hurston theorem},
\newblock {\bf Geom. Topol.}, 25:643--710, 2021.
        
\bibitem{MillerE1986}
{E.Y.~Miller},
\newblock {\it The homology of the mapping class group},
\newblock {\bf J. Differential Geom}, 24:1--14, 1986.

\bibitem{Milnor1983}
{J.~Milnor},
\newblock {\it On the homology of {L}ie groups made discrete},
\newblock {\bf Comment. Math. Helv.}, 58:72--85, 1983.

\bibitem{MitsumatsuVogt2017}
{Y.~Mitsumatsu and E.~Vogt},
\newblock {\it Thurston's h-principle for 2-dimensional foliations of codimension greater than one},
\newblock in {\bf Geometry, dynamics, and foliations 2013}, 
\newblock {Math. Soc. Japan, Tokyo}, 2017, pages 181--209.

\bibitem{Moerdijk1991}
{I.~Moerdijk},
\newblock {\it Classifying toposes and foliations},
\newblock {\bf Ann. Inst. Fourier (Grenoble)}, 41:189--209, 1991.

 \bibitem{MonodNariman2023}
{N.~Monod and S.~Nariman},
\newblock {\it Bounded and unbounded cohomology of homeomorphism and diffeomorphism groups},
\newblock   {\bf Invent. Math.},  232:1439--1475, 2023.  

\bibitem{Morita1977}
{S.~Morita},
\newblock {\it A remark on the continuous variation of secondary characteristic classes for foliations},
\newblock {\bf J. Math. Soc. Japan}, 29:253--260, 1977. 
 
\bibitem{Morita1984a}
{S.~Morita},
\newblock {\it Nontriviality of the {G}elfand-{F}uchs characteristic classes for flat {$S^{1}$}-bundles},
\newblock   {\bf Osaka J. Math.},  21:545--563, 1984.  

\bibitem{Morita1984b}
{S.~Morita},
\newblock {\it Characteristic classes of surface bundles},
\newblock   {\bf Bull. Amer. Math. Soc. (N.S.)},  11:386--388, 1984.  

\bibitem{Morita1985}
{S.~Morita},
\newblock {\it Discontinuous invariants of foliations},
\newblock in {\bf Foliations (Tokyo, 1983)},
\newblock {Advanced Studies in Pure Math.  Vol. 5},
\newblock {North-Holland, Amsterdam}, 1985:169--193.

\bibitem{Morita1987}
{S.~Morita},
\newblock {\it Characteristic classes of surface bundles},
\newblock   {\bf Invent. Math.},  90:551--577, 1987.  

\bibitem{Morita2001}
{S.~Morita},
\newblock {\bf Geometry of characteristic classes},
\newblock {American Mathematical Society, Providence, R.I.}, 2001.

\bibitem{Morita2006}
{S.~Morita},
\newblock {\it Cohomological structure of the mapping class group and beyond},
\newblock in {\bf  Problems on mapping class groups and related topics}, 
{Proc. Sympos. Pure Math.}, Vol. 74, {American Mathematical Society, Providence, R.I.}, 2006, pages 329--354.

\bibitem{Moskowitz1985}
{I.~Moskowitz},
\newblock {\it A note on the {B}ott vanishing theorem},
\newblock   {\bf Proc. Amer. Math. Soc.}, V94:529--530, 1985.

\bibitem{MostowM1976}
{M.~Mostow},
\newblock {\it Continuous cohomology of spaces with two topologies},
\newblock   {\bf Mem. Amer. Math. Soc.}, Vol. 7, 1976.

\bibitem{Nariman2017a}
{S.~Nariman},
\newblock {\it Homological stability and stable moduli of flat manifold bundles},
\newblock   {\bf Adv. Math.},  320:1227--1268, 2017.  
 
\bibitem{Nariman2017b}
{S.~Nariman},
\newblock {\it Stable homology of surface diffeomorphism groups made discrete},
\newblock   {\bf Geom. Topol.},  5:3047--3092, 2017.  
 
\bibitem{Nariman2020}
{S.~Nariman},
\newblock {\it A local to global argument on low dimensional manifolds},
\newblock   {\bf Trans. Amer. Math. Soc.},  373:1307--1342, 2020.  
 	
 \bibitem{Nariman2022}
{S.~Nariman},
\newblock {\it On flat manifold bundles and the connectivity of Haefliger's classifying spaces},
\newblock   {\rm arXiv:2202.00052}.  
  	
\bibitem{Nariman2023b}
{S.~Nariman},
\newblock {\it On invariants of foliated sphere bundles},
\newblock   {\rm arXiv:2308.16310}.
 
 \bibitem{Pelletier1983}
{W.~Pelletier},
\newblock {\it The secondary characteristic classes of solvable foliations},
\newblock {\bf  Proc. Amer. Math. Soc.}, 88:651--659, 1983.

 \bibitem{Phillips1968}
{A.~Phillips},
\newblock {\it Foliations on open manifolds. {I}},
\newblock {\bf  Comment. Math. Helv.}, 43:204--211, 1968.

 \bibitem{Phillips1969}
{A.~Phillips},
\newblock {\it Foliations on open manifolds. {II}},
\newblock {\bf  Comment. Math. Helv.}, 44:367--370, 1969.

\bibitem{Pittie1976a}
{H.V.~Pittie},
\newblock {\bf Characteristic classes of foliations},
\newblock  {Research Notes in Mathematics, No. 10},
\newblock {Pitman Publishing, London-San Francisco, Calif.-Melbourne}, 1976.

\bibitem{Pittie1979}
{H.V.~Pittie},
\newblock {\it The secondary characteristic classes of parabolic foliations},
\newblock {\bf Comment. Math. Helv.}, 54:601--614, 1979.

\bibitem{Rasmussen1980}
{O.H.~Rasmussen},
\newblock {\it Continuous variation of foliations in codimension two},
\newblock {\bf Topology}, 19:335--349, 1980.

\bibitem{Samelson1953}
{H.~Samelson},
\newblock {\it A connection between the {W}hitehead and the {P}ontryagin product},
\newblock {\bf Amer. J. Math.}, 75:744--752, 1953.

\bibitem{SchweitzerWhitman1978}
{P.A.~Schweitzer and A.~Whitman},
\newblock {\it Pontryagin polynomial residues of isolated foliation singularities},
\newblock In {\bf Differential topology, foliations and Gelfand-Fuks cohomology (Proc. Sympos., Pontif\'\i cia Univ. Cat\'olica, Rio de Janeiro, 1976)},
\newblock {Lect. Notes in Math. Vol.  652},   {Springer--Verlag, Berlin}, 1978:95--103.

\bibitem{Segal1978}
{G.~Segal}, 
\newblock {\it Classifying spaces related to foliations}, 
\newblock {\bf Topology}, 17:367--382, 1978.

\bibitem{Shibata1984a}
{K.~Shibata},
\newblock {\it Sullivan-{Q}uillen mixed type model for fibrations and the {H}aefliger model for the {G}el'fand-{F}uks cohomology}, 
\newblock in {\bf Algebraic homotopy and local algebra ({L}uminy, 1982)},  
\newblock {Ast\'{e}risque}, Vol.  113,   {Soc. Math. France, Paris}, 1984:292--297.

\bibitem{Shibata1984b}
{K.~Shibata},
\newblock {\it Sullivan-{Q}uillen mixed type model for fibrations}, 
\newblock {\bf J. Math. Soc. Japan},   36:221--242, 1984.

\bibitem{Stasheff1978}
{J.D.~Stasheff}, 
\newblock {\it Continuous cohomology of groups and classifying spaces}, 
\newblock {\bf Bull. Amer. Math. Soc.}, 84:513--530, 1978.

\bibitem{Sullivan1976}
{D.~Sullivan},
\newblock {\it Cartan-de {R}ham homotopy theory}, 
\newblock in {\bf Colloque ``{A}nalyse et {T}opologie'' en l'{H}onneur de {H}enri {C}artan ({O}rsay, 1974)} 
\newblock {Ast\'{e}risque}, No. 32-33, 1976, pages 227--254.

\bibitem{SVP1976}
{D.~Sullivan and M.~Vigu\'{e}-Poirrier},
\newblock {\it The homology theory of the closed geodesic problem},
\newblock {\bf J. Differential Geometry}, 11:633--644, 1976.

\bibitem{Sullivan1978}
{D.~Sullivan},
\newblock {\it Infinitesimal computations in topology}, 
\newblock {\bf Inst. Hautes \'{E}tudes Sci. Publ. Math.} 47:269--331, 1978.
   
\bibitem{Tabachnikov1984}
{S.~Tabachnikov},
\newblock {\it  Characteristic classes of homogeneous foliations},
\newblock {\bf Uspekhi Mat. Nauk}, 39:189--190, 1984.
\newblock Translation: {\bf Russian Math. Surveys}, 39:203-204, 1984.

\bibitem{Tabachnikov1985}
{S.~Tabachnikov},
\newblock {\it  Characteristic classes of parabolic foliations and symmetric functions},
\newblock {\bf Serdica},
\newblock   {Serdica. Bulgaricae Mathematicae Publicationes}, 11:86--95, 1985.

\bibitem{Tabachnikov1986}
{S.~Tabachnikov},
\newblock {\it  Characteristic classes of parabolic foliations of series {$B$}, {$C$}, and {$D$} and degrees of isotropic {G}rassmannians},
\newblock {\bf Funktsional. Anal. i Prilozhen.}, 20:84--85, 1986.
\newblock Translation: {\bf Functional Anal. Appl.}, 20:158-160, 1986.

\bibitem{Thurston1972}
{W.P.~Thurston},
\newblock {\it Noncobordant foliations of {$S^{3}$}}
\newblock {\bf Bull. Amer. Math. Soc.} 78:511--514, 1972.

\bibitem{Thurston1974}
{W.P.~Thurston},
\newblock {\it Foliations and groups of diffeomorphisms}
\newblock {\bf Bull. Amer. Math. Soc.} 80:304--307, 1974.

\bibitem{Thurston1974b}
{W.P.~Thurston},
\newblock {\it The theory of foliations of codimension greater than one}
\newblock {\bf Comment. Math. Helv.} 49:214--231, 1974.

\bibitem{Tillmann2012}
{U.~Tillmann},
\newblock {\it Spaces of graphs and surfaces: on the work of {S}\o ren {G}alatius}
\newblock {\bf Bull. Amer. Math. Soc. (N.S.)} 49:73--90, 2012.

\bibitem{Tsuboi1984}
{T.~Tsuboi},
\newblock {\it Foliated cobordism classes  of certain foliated $\mS^1$-bundles over surfaces},
\newblock {\bf Topology}, 23:233--244, 1984.
 

\bibitem{Tsuboi1998}
{T.~Tsuboi},
\newblock {\it Irrational foliations of {$S^3\times S^3$}},
\newblock {\bf Hokkaido Math. J.}, 27:605--623, 1998.


\bibitem{Tsujishita1977}
{T.~Tsujishita},
\newblock {\it On the continuous cohomology of the {L}ie algebra of vector fields},
\newblock {\bf Proc. Japan Acad. Ser. A Math. Sci.}, 53:134--138, 1977.

\bibitem{Tsujishita1981}
{T.~Tsujishita},
\newblock {\it Continuous cohomology of the {L}ie algebra of vector fields},
\newblock {\bf Mem. Amer. Math. Soc.}, Vol. 34, 1981.

\bibitem{Weibel2013}
{C.~Weibel},
\newblock {\bf The {$K$}-book},
\newblock    {Graduate Studies in Mathematics, Vol. 145}, {American Mathematical Society, Providence, RI},  2013.
  
\bibitem{WW2001}
{M.~Weiss and B.~Williams},
\newblock {\it Automorphisms of manifolds},
\newblock {\bf  Surveys on surgery theory, Vol. 2},  {Ann. of Math. Stud., vol. 149}, {Princeton Univ. Press, Princeton, NJ}, 2001, pages 165?220.  {\rm arXiv:0012101}.

\bibitem{Yamato1975}
{K.~Yamato},
\newblock {\it Examples of foliations with non trivial exotic characteristic classes}
\newblock {\bf Osaka J. Math.} 12:401--417, 1975.


\end{thebibliography}
\end{document}